\documentclass[12pt]{article}
\usepackage[usenames,dvipsnames,svgnames,table]{xcolor} 
\usepackage{kotex}
\usepackage{graphicx}
\usepackage{epsfig}
\usepackage{amssymb, amsmath, mathtools}
\usepackage{verbatim}
\usepackage{natbib}
\usepackage[affil-it]{authblk}
\usepackage{mathrsfs}
\usepackage{here}
\usepackage{makecell}
\usepackage{tikz}
\usepackage{enumerate}
\usepackage[shortlabels]{enumitem}
\usepackage{yfonts}
\usepackage{comment}
\usepackage{longtable}
\usepackage{tikz,amsmath, amssymb,bm,color}
\usetikzlibrary{circuits.logic.US,circuits.logic.IEC,fit}
\usepackage{algorithm2e}
\usepackage{multirow}
\usepackage{physics}
\usepackage{amsthm}
\usepackage{makecell}
\usepackage{ulem}
\usepackage{proofread}
\usepackage{booktabs}
\usepackage{array}
\usepackage{xr}
\usepackage{adjustbox}
\usepackage{lineno}
\usepackage{import}
\usepackage{subcaption}

\usepackage[colorlinks]{hyperref} 
\usepackage[colorinlistoftodos, textsize=scriptsize]{todonotes} 

\usepackage[framemethod=TikZ]{mdframed}
\mdfdefinestyle{JL}{%
	linecolor=blue,
	outerlinewidth=1pt,
	roundcorner=2pt,
	innertopmargin=0.1\baselineskip,
	innerbottommargin=0.1\baselineskip,
	innerrightmargin=2pt,
	innerleftmargin=2pt,
	backgroundcolor=orange!100!white}

\hypersetup{
	colorlinks,
	citecolor=Black,
	linkcolor=Red,
	urlcolor=Blue}

\includecomment{note} 

\newtheorem{theorem}{Theorem}[section]
\newtheorem{lemma}[theorem]{Lemma}
\newtheorem{definition}[theorem]{Definition}
\newtheorem{proposition}[theorem]{Proposition}
\newtheorem{corollary}[theorem]{Corollary}

\newtheorem{remark}[theorem]{Remark}

\DeclareMathOperator{\argmin}{argmin}

\everymath{\displaystyle}

\newcommand{\sse}[1]{\subsection{#1}}

\newcommand{\be}{\begin{equation}}
	\newcommand{\ee}{\end{equation}}
\newcommand{\bea}{\begin{eqnarray*}}
	\newcommand{\eea}{\end{eqnarray*}}
\newcommand{\bean}{\begin{eqnarray}}
	\newcommand{\eean}{\end{eqnarray}}
\newcommand{\ben}{\begin{enumerate}}
	\newcommand{\een}{\end{enumerate}}
\newcommand{\bi}{\begin{itemize}}
	\newcommand{\ei}{\end{itemize}}
\newcommand{\brem}{\begin{remark}}
	\newcommand{\erem}{\end{remark}}
\newcommand{\bcen}{\begin{center}}
	\newcommand{\ecen}{\end{center}}
\newcommand{\bsv}{\begin{semiverbatim}}
	\newcommand{\esv}{\end{semiverbatim}}

\newcommand{\bt}{\begin{theorem}}
	\newcommand{\et}{\end{theorem}}
\newcommand{\bl}{\begin{lemma}}
	\newcommand{\el}{\end{lemma}}
\newcommand{\bd}{\begin{definition}}
	\newcommand{\ed}{\end{definition}}
\newcommand{\bc}{\begin{corollary}}
	\newcommand{\ec}{\end{corollary}}
\newcommand{\bp}{\begin{proposition}}
	\newcommand{\ep}{\end{proposition}}

\newcommand{\bfX}{ \mathbf{X}}

\newcommand{\bbE}{{ \mathbb{E}}}

\newcommand{\bbR}{ \mathbb{R}}

\newcommand{\rmU}{\mathrm{U}}

\usepackage{xr}
\makeatletter
\newcommand*{\addFileDependency}[1]{
	\typeout{(#1)}
	\@addtofilelist{#1}
	\IfFileExists{#1}{}{\typeout{No file #1.}}
}
\makeatother

\modulolinenumbers[5]

\title{ Eigenstructure inference for high-dimensional covariance with generalized shrinkage inverse-Wishart prior }
\author[1]{Seongmin Kim}
\author[2]{Kwangmin Lee}
\author[3]{Sewon Park}
\author[1]{Jaeyong Lee}
\affil[1]{Department of Statistics, Seoul National University}
\affil[2]{Department of Big Data Convergence, Chonnam National University}
\affil[3]{Department of Statistics, Sookmyung Women's University}

\allowdisplaybreaks
\begin{document}
	
\maketitle

\begin{abstract}

In multivariate statistics, estimating the covariance matrix is essential for understanding the interdependence among variables. In high-dimensional settings, where the number of covariates increases with the sample size, it is well known that the eigenstructure of the sample covariance matrix is inconsistent. The inverse-Wishart prior, a standard choice for covariance estimation in Bayesian inference, also suffers from posterior inconsistency.

To address the issue of eigenvalue dispersion in high-dimensional settings, the shrinkage inverse-Wishart (SIW) prior has recently been proposed. Despite its conceptual appeal and empirical success, 
the asymptotic justification for the SIW prior has remained limited. 

In this paper, we propose a generalized shrinkage inverse-Wishart (gSIW) prior for high-dimensional covariance modeling. By extending the SIW framework, the gSIW prior accommodates a broader class of prior distributions and facilitates the derivation of theoretical properties under specific parameter choices.

In particular, under the spiked covariance assumption, we establish the asymptotic behavior of the posterior distribution for both eigenvalues and eigenvectors by directly evaluating the posterior expectations for two sets of parameter choices. This direct evaluation provides insights into the large-sample behavior of the posterior that cannot be obtained through general posterior asymptotic theorems.

Finally, simulation studies illustrate that the proposed prior provides accurate estimation of the eigenstructure, particularly for spiked eigenvalues, achieving narrower credible intervals and higher coverage probabilities compared to existing methods. For spiked eigenvectors, the performance is generally comparable to that of competing approaches, including the sample covariance.

\end{abstract}


\section{Introduction} 
Suppose $X_1,\ldots,X_n$ are $p$-dimensional independent samples from $N(0,\Sigma)$, where $\Sigma$ is a $p\times p$ covariance matrix. In the study of multivariate statistics, estimating the covariance matrix is crucial for understanding the interdependence among variables. Extensive research has been conducted on covariance estimation, including \cite{efron1976multivariate}, \cite{dey1985estimation}, and \cite{daniels2001shrinkage}.

In the high-dimensional settings, where the number of covariates increases with the sample size, \cite{johnstone2001distribution} and \cite{johnstone2009consistency} demonstrated the inconsistency of the largest eigenvalue and eigenvector of sample covariance matrix.

To address the challenges in high-dimensional settings, several shrinkage estimators have been proposed. \cite{stein1975estimation}, \cite{haff1980empirical}, \cite{ledoit2004well} and \cite{bodnar2014strong} suggested shrinkage estimators for large covariances by minimizing specific loss functions. Additionally, \cite{karoui2008spectrum}, \cite{ledoit2012nonlinear}, and \cite{lam2016nonparametric} proposed orthogonal invariant shrinkage estimators based on random matrix theory.  
These methods require the ratio $p/n$ to be bounded, which limits their applicability in ultra-high dimensional settings. In contrast, we consider the regime where $p/n \rightarrow \infty$.

Some structural assumptions are necessary to address the problem of estimating high-dimensional covariance where $p/n \rightarrow \infty$. The first structural assumption is the sparsity, either in the covariance or precision matrix. In the Frequentist literature, \cite{bickel2008covariance}, \cite{cai2012optimal}, and \cite{cai2016estimating} suggested sparse covariance estimation using thresholding methods. In the Bayesian literature, \cite{banerjee2015bayesian} suggested sparse precision estimation using the Gaussian graphical model. \cite{lee2022beta} introduced a beta-mixture shrinkage prior for sparse covariance and showed the posterior has nearly optimal contraction rate. Recently, \cite{lee2023post_sparse} suggested post-processed posterior for sparse covariance model. Another structure assumption is bandedness. In the Frequentist literature, \cite{bickel2008regularized}, \cite{cai2010optimal} and \cite{bien2016convex} proposed banded covariance using tapering or banding method. In the Bayesian literature, \cite{banerjee2014posterior} introduced banded precision estimation by using the Gaussian graphical model. Recently, \cite{lee2023post_banded} suggested post-processed posterior for the banded covariance model. These assumptions of sparsity or bandedness require prior knowledge of the covariance matrix structure, which may not always be available.

Without requiring sparsity or bandedness assumptions on the covariance,  the spiked covariance assumption was initially proposed by \cite{johnstone2001distribution}. The spiked covariance $\Sigma$ is composed of $k$ large eigenvalues, while the remaining eigenvalues are relatively small. \cite{paul2007asymptotics}, \cite{shen2016statistics}, and \cite{wang2017asymptotics} investigated the asymptotic properties of the eigenstructure of sample covariances under the assumption that the true covariance is spiked covariance. Further, \cite{wang2017asymptotics} suggested an estimator of covariance based on the asymptotic properties of sample covariance. In the Bayesian study, \cite{ning2021spike} and \cite{xie2022bayesian} proposed a prior for the spiked covariance and verified the posterior convergence rate of covariance. 

In Bayesian inference, the inverse-Wishart prior is commonly used for covariance estimation. Surprisingly, there has been little prior research in Bayesian statistics on the eigenstructures of unconstrained covariance matrices. Recently, \cite{lee2024posterior} explored the asymptotic properties of the eigenstructure of covariance matrices under the inverse-Wishart prior, marking the first time such properties have been verified in a Bayesian context. Additionally, the shrinkage inverse-Wishart (SIW) prior was proposed by \cite{berger2020bayesian} to address some of the limitations of the inverse-Wishart prior. However,  despite its conceptual appeal and empirical success, 
the asymptotic justification for the SIW prior has remained limited. 

In this paper, we propose a generalized shrinkage inverse-Wishart (gSIW) prior that accommodates a broader class of prior distributions and facilitates theoretical analysis of the eigenstructure of high-dimensional covariance matrices under specific parameter choices. It is known that the eigenvalues of the sample covariance matrix are biased (\citealt{wang2017asymptotics}), and that the inverse-Wishart posterior is also biased for the eigenvalues (\citealt{lee2024posterior}). When parameters are chosen appropriately, the proposed gSIW prior offers an asymptotically unbiased estimator for the eigenstructure and enables rapid posterior computation.

There has been extensive research on posterior convergence rates for the covariance matrix in high-dimensional settings, including \cite{banerjee2014posterior}, \cite{xiang2015high}, \cite{lee2022beta}, and \cite{lee2023post_sparse}. These studies primarily focus on covariance or the precision matrix itself, with little attention given to the eigenstructures of unconstrained covariance matrices in the Bayesian framework.


Commonly, to investigate the asymptotic properties of the posterior, researchers resort to general theorems on posterior convergence rates (e.g., \cite{ghosal2000convergence}; \cite{ghosal2007convergence}). However, in this paper, we directly evaluate the posterior expectations of eigenvalues and eigenvectors.
Our approach offers a couple of advantages. In particular, we derive an explicit expression for the bounds of the posterior expectations, which provides deeper statistical insights into the behavior of the posterior. Moreover, our method is applicable to general functionals of the covariance matrix.

The rest of the paper is structured as follows. In Section \ref{sec:bayes_inf}, we propose a generalization of the shrinkage inverse-Wishart introduced by \cite{berger2020bayesian}. We also suggest a sampling method for the posterior. Section \ref{sec:main_results} presents the posterior convergence rates for the eigenvalues and eigenvectors of the covariance matrix. Simulation studies are discussed in Section \ref{sec:simul_studies}, followed by real data analysis in Section \ref{sec:real_data}. A conclusion is provided in Section \ref{sec:conclusion}. The proofs of the theorems and related lemmas are provided in the Appendix.


\section{Generalized shrinkage inverse-Wishart prior and its posterior}\label{sec:bayes_inf}

\sse{Notation}
For any positive sequence $a_n$ and $b_n$, we write $a_n=o (b_n)$ or equivalently $a_n\prec b_n$, if $\dfrac{a_n}{b_n}\rightarrow 0$ as $n \to \infty$.  We write $a_n=O (b_n)$ or equivalently $a_n\preccurlyeq b_n$ if  $\abs{\dfrac{a_n}{b_n}}\leq c$ for some constant $c > 0$. For real constants $a$ and $b$, we use $a\wedge b$ and $a\vee b$ to denote the minimum and maximum values between $a$ and $b$, respectively. Let $\abs{\abs{A}}_F=\sqrt{tr(A^TA)}$ denote the Frobenius norm, and $\abs{\abs{A}}$ denote the spectral norm, i.e., the largest singular value of $A$. For a squared matrix $A$, we use  $etr(A)$  to represent  $exp(tr(A))$. We define the set of $p$ dimensional orthogonal matrices as $O(p) :=\{\rmU\in \bbR^{p\times p}: \rmU^T\rmU =I_p\}$ and the set of  $p\times p$ positive definite matrices as $\mathcal{C}_p :=\{ A \in \bbR^{p\times p}: $A$ \text{ is positive definite } \} $. 

We write $X \sim \text{InvGam}(\alpha, \beta)$ if $X$ follows the inverse-gamma distribution with shape parameter $\alpha>0$ and scale parameter $\beta>0$, having density proportional to $x^{-\alpha - 1} \exp(-\beta / x)$ for $x > 0$. The notation $X_i \overset{\text{iid}}{\sim} P$ indicates that the random variables $X_i$ are independent and identically distributed according to distribution $P$, while $X_i \overset{\text{ind}}{\sim} P_i$ denotes that the $X_i$'s are independent but not necessarily identically distributed.

Suppose $\Sigma = (\sigma_{ij}, 1 \leq i, j \leq p)$ is a $p \times p$ symmetric matrix, $\Lambda = \text{diag}(\lambda_1, \ldots, \lambda_p)$ is a $p \times p$ diagonal matrix, and $\rmU = [u_1, \ldots, u_p] \in O(p)$, where $u_i$ is the $i$-th column of $\rmU$. We define the matrix differential forms of these matrices as follows. Let $(d\Sigma) := \bigwedge_{i \leq j} d\sigma_{ij}$ denote the exterior product of the distinct elements of $\Sigma$, $(d\Lambda) := \bigwedge_{i=1}^p d\lambda_i$, and $(d\rmU) := \bigwedge_{1 \leq i \leq j \leq p} u_j^T du_i$, which represents the differential form corresponding to the invariant measure on $O(p)$. Note that the volume of a set in $O(p)$ can be obtained by integrating $(d\rmU)$ over the set, i.e.,  
\[
\text{volume}(D) = \int_D (d\rmU), \quad D \subset O(p).
\]  
For a detailed discussion on differential forms related to matrices, see \cite{muirhead2009aspects}.


\subsection{Generalized shrinkage inverse-Wishart prior}
Suppose $X_1,\ldots,X_n$ are independent random vectors following a $p$-dimensional multivariate normal distribution with covariance matrix $\Sigma \in \mathcal{C}_p$:
\begin{equation}\label{eq:model}
    X_i\overset{iid}{\sim}N(0,\Sigma),\quad\forall i=1,\ldots,n.
\end{equation}

\cite{berger2020bayesian} proposed the shrinkage inverse-Wishart (SIW) prior, a new class of priors for $\Sigma$ whose densities are given by:
$$\pi^{SIW}(\Sigma) (d\Sigma) = \pi^{SIW}(\Sigma\vert a,b,H) (d\Sigma) \propto \dfrac{etr(-\dfrac{1}{2}\Sigma^{-1}H)}{\abs{\Sigma}^a[{\prod\limits_{i<j}(\lambda_i-\lambda_j)]}^b} (d\Sigma),$$ 
where $a \geq 0$, $b  \in [0,1]$, $H \in \mathcal{C}_p$ and $\lambda_1,\ldots,\lambda_p$ are the ordered eigenvalues of $\Sigma$. This prior is called the shrinkage inverse-Wishart prior, as it induces shrinkage of the eigenvalues by the term $\prod\limits_{i<j}(\lambda_i-\lambda_j)$ compared to the inverse-Wishart prior. When $b = 0$, the SIW prior reduces to the inverse-Wishart prior.

When $b<0$, the eigenvalue spreading problem, which is already problematic under the inverse-Wishart prior in high-dimensional settings, becomes even more severe. On the other hand, when $b>1$, the prior forces the eigenvalues to become nearly identical. Accordingly, as proposed in \cite{berger2020bayesian}, we restrict $b$ to the range $[0,1]$ to avoid these undesirable behaviors.

Consider the following spectral decomposition of $\Sigma= \rmU\Lambda \rmU^T$, where $\Lambda=\text{diag}(\lambda_1,\ldots,\lambda_p)$ and $\rmU =[u_1, \ldots, u_p]$ is a $p\times p$ matrix  whose $i$-th column $u_i$ is the eigenvector corresponding to the $i$-th eigenvalue $\lambda_i$. We can rewrite $(d\Sigma)$ in terms of $(d\rmU)$ and $(d\Lambda)$: 
\begin{equation} \label{eq:jacobian}
(d\Sigma) =\prod\limits_{i<j}(\lambda_i-\lambda_j)(d\Lambda) (d\rmU).
\end{equation}
See \cite{farrell2012multivariate} and \cite{muirhead2009aspects}.

Using equation \eqref{eq:jacobian}, the density of the shrinkage inverse-Wishart prior can be written as 
$$\pi^{SIW}(\Lambda,\rmU|a,b,H) (d\Lambda) (d\rmU) \propto \dfrac{etr(-\dfrac{1}{2}\Sigma^{-1}H)}{\abs{\Sigma}^a[\prod\limits_{i<j}(\lambda_i-\lambda_j)]^{b-1}}(d\Lambda) (d\rmU),$$
where $\lambda_1,\ldots,\lambda_p$ are the ordered eigenvalues of $\Sigma$. Since the SIW prior is symmetric with respect to the eigenvalues, the ordering constraint can be removed without affecting the distribution. The SIW prior has demonstrated  empirically successful performance but is not supported by comprehensive theoretical analysis. In this paper, we propose a generalized Shrinkage inverse-Wishart (gSIW) prior whose density is defined as 
\begin{align}\label{eq:SIW}
    \pi^{gSIW}(\Lambda,\rmU|a_1,\ldots,a_p,b,H) (d\Lambda) (d\rmU) \propto   
    \dfrac{etr(-\dfrac{1}{2}\rmU\Lambda^{-1} \rmU^T H)}{\prod\limits_{i=1}^p\lambda_i^{a_i}\prod\limits_{i<j}\abs{\lambda_i-\lambda_j}^{b-1}}(d\Lambda) (d\rmU),
\end{align} 
where $a_1,\ldots,a_p>0$, $b\in[0,1]$, and $\lambda_1,\cdots,\lambda_p$ are unordered eigenvalues of $\Sigma$. By allowing each $a_i$ to be defined differently for each $\lambda_i$, we can consider a prior distribution that is more general than the SIW prior. In particular, when $a_1 = \cdots = a_p = a$, the gSIW prior reduces to the SIW prior. When $b=1$, the term $\abs{\lambda_i-\lambda_j}^{1-b}$ disappears from the prior distribution of $(\Lambda,\rmU)$, which facilitates faster sampling and simplifies theoretical analysis. Therefore, we assume $b=1$ in the following sections to establish the theoretical properties.

The gSIW prior is conjugate to a multivariate normal distribution. The posterior density of $(\Lambda,\rmU)$ is given by:

\begin{equation*}
    \pi(\Lambda,\rmU\vert \bfX_n) (d\Lambda) (d\rmU)\propto \prod\limits_{i=1}^p \lambda_i^{-a_i-n/2}etr(-\dfrac{1}{2}\Lambda^{-1}\rmU^T(H+nS)\rmU)(d\Lambda) (d\rmU),
\end{equation*}

where $\bfX_n = (X_1, \ldots, X_n)$, and $S$ is the sample covariance. Consider the following spectral decomposition of $nS= QWQ^T$, where $W=\text{diag}(n\hat{\lambda}_1,\ldots,n\hat{\lambda}_{n\wedge p},0,\ldots,0)$ and $Q$ is a $p\times p$ matrix whose $i$-th column is the eigenvector corresponding to the $i$-th eigenvalue. If we denote $\Gamma = Q^T\rmU$ and $H=hI_p$ for some positive $h>0$, then  $(d\rmU) =(d\Gamma)$ and  we get 
\begin{equation*}
    \pi(\Lambda,\Gamma\vert \bfX_n) (d\Lambda) (d\Gamma)\propto \prod\limits_{i=1}^p \lambda_i^{-a_i-n/2}etr(-\dfrac{1}{2}\Lambda^{-1}\Gamma^T(hI_p+W)\Gamma)(d\Lambda) (d\Gamma).
\end{equation*} 

To obtain the posterior convergence rate, we need the following expression, 
\begin{align*}
    \pi({\Gamma}\vert \bfX_n)( d{\Gamma}) &\propto \bigg(\displaystyle\int \pi(\Lambda,\Gamma\vert\bfX_n) (d\Lambda) \bigg)(d{\Gamma})\\
    &\propto \prod\limits_{i=1}^p (2^{a_i+n/2-1}\Gamma(a_i+n/2-1))\cdot  \prod\limits_{i=1}^p  c_i^{-a_i-n/2+1}  (d{\Gamma}),
\end{align*}
where $c_i$ is the $(i,i)$-th element of ${\Gamma}^T(hI_p+W){\Gamma}$.

\subsection{Sampling Method}\label{subsec:sampling}
3  Since the structure of the gSIW prior is closely related to the SIW prior, we employ the Gibbs sampling algorithm proposed by \cite{berger2020bayesian}. For completeness, we provide the sampling algorithm below using our notation. 

The posterior density is given by
\begin{equation*}
    \pi(\Lambda,\Gamma\vert \bfX_n)\propto \prod\limits_{i=1}^p \lambda_i^{-a_i-n/2}etr(-\dfrac{1}{2}\Lambda^{-1}\Gamma^T(hI_p+W)\Gamma),
\end{equation*}
where $\lambda_1,\ldots,\lambda_p$ are unordered diagonal elements of $\Lambda$.

\begin{enumerate}
    \item Sampling $\Lambda$ given $(\Gamma, \mathbf{X}_n)$  
    $$
    \pi(\Lambda \vert \Gamma, \mathbf{X}_n) \propto \prod_{i=1}^p \lambda_i^{-a_i - n/2} \exp\left(-\frac{c_i}{2\lambda_i}\right),
    $$  
    where $c_i$ is the $(i,i)$-th element of $\Gamma^\top(h I_p + W) \Gamma$. Each $\lambda_i$ can be sampled independently from $\lambda_i \sim \mathrm{InvGam}(a_i + n/2 - 1, c_i/2)$.

    \item Sampling $\Gamma$ given $(\Lambda, \mathbf{X}_n)$  
    $$\pi(\Gamma \vert \Lambda, \mathbf{X}_n) \propto etr(-\frac{1}{2} \Lambda^{-1} \Gamma^T(h I_p + W)\Gamma).$$  
    We adopt the Gibbs sampling scheme proposed in \cite{berger2020bayesian},  
    which iteratively updates pairs of rows in $\Gamma$. The procedure for updating the first two rows is as follows:

    \begin{enumerate}[label=Step \Roman*.]
    \item \text{Rotation transformation.}\\
    To update the first two rows of $\Gamma$ while preserving orthogonality, we apply a signed Givens rotation
    $$\Gamma_{\text{new}} = diag(R, I_{p-2}) 
    \begin{pmatrix}
    \Gamma_{1} \\
    \Gamma_{2}
    \end{pmatrix},$$  
    where $\Gamma_{1}$ and $\Gamma_{2}$ denote the first two and remaining $p - 2$ rows of $\Gamma$, respectively. The matrix $R$ is a signed rotation matrix defined as  
    $$ R = 
    \begin{pmatrix}
    \epsilon_1 & 0 \\
    0 & \epsilon_2
    \end{pmatrix} R_\theta, \quad 
    R_\theta = 
    \begin{pmatrix}
    \cos \theta & -\sin \theta \\
    \sin \theta & \cos \theta
    \end{pmatrix},$$  
    where $\epsilon_1, \epsilon_2 \in \{-1, 1\}$ and $\theta \in (-\pi/2, \pi/2]$.

    \item \text{Conditional distribution of $\theta$.}\\
    Since $H_0 = h I_p + W$ is diagonal, it can be partitioned as  
    $$ H_0 = 
    \begin{pmatrix}
    H_1 & 0 \\
    0 & H_2
    \end{pmatrix}, $$  
    where $H_1 \in \mathbb{R}^{2 \times 2}$. The conditional distribution of $\theta$ becomes
    $$\pi(\theta \vert \Gamma_{1}, \Gamma_{2}, \Lambda,\bfX_n) 
    \propto etr(-\dfrac{1}{2} H_1 R_\theta \Gamma_{1} \Lambda^{-1} \Gamma_{1}^T R_\theta^T).$$
    
    Let the spectral decomposition of $\Gamma_{1} \Lambda^{-1} \Gamma_{1}^T$ be
    $$ \Gamma_{1} \Lambda^{-1} \Gamma_{1}^T = 
    \begin{pmatrix}
    \cos \omega & -\sin \omega \\
    \sin \omega & \cos \omega
    \end{pmatrix}
    \begin{pmatrix}
    s_1 & 0 \\
    0 & s_2
    \end{pmatrix}
    \begin{pmatrix}
    \cos \omega & \sin \omega \\
    -\sin \omega & \cos \omega
    \end{pmatrix},$$  
    with $s_1 > s_2$ and $\omega \in (-\dfrac{\pi}{2}, \dfrac{\pi}{2}]$. Then, the conditional distribution of $\theta$ simplifies to
    $$\pi(\theta \vert \Gamma_{1}, \Gamma_{2}, \Lambda,\bfX_n) \propto \exp\left(c \cos^2(\theta + \omega)\right),$$  
    where $c = -\frac{1}{2}\abs{(s_1 - s_2)(h_1 - h_2)}$.

    \item \text{Sampling via transformation.}\\
    Let $\alpha = \cos^2(\theta + \omega)$, then we obtain
    $$\pi(\alpha \vert \Gamma_{1}, \Gamma_{2}, \Lambda,\bfX_n) \propto 
    \exp(c \alpha)  \alpha^{-1/2} (1 - \alpha)^{-1/2}, \quad \alpha \in (0, 1),$$  
    which corresponds to a tilted $Beta(1/2, 1/2)$ distribution. 
    We sample $\alpha$ using rejection sampling with a Beta proposal. Once $\alpha$ is drawn, the corresponding pair of rows in $\Gamma$ is updated accordingly,  
    and this procedure is repeated for all row pairs. Finally, since $\Gamma = Q^T U$, the eigenvector matrix of $\Sigma$ is recovered as $U = Q \Gamma$.
    \end{enumerate}
\end{enumerate}

\subsection{Choice of the number of spiked eigenvalues $k$}
To implement the proposed gSIW prior, the number of spiked eigenvalues k must be specified in advance.  
To select $k$, we consider three methods.
First, we adopt the Watanabe-Akaike Information Criterion (WAIC; \citealt{watanabe2010asymptotic}),  
defined for each candidate value k as
$$
\mathrm{WAIC}(k) = -2 \sum_{i=1}^n \log \mathbb{E}_{\Sigma \vert \mathbf{X}_n^{(k)}}[p(X_i \vert \Sigma)] + 2 \sum_{i=1}^n \mathrm{Var}_{\Sigma \vert \mathbf{X}_n^{(k)}}[\log p(X_i \vert \Sigma)],
$$
where the posterior distribution is derived under the gSIW prior assuming k spiked eigenvalues.  
In computing the WAIC, we follow the second approach proposed by \cite{gelman2014understanding},
which estimates the effective number of parameters using the variance of the log predictive density.
The value of k is then selected by minimizing WAIC over $k = 1,\ldots, k_{max}$, where $k_{max}$ is a suitably chosen upper bound.
Second, we consider the Growth Ratio (GR) method proposed by \cite{ahn2013eigenvalue},  
which is based on the eigenvalues of sample covariance.  
For each k, the GR is defined as
$$
\mathrm{GR}(k) = \frac{\log\left(1 + \frac{\hat{\lambda}_k}{V(k)}\right)}{\log\left(1 + \frac{\hat{\lambda}_{k+1}}{V(k+1)}\right)},
$$
where $ V(k) = \sum_{j=k+1}^{n \wedge p} \hat{\lambda}_j $, and $ \hat{\lambda}_j $ denotes the $j$-th largest eigenvalue of the sample covariance matrix.  
The value of k is selected by maximizing GR over $k = 1,\ldots, k_{max}$.
Finally, we consider the information criterion $\mathrm{IC}_{p3}$ proposed by \cite{bai2002determining}, among several criteria introduced in their work, as it showed relatively strong empirical performance under our model setting. It is defined as 
$$
\mathrm{IC}_{p3}(k) = \log(\dfrac{1}{np}\abs{\abs{ X - X U_1 (U_1^\top U_1)^{-1} U_1^T }}_F^2 ) + k  \dfrac{\log(n \wedge p)}{n \wedge p},
$$
where $ X \in \mathbb{R}^{n \times p} $ denotes the observed data matrix, and $ U_1 \in \mathbb{R}^{p \times k} $ denotes the matrix of the estimated leading k eigenvectors. The value of k is selected by minimizing $\mathrm{IC}_{p3}$ over $k = 1,\ldots, k_{max}$.

In practice, all three methods provide good estimates when the spiked eigenvalues are large, and WAIC generally performs better even when the spiked eigenvalues are relatively small. Overall, WAIC exhibits more robust performance across a range of settings. Therefore, we recommend using WAIC for selecting $k$.

The code for the Gibbs sampling procedure and choice of $k$ are publicly available at \url{https://github.com/swpark0413/besiw}.

\section{Main Results}\label{sec:main_results}
Consider $n$ independent samples $X_1,\ldots,X_n$ from $N(0,\Sigma)$ where $\Sigma$ is a $p\times p$ positive definite matrix. We assume the gSIW prior as defined in \eqref{eq:SIW} with $b=1$. Under the gSIW prior, the posterior is given by
\begin{equation*}
    \pi(\Lambda,\Gamma\vert \bfX_n) (d\Lambda) (d\Gamma)\propto \prod\limits_{i=1}^p \lambda_i^{-a_i-n/2}etr(-\dfrac{1}{2}\Lambda^{-1}\Gamma^T(hI_p+W)\Gamma)(d\Lambda) (d\Gamma),
\end{equation*}
where $\lambda_1,\ldots,\lambda_p$ are the unordered eigenvalues of $\Sigma$. Suppose that $\Sigma$ is a spiked covariance with $k$ spiked eigenvalues. Let $\lambda_{0,1},\ldots,\lambda_{0,p}$ be the ordered eigenvalues of the true covariance $\Sigma_0$. To obtain the posterior convergence rates for eigenstructure, we assume the following conditions:

\begin{enumerate}[label={A\arabic*.}]
    \item We consider high-dimensional settings where $n/p\rightarrow 0$.
    \item There exist positive constants $c_0$ and $C_0$ such that, for all $n$,
    $$C_0>\lambda_{0,k+1}>\cdots>\lambda_{0,p}>c_0.$$
    \item The $k$ largest eigenvalues are sufficiently separated by a constant value $\delta_0>0$:
$$\dfrac{\lambda_{0,j}-\lambda_{0,j+1}}{\lambda_{0,j}}\geq\delta_0,\quad\forall j=1,\cdots ,k.$$
    \item For $k$ spiked eigenvalues,  the quantities $d_j=\dfrac{p}{n \lambda_{0,j}}$ are bounded above.
    \item We set the hyperparameters of the prior as follows:
    $$a_1\leq \cdots\leq a_k\leq a_{k+1} = \cdots = a_n\leq a_{n+1}=\cdots = a_p,\quad a_k\prec n,\quad H=hI_p,\;\text{with}\; h<n.$$
\end{enumerate}

The conditions $A1-A4$ are used to establish the asymptotic properties of the sample covariance as in \cite{wang2017asymptotics}. The condition $A2$ implies that the non-spiked eigenvalues are bounded away from zero and from above. The condition $A4$ requires that $\lambda_{0,j}$ is at least of order $\dfrac{p}{n}$. The conditions we impose differ slightly from those considered in \cite{lee2024posterior}, which assume that $\dfrac{k^3}{n} \rightarrow 0$ and $\dfrac{\lambda_{0,k+1}p}{\lambda_{0,k} n} \rightarrow 0$. In their setting, the number of spiked eigenvalues $k$ is allowed to diverge, and the non-spiked eigenvalues are not necessarily bounded by a constant. In contrast, we assume that $k$ is fixed and impose a boundedness condition on the non-spiked eigenvalues. The condition $A5$ specifies the hyperparameters of the proposed prior and serves as a sufficient condition for verifying posterior convergence.

We demonstrate the posterior convergence rates for eigenvalues and eigenvectors when using the gSIW prior. The proof relies on the asymptotic properties of sample covariance of \cite{wang2017asymptotics} and will be provided in the Appendix.

To sketch the proof, we first introduce the shrinkage behavior of the posterior distribution. Define the subset $D_\epsilon$ of $O(p)$ as
\begin{align}\label{eq:D_epsilon}
    D_{\epsilon} = \bigg\{\Gamma\in O(p):\displaystyle\inf_{Q_2\in O_{p-k}} \abs{\abs{\begin{pmatrix}
       I_k & 0 \\ 
        0& Q_2
    \end{pmatrix}-\Gamma}}_F <\epsilon\bigg\}.
\end{align}
This subset $D_\epsilon$ represents a neighborhood around matrices $\Gamma$ whose first $k$ columns are close to $(I_k, 0)^T \in \mathbb{R}^{p \times k}$. Recalling that the eigenvectors of $\Sigma$ can be recovered as $U = Q\Gamma$, the set $D_\epsilon$ includes matrices whose leading $k$ columns in $U$ are approximately aligned with the top $k$ eigenvectors of the sample covariance matrix.

Under suitable conditions, we show that
$$\displaystyle\int_{{D_\epsilon^c}} \pi(\Gamma\vert \bfX_n)(d{\Gamma}),$$
is sufficiently small, indicating that the posterior distribution for $\Gamma$ is concentrated near $D_\epsilon$. This is formalized in Lemma \ref{lem:shrink_post_diff_a}.

In summary, when the joint posterior $\pi(\Lambda, \Gamma \vert \mathbf{X}_n)$ is integrable with respect to $\Lambda$, 
posterior expectations can be approximated by identifying a subset $A \subset O(p)$ for which
$$\int_A \int \pi(\Lambda,\Gamma\vert \bfX_n) (d\Lambda) (d\Gamma) \approx 1.$$

This idea also extends to other models. For example, consider the model
\begin{align*}
    X_1, \ldots, X_n &\overset{iid}{\sim} N(\mu, \sigma^2)\\
    \mu &\sim N(\mu_0, \sigma_0^2)\\
    \sigma^2 &\sim \pi(\sigma^2),
\end{align*}
where $\mu_0\in\bbR$ and $\sigma_0^2>0$. The prior $\pi(\sigma^2)$ is supported on $(0, \infty)$. Then the marginal posterior $\pi(\mu, \sigma^2 \vert \mathbf{X}_n)$ is integrable in $\mu$, 
so for a set $A \subset (0, \infty)$,
$$\int_A \int \pi(\mu,\sigma^2\vert \bfX_n) d\mu d\sigma^2 \approx 1,$$
allowing direct approximation of the posterior expectation.

\begin{lemma}\label{lem:shrink_post_diff_a} 
Under model \eqref{eq:model} and prior \eqref{eq:SIW}, assume that conditions $A1-A5$ hold. Let $\epsilon > 0$, and define $\kappa = \min\limits_{l<k}(a_{l+1} - a_l) \cdot \min\limits_{l<k}\log(\dfrac{\lambda_{0,l}}{\lambda_{0,l+1}})$. Suppose that $\dfrac{np}{a_n}\prec \dfrac{\lambda_{0,k}}{n\lambda_{0,1}}\dfrac{\kappa}{2}\epsilon$ and $\epsilon\succ \sqrt{\dfrac{np}{a_n}}\vee \kappa^{-1}$. Then,
    \begin{equation}\label{eq:shrink_post}
        \displaystyle\int_{{D_\epsilon^c}} \pi(\Gamma\vert \bfX_n)(d{\Gamma})= O\Big(\exp(-\dfrac{\kappa}{2}\epsilon)  \Big)  + O\Big(\big(\dfrac{n\lambda_{0,k}+p}{p}\big)^{-\epsilon^2 a_n}\Big),
    \end{equation}
    where $D_\epsilon$ is defined in \eqref{eq:D_epsilon}.
\end{lemma}
\hfill\break


While Lemma \ref{lem:shrink_post_diff_a} is useful when the prior parameters satisfy $a_1 < \cdots < a_k$, it becomes less informative when $a_1, \ldots, a_k$ are nearly identical. In such cases, the quantity $\kappa$ approaches zero, thereby weakening the exponential decay term in \eqref{eq:shrink_post}, and hence failing to adequately describe posterior contraction.

To address this limitation, we consider two prior parameter choices for the asymptotic analysis:

\begin{itemize}
    \item When $a_1 < \cdots < a_k$, Lemma \ref{lem:shrink_post_diff_a} provides a direct bound on the posterior mass outside $D_\epsilon$.
    \item When $a_1, \ldots, a_k$ are approximately equal, we reparametrize $(\Lambda, \Gamma)$ and analyze the posterior concentration around certain structured subsets (to be defined later as $D_{\epsilon,l}$), formalized in Lemma \ref{lem:shrink_post_equal_a}.
\end{itemize}

Assume that $a_1, \ldots, a_k$ are approximately equal. For $\Gamma\in O(p)$, we determine the $k\times k$ permutation matrix $P_l$ by solving the following optimization problem
$$P_l = \argmin\limits_{P\in\mathbf{P}}\Bigg[\displaystyle\inf_{Q_2\in O_{p-k}} \abs{\abs{\begin{pmatrix}
           P & 0 \\ 
            0& Q_2
        \end{pmatrix}-\Gamma}}_F \Bigg],$$
where $\mathbf{P}$ is the set of all $k\times k$ permutation matrices. Once $P_l$ is determined, we consider the following reparametrization on $\Lambda$ and $\Gamma$,
\begin{align}\label{eq:spiked_repara}
    \tilde{\Lambda} &= \begin{pmatrix}
           P_l & 0 \\ 
            0& I_{p-k}
        \end{pmatrix}\Lambda \begin{pmatrix}
           P_l & 0 \\ 
            0& I_{p-k}
        \end{pmatrix}^T, \quad 
    \tilde{\Gamma} = \Gamma  \begin{pmatrix}
           P_l & 0 \\ 
            0 & I_{p-k}
        \end{pmatrix}^T.
\end{align}

Since the transformation preserves the measure, it follows that
$$(d\Gamma) = (d\tilde{\Gamma}).$$

By reparametrizing the variables, the posterior distribution becomes
$$\pi(\tilde{\Lambda},\tilde{\Gamma}\vert\bfX_n)(d\tilde{\Lambda})(d\tilde{\Gamma}) = \prod\limits_{i=1}^p\tilde{\lambda}_i^{-a_i-n/2} etr(-\dfrac{1}{2}\tilde{\Lambda}^{-1}\tilde{\Gamma}^T(hI_p + W)\tilde{\Gamma})(d\tilde{\Lambda})(d\tilde{\Gamma}).$$

For each $\Gamma \in O(p)$, the optimal permutation $P_l$ can be selected, allowing us to define the subset $D_{\epsilon,l} \subset O(p)$ as
\begin{align}\label{eq:D_epsilon_perm}
    D_{\epsilon,l} = \bigg\{\Gamma\in O(p):\displaystyle\inf_{Q_2\in O_{p-k}} \abs{\abs{\begin{pmatrix}
       P_l & 0 \\ 
    0& Q_2
    \end{pmatrix}-\Gamma}}_F <\epsilon\bigg\},
\end{align}
where $P_1,\ldots,P_L$ denotes all $k\times k$ permutation matrices.

Consequently, we have
$$\displaystyle\int_{\bigcup_{l=1}^L D_{\epsilon,l}} \pi(\Gamma\vert\bfX_n)(d{\Gamma}) \approx \displaystyle\int_{D_\epsilon} \pi(\tilde{\Gamma}\vert\bfX_n)(d{\tilde{\Gamma}}).$$

Therefore, even when $a_1, \ldots, a_k$ are approximately equal, we can analyze the posterior probability over $(\cup_{l=1}^L D_{\epsilon,l})^c$ analogously to Lemma \ref{lem:shrink_post_diff_a}. Under suitable conditions, we can still show that
$$\displaystyle\int_{(\bigcup_{l=1}^L D_{\epsilon,l})^c} \pi(\Gamma\vert \bfX_n)(d{\Gamma}),$$
is sufficiently small, implying that the posterior for $\Gamma$ is well concentrated near 
$\cup_{l=1}^L D_{\epsilon,l}$. This is formalized in the Lemma \ref{lem:shrink_post_equal_a}.

\begin{lemma}\label{lem:shrink_post_equal_a}
    Under model \eqref{eq:model} and prior \eqref{eq:SIW}, assume that conditions $A1-A5$ hold, and $\max_{i\leq k} a_i - \min_{i\leq k }a _i \leq q$ for some positive $q$, which may depend on $n$ and $p$. Let $\epsilon>0$, and define $\tau = \min\limits_{l<k}\log(\dfrac{\lambda_{0,l}}{\lambda_{0,l+1}})$. Suppose that $\dfrac{np}{a_n}\prec  \dfrac{\lambda_{0,k}}{n\lambda_{0,1}}\dfrac{\tau}{2}n\epsilon^2$
    and $\epsilon\succ \sqrt{\dfrac{np}{a_n}}\vee (n\tau)^{-1/2}$. Then,
    \begin{equation}\label{eq:shrink_post2}
        \displaystyle\int_{(\bigcup_{l=1}^L D_{\epsilon,l})^c} \pi(\Gamma\vert\bfX_n)(d{\Gamma})=  O\Big(\exp(-\dfrac{\tau}{2}n\epsilon^2)\cdot\Big(\dfrac{2\lambda_{0,1}}{\lambda_{0,k}}\Big)^{kq}\Big)  + O\Big(\big(\dfrac{n\lambda_{0,k}+p}{p}\big)^{-\epsilon^2 a_n}\Big),
    \end{equation}
    where $D_{\epsilon,l}$ is defined in \eqref{eq:D_epsilon_perm}.
\end{lemma}
\hfill\break

Lemma \ref{lem:shrink_post_diff_a} and \ref{lem:shrink_post_equal_a} describe the degree of posterior shrinkage on the sets of orthogonal matrices $D_\epsilon$ and $\cup_{l=1}^L D_{\epsilon,l}$, respectively. These lemmas provide the foundation for analyzing posterior expectation. Consider the posterior expectation of a general nonnegative function $f(\Gamma, \Lambda)$, defined as
\begin{align}\label{eq:post_exp}
    \bbE\big[f(\Gamma,\Lambda)\vert\bfX_n\big]=\dfrac{\displaystyle\int f'(c,\Gamma)\prod c_{i}^{-a_i-\frac{n}{2}+1}(d{\Gamma})}{\displaystyle\int\prod c_{i}^{-a_i-\frac{n}{2}+1} (d{\Gamma})},
\end{align}
where $f'(c,\Gamma) = \mathbb{E}_\Lambda[f(\Gamma,\Lambda)]$ and $\lambda_i \overset{\text{ind}}{\sim} \mathrm{InvGam}(a_i + n/2 - 1, c_i/2)$.

In the case where $a_1 < \cdots < a_k$, an upper bound on the posterior expectation can be derived as
\begin{align*}
    \eqref{eq:post_exp} & = \dfrac{\displaystyle\int_{D_\epsilon} f'(c,\Gamma)\prod c_{i}^{-a_i-\frac{n}{2}+1}(d{\Gamma})}{\displaystyle\int\prod c_{i}^{-a_i-\frac{n}{2}+1} (d{\Gamma})} + \dfrac{\displaystyle\int_{D_\epsilon^c} f'(c,\Gamma)\prod c_{i}^{-a_i-\frac{n}{2}+1}(d{\Gamma})}{\displaystyle\int\prod c_{i}^{-a_i-\frac{n}{2}+1} (d{\Gamma})}\\
    &\leq \sup_{\Gamma\in D_\epsilon} f'(c,\Gamma) + \sup_{\Gamma\in O(p)} f'(c,\Gamma) \cdot \displaystyle\int_{{D_\epsilon^c}} \pi(\Gamma\vert \bfX_n)(d{\Gamma}).
\end{align*}

Similarly, a lower bound is given by
\begin{align*}
    \eqref{eq:post_exp}  & = \dfrac{\displaystyle\int_{D_\epsilon} f'(c,\Gamma)\prod c_{i}^{-a_i-\frac{n}{2}+1}(d{\Gamma})}{\displaystyle\int\prod c_{i}^{-a_i-\frac{n}{2}+1} (d{\Gamma})} + \dfrac{\displaystyle\int_{D_\epsilon^c} f'(c,\Gamma)\prod c_{i}^{-a_i-\frac{n}{2}+1}(d{\Gamma})}{\displaystyle\int\prod c_{i}^{-a_i-\frac{n}{2}+1} (d{\Gamma})}\\
    &\geq \inf_{\Gamma\in D_\epsilon} f'(c,\Gamma) + \inf_{\Gamma\in O(p)}  f'(c,\Gamma) \cdot \displaystyle\int_{{D_\epsilon^c}} \pi(\Gamma\vert \bfX_n)(d{\Gamma}).
\end{align*}

When $a_1, \ldots, a_k$ are approximately equal, the same logic holds with $D_\epsilon$ replaced by $\cup_{l=1}^L D_{\epsilon,l}$.

Suppose the function $f'$ is sufficiently bounded such that 
$$\sup\abs{f'}\prec
\begin{cases}
     \exp(\dfrac{\tau}{2}n\epsilon^2)\cdot\Big(\dfrac{2\lambda_{0,1}}{\lambda_{0,k}}\Big)^{-kq}\wedge \Big(\dfrac{n\lambda_{0,k}+p}{p}\Big)^{\epsilon^2 a_n},\quad &\text{if }a_1\approx \cdots\approx a_k,\\
     \exp(\dfrac{\kappa}{2}\epsilon)\wedge \Big(\dfrac{n\lambda_{0,k}+p}{p}\Big)^{\epsilon^2 a_n},\quad &\text{if }a_1<\cdots<a_k.
\end{cases}$$
Under this condition, Lemma \ref{lem:shrink_post_diff_a} or \ref{lem:shrink_post_equal_a} can be applied to derive the posterior expectation, thereby establishing the posterior convergence rate. The following corollaries demonstrate specific cases where the posterior expectation can be explicitly derived.

\begin{theorem}\label{cor:conv_rate_eigenvalue_equal_a}
Under the assumptions of Lemma \ref{lem:shrink_post_equal_a}, 
let $\lambda_{(i)}$ denote the $i$-th largest eigenvalue of $\Sigma$. 
Suppose that the ratio $\lambda_{0,1}/\lambda_{0,k}$ is bounded by a positive constant, 
$\epsilon \asymp n^{-1/4}$,  $a_1,\ldots,a_k \preccurlyeq n^{1/2}$, and $a_n\succ n^{3/2}p$. 
Then, 
$$\bbE\Big[\dfrac{\lambda_{(i)}-\lambda_{0,i}}{{\lambda_{0,i}}}\big\vert \bfX_n \Big] = O(\dfrac{p}{n\lambda_{0,i}})+ O(n^{-1/2+\delta}),$$
for $i = 1,\ldots,k$ and for all small $\delta > 0$.
\end{theorem}
\hfill\break

The convergence rate stated in Theorem \ref{cor:conv_rate_eigenvalue_equal_a} is equivalent to the rate achieved by the inverse-Wishart prior, as established by \cite{lee2024posterior}. To achieve faster rates than those given in Theorem \ref{cor:conv_rate_eigenvalue_equal_a}, we set the hyperparameter $a_i$ depending on the eigenvalues of the sample covariance.

\begin{theorem}\label{cor:conv_rate_eigenvalue_diff_a}
Under the assumptions of Lemma \ref{lem:shrink_post_diff_a} and \ref{lem:shrink_post_equal_a}, let $\lambda_{(i)}$ denote the $i$-th largest eigenvalue of $\Sigma$. 
Suppose that the ratio $\lambda_{0,1} / \lambda_{0,k}$ is bounded by a positive constant, 
$\epsilon \asymp n^{-1/4}$, 
and the parameters satisfy $2a_i - 4 = \frac{nt}{\hat{\lambda}_i - t}, \quad \text{for } i = 1, \ldots, k,$ with $t \in [\hat{\lambda}_{k+1}, \hat{\lambda}_n]$, where $\hat{\lambda}_i$ denotes the $i$-th eigenvalue of the sample covariance matrix, and $a_n \succ n^{3/2}p$. Then,
  $$\bbE\Big[\dfrac{\lambda_{(i)}-\lambda_{0,i}}{{\lambda_{0,i}}}\big\vert \bfX_n \Big] 
  =  O(\dfrac{1}{\lambda_{0,i}}\sqrt{\dfrac{p}{n}})+ O(n^{-\frac{1}{2}+\delta}),$$
for  $i = 1,\ldots,k$ and for all small $\delta > 0$.
\end{theorem} 
\hfill\break

Theorem \ref{cor:conv_rate_eigenvalue_diff_a} provides a posterior convergence rate that is better or at least equivalent compared to that described in Theorem \ref{cor:conv_rate_eigenvalue_equal_a}. If the condition $\lambda_{0,i} \prec \dfrac{p}{\sqrt{n}}$ holds, then the convergence rate given in Theorem \ref{cor:conv_rate_eigenvalue_diff_a} is faster than that of the prior described in Theorem \ref{cor:conv_rate_eigenvalue_equal_a}. The rate coincides with the shrinkage eigenvalues proposed by \cite{wang2017asymptotics}, given by:
$$\hat{\lambda}_j^S = \hat{\lambda}_j - \dfrac{p}{np-nk-pk} \Big(tr(S)-\sum\limits_{i=1}^k\hat{\lambda}_i\Big).$$

Let $\xi_j$ be the eigenvector corresponding to $\lambda_j$, the $j$-th eigenvalue of the covariance $\Sigma$. The eigenvector corresponding to the true covariance $\lambda_{0,i}$ is denoted by $\xi_{0,j}$. Theorem \ref{cor:conv_rate_eigenvector} presents the posterior convergence rate for the first $k$ eigenvectors.

\begin{theorem}\label{cor:conv_rate_eigenvector}
Under the assumptions of Theorem \ref{cor:conv_rate_eigenvalue_diff_a}, let $\xi_{(j)}$ denote the eigenvector corresponding to the $j$-th largest eigenvalue $\lambda_{(j)}$ of $\Sigma$. Then,
$$\bbE\Big[1-\abs{\xi_{0,j}^T\xi_{(j)}}^2 \big\vert \bfX_n \Big] = O (\dfrac{p}{n\lambda_{0,j}})+ O_p(\zeta_j),$$
for $j=1,\ldots,k$, where $\zeta_j = \dfrac{1}{\lambda_{0,j}}\sqrt{\dfrac{p}{n}}+ \dfrac{p}{n^{3/2}\lambda_{0,j}}+\dfrac{1}{n}$.
\end{theorem}

\hfill\break

This convergence rate is equivalent to the rate of the sample eigenvector as an estimator of the first $k$ eigenvectors. According to the study of \cite{cai2013sparse}, this rate is proven to be optimal. 

\begin{corollary}\label{cor:conv_rate_cov}
Under the assumptions of Theorem \ref{cor:conv_rate_eigenvalue_diff_a},
\begin{align}\label{eq:conv_rate_cov}
    \bbE\Big[\dfrac{\abs{\abs{\Sigma-\Sigma_0}}_F^2}{\abs{\abs{\Sigma_0}}_F^2}\big\vert \bfX_n \Big] = O\Big(\dfrac{p }{\lambda_{0,1}^2 + p} + n^{-1+2\delta}\dfrac{\lambda_{0,1}^2}{\lambda_{0,1}^2 + p}\Big)+ O \Big(\dfrac{p}{n\lambda_{0,k}}\dfrac{\lambda_{0,1}^2}{\lambda_{0,1}^2 + p}\Big)+O_p(\dfrac{\lambda_{0,1}^2}{\lambda_{0,1}^2 + p}\zeta_k),
\end{align}
for all small $\delta>0$, where $\zeta_k = \dfrac{1}{\lambda_{0,k}}\sqrt{\dfrac{p}{n}}+ \dfrac{p}{n^{3/2}\lambda_{0,k}}+\dfrac{1}{n}$. Furthermore, for the sample covariance, we obtain the following rate
\begin{align}\label{eq:conv_rate_samplecov}
\dfrac{\abs{\abs{S-\Sigma_0}}_F^2}{\abs{\abs{\Sigma_0}}_F^2} = O\Big(\dfrac{p^2}{n(\lambda_{0,1}^2 + p)}+n^{-1+2\delta}\dfrac{ \lambda_{0,1}^2}{\lambda_{0,1}^2+p}\Big) + O\Big(\dfrac{p}{n\lambda_{0,k}} \dfrac{\lambda_{0,1}^2}{\lambda_{0,1}^2 + p} \Big) + O_p\Big(\dfrac{\lambda_{0,1}^2}{\lambda_{0,1}^2 + p}\zeta_k\Big).
\end{align}
\end{corollary}

The first term on the right-hand side of \eqref{eq:conv_rate_cov} corresponds to the convergence rate of the eigenvalues, while the second and third terms are attributed to the convergence of the eigenvectors. Comparing the rate in \eqref{eq:conv_rate_cov} with the sample covariance rate in \eqref{eq:conv_rate_samplecov}, we observe that the first term in the sample covariance rate is larger.

\section{Simulation Studies}\label{sec:simul_studies}

To evaluate the performance of the proposed gSIW prior, we consider observations sampled from a multivariate normal distribution $N(0, \Sigma_0)$. We design two types of simulation settings and compare six different methods. Throughout the simulations, we assume the true number of spiked eigenvalues is known and focus on estimating the spiked eigenvalues and their corresponding eigenvectors. Each simulation is repeated $100$ times.

\begin{enumerate}[label={Method \Roman*.}]
    \item Sample covariance (Sample).
    \item Generalized shrinkage inverse-Wishart prior (gSIW) :
    $$\pi(\Sigma\vert a,H) \propto \dfrac{etr(-\dfrac{1}{2}\Gamma\Lambda^{-1} \Gamma^TH)}{\prod\limits_{i=1}^p\lambda_i^{a_i}\prod\limits_{i<j}\abs{\lambda_i-\lambda_j}}.$$
    \item Generalized inverse-Wishart prior (gIW) : 
    $$\pi(\Sigma\vert a,H) \propto \dfrac{etr(-\dfrac{1}{2}\Gamma\Lambda^{-1} \Gamma^TH)}{\prod\limits_{i=1}^p\lambda_i^{a_i}},$$
    where the hyperparameters $a$ and $H$ are set identically for both gSIW and gIW. In particular, we use $a_i = \dfrac{nt}{2(\hat{\lambda}_i - t)} + 2$ for $i = 1, \ldots, k$, where $t = \dfrac{1}{n-k} \sum_{i=1}^{n-k} \hat{\lambda}_i$, $a_n = p/2$, $a_p = 2p$, and $H = 4I_p$.
    \item Shrinkage inverse-Wishart prior proposed by \cite{berger2020bayesian} (SIW) :
    $$\pi(\Sigma\vert a,H) \propto \dfrac{etr(-\dfrac{1}{2}\Gamma\Lambda^{-1} \Gamma^TH)}{\prod\limits_{i=1}^p\lambda_i^{a}\prod\limits_{i<j}\abs{\lambda_i-\lambda_j}},$$
    where $a=4$ and $H=4I_p$.
    \item Inverse-Wishart prior (IW) :
    \begin{align*}
        \pi(\Sigma\vert a,H) \propto \dfrac{etr(-\dfrac{1}{2}\Gamma\Lambda^{-1} \Gamma^TH)}{\prod\limits_{i=1}^p\lambda_i^{a}},
    \end{align*}
    where $a = p+1$ and $H = I_p$, as specified by \cite{lee2024posterior}.
    \item Shrinkage POET estimator proposed by \cite{wang2017asymptotics} (S-POET).
\end{enumerate}

To evaluate performance, we use the following errors for spiked eigenvalues and eigenvectors:
$$\text{Err}_\lambda = \dfrac{\abs{\lambda_{0,i}-\lambda_i}}{\lambda_{0,i}},\quad \text{Err}_\xi= 1-(\xi_{i}^T\xi_{i,0})^2.$$
For the Bayesian methods, the posterior mean is used as the point estimator. Eigenvector estimates are obtained by aligning signs, taking the Euclidean average, and then normalizing. We also report the $95\%$ credible (or confidence) interval length (IR) and the coverage probability (CP) for each spiked eigenvalue. The coverage probability is defined as the proportion of repetitions in which the true value lies within the credible (or confidence) interval.

For the S-POET method, the confidence intervals are derived from the asymptotic distribution
$$\sqrt{n}\Big(\dfrac{\hat{\lambda}_i^{S}}{\lambda_{0,i}}-1\Big)\overset{d}{\Rightarrow} N(0,2),\quad \text{if } \sqrt{p}=o(\lambda_{0,i}),$$
where $\hat{\lambda}_i^S$ is the shrinkage eigenvalue, as established by \cite{wang2017asymptotics}.

\subsection{Case 1}
In this simulation, we set $(n,p) = (50,500)$ with $\Sigma_0 = \text{diag}(50,20,10,1,\ldots,1)$. For the $\Sigma_0$, we obtain $d_1 = 0.2$, $d_2=0.5$, and $d_3=1$ where $d_j =\dfrac{p}{n\lambda_{0,j}}$. This setting satisfies assumption $A4$, as the values of $d_j$ fall within the specified range. The goal of this setting is to assess the performance of the methods when the dimension $p$ is large and the spiked eigenvalues are significantly larger than the others.

\begin{table}[t!]
\begin{center}
{\scriptsize
\setlength{\tabcolsep}{4pt}
\begin{tabular}{cc ccc ccc ccc ccc ccc}
\hline
\multirow{2}{*}[-3pt]{Eigenvalue} & \multicolumn{1}{c}{Sample} 
& \multicolumn{3}{c}{gSIW} & \multicolumn{3}{c}{gIW} 
& \multicolumn{3}{c}{SIW} & \multicolumn{3}{c}{IW} 
& \multicolumn{3}{c}{S-POET} \\
\cmidrule(lr){2-2} \cmidrule(lr){3-5} \cmidrule(lr){6-8}
\cmidrule(lr){9-11} \cmidrule(lr){12-14} \cmidrule(lr){15-17}
& Err$_\lambda$ 
& Err$_\lambda$ & CP & IL 
& Err$_\lambda$ & CP & IL 
& Err$_\lambda$ & CP & IL 
& Err$_\lambda$ & CP & IL 
& Err$_\lambda$ & CP & IL\\
\hline
$\lambda_1$ & 0.258 & \textbf{0.164} & 0.88 & 36.1 & 5.28 & 0.01 & 608 & 0.173 & 0.85 & \textbf{35.3} & 0.589 & 0.06 & 56.2 & 0.170 & \textbf{0.92} & 48.3\\
$\lambda_2$ & 0.517 & 0.149 & 0.95 & \textbf{13.7} & 2.23 & 0.01 & 77.7 & 0.165 & 0.94 & 14.6 & 1.59 & 0 & 26.7 & \textbf{0.128} & \textbf{0.98} & 19.5\\
$\lambda_3$ & 1.05 & \textbf{0.141} & 0.91 & \textbf{6.14} & 1.28 & \textbf{0.05} & 19.5 & 0.289 & 0.65 & 6.93 & 3.13 & 0 & 17.2 & 0.174 & \textbf{0.94} & 10.7\\
\hline
\end{tabular}
}
\caption{Average errors (Err$_\lambda$), coverage probabilities (CP), and interval lengths (IL) for the estimated spiked eigenvalues.}
\label{tbl:eigenval_case1}
\end{center}
\end{table}

Table \ref{tbl:eigenval_case1} presents the average errors (Err$_\lambda$), coverage probabilities (CP), and credible (or confidence) interval length (IL) for the estimated spiked eigenvalues across different methods. The IW and gIW methods perform poorly in both estimation accuracy and coverage. In contrast, SIW, S-POET, and gSIW outperform the other methods.  Among these, the gSIW achieves the lowest errors for the first and third eigenvalues, while S-POET performs best for the second eigenvalue. 
The gSIW, SIW, and S-POET methods exhibit comparable coverage probabilities. Although SIW performs slightly worse and S-POET slightly better than gSIW, the credible intervals of gSIW are narrower than those of S-POET. Furthermore, gSIW achieves lower estimation errors than SIW for all spiked eigenvalues.

\begin{table}[t!]
\begin{center}
{\scriptsize
\setlength{\tabcolsep}{4pt}
\begin{tabular}{ccccccc}
\hline
Eigenvector & Sample & gSIW & gIW & SIW & IW & S-POET \\
\hline
$\xi_1$ & \textbf{0.202} & 0.227 & 0.223 & 0.203 & 0.221 & 0.206 \\
$\xi_2$ & \textbf{0.419} & 0.443 & 0.466 & 0.427 & 0.673 & 0.428 \\
$\xi_3$ & 0.639 & 0.642 & 0.658 & 0.679 & 0.921 & \textbf{0.620} \\
\hline
\end{tabular}
}
\caption{Average errors (Err$_\xi$) for estimated eigenvectors.}
\label{tbl:eigenvec_case1}
\end{center}
\end{table}

Table~\ref{tbl:eigenvec_case1} presents the average errors (Err$_\xi$) for the estimated eigenvectors.  
The sample covariance achieves the lowest errors for the first and second spiked eigenvectors, while the S-POET method demonstrates the best performance for the third spiked eigenvector. Overall, the differences among the methods are relatively small.

\subsection{Case 2}
In this simulation, we vary the sample size $n$ and $p$ as 
$$(n,p) = (20,100), (40,100), \ldots, (80,100),$$ 
while fixing the covariance as $\Sigma_0 = \text{diag}(5, 4, 3, 1, \ldots, 1)$. 
We assume that the true number of spiked eigenvalues is 3. As n increases, the quantity $d_j = \dfrac{p}{n\lambda_{0,j}}$ decreases, in accordance with assumption A4, which holds more strongly for larger $n$. The goal of this setting is to evaluate the performance of the methods when the spiked eigenvalues are relatively weak.

\begin{table}[t!]
\begin{center}
{\scriptsize
\setlength{\tabcolsep}{4pt}
\begin{tabular}{cc ccc ccc ccc ccc ccc ccc}
\hline
\multirow{2}{*}[-3pt]{Eigenvalue} & \multirow{2}{*}[-3pt]{n}
& \multicolumn{1}{c}{Sample} 
& \multicolumn{3}{c}{gSIW} & \multicolumn{3}{c}{gIW} 
& \multicolumn{3}{c}{SIW} & \multicolumn{3}{c}{IW} 
& \multicolumn{3}{c}{S-POET} \\
\cmidrule(lr){3-3} \cmidrule(lr){4-6} \cmidrule(lr){7-9}
\cmidrule(lr){10-12} \cmidrule(lr){13-15} \cmidrule(lr){16-18}
& & Err$_\lambda$ 
& Err$_\lambda$ & CP & IL 
& Err$_\lambda$ & CP & IL 
& Err$_\lambda$ & CP & IL 
& Err$_\lambda$ & CP & IL 
& Err$_\lambda$ & CP & IL\\
\hline
\multirow{4}{*}{$\lambda_1$} & 20 & 1.45 & 0.422 & 0.69 & 6.27 & 2.80 & 0 & 29.3 & \textbf{0.110} & \textbf{1} & \textbf{5.37} & 4.23 & 0 & 27.4 & 0.670 & 0.42 & 16.8\\
& 40 & 0.732 & \textbf{0.152} & \textbf{0.96} & \textbf{3.81} & 3.00 & 0 & 27.4 & 0.154 & \textbf{0.96} & 3.90 & 2.16 & 0 & 10.8 & 0.402 & 0.55 & 7.60\\
& 60 & 0.499 & \textbf{0.102} & \textbf{0.98} & \textbf{3.11} & 3.64 & 0 & 30.9 & 0.132 & 0.96 & 3.50 & 1.41 & 0 & 6.66 & 0.307 & 0.64 & 5.36\\
& 80 & 0.387 & \textbf{0.0958} & \textbf{0.97} & \textbf{2.75} & 4.40 & 0 & 35.4 & {0.113} & {0.96} & {3.09} & 1.03 & 0 & 4.89 & 0.254 & 0.70 & 4.30\\
\hline
\multirow{4}{*}{$\lambda_2$} & 20 & 1.54 & 0.292 & 0.81 & 3.55 & 1.80 & 0 & 11.6 & \textbf{0.0793} & \textbf{1} & \textbf{2.72} & 3.23 & 0 & 12.7 & 0.555 & 0.64 & 12.5\\
& 40 & 0.783 & \textbf{0.105} & \textbf{0.97} & 2.43 & 2.10 & 0 & 12.0 & 0.165 & 0.88 & \textbf{2.08} & 1.95 & 0 & 6.08 & 0.390 & 0.65 & 6.03\\
& 60 & 0.525 & \textbf{0.0791} & \textbf{0.98} & \textbf{2.04} & 2.54 & 0 & 13.7 & 0.184 & 0.8 & 2.06 & 1.36 & 0 & 3.92 & 0.314 & 0.66 & 4.31\\
& 80 & 0.395 & \textbf{0.0895} & \textbf{0.98} & \textbf{1.86} & 3.07 & 0 & 15.6 & 0.163 & 0.78 & 2.10 & 1.03 & 0 & 2.98 & 0.262 & 0.72 & 3.46\\
\hline
\multirow{4}{*}{$\lambda_3$} & 20 & 1.98 & 0.322 & 0.84 & 2.83 & 1.54 & 0 & 7.41 & \textbf{0.0913} & \textbf{1} & \textbf{1.93} & 3.16 & 0 & 7.76 & 0.655 & 0.47 & 10.0\\
& 40 & 1.05 & 0.123 & 0.99 & 2.02 & 1.80 & 0 & 7.88 & \textbf{0.0544} & \textbf{1} & \textbf{1.43} & 2.18 & 0 & 4.16 & 0.533 & 0.17 & 4.99\\
& 60 & 0.700 & \textbf{0.0758} & \textbf{0.99} & 1.70 & 2.11 & 0 & 8.83 & 0.101 & 0.98 & \textbf{1.31} & 1.65 & 0 & 2.84 & 0.453 & 0.14 & 3.58\\
& 80 & 0.528 & \textbf{0.0812} & \textbf{0.97} & 1.54 & 2.49 & 0 & 9.92 & 0.145 & 0.88 & \textbf{1.30} & 1.31 & 0 & 2.18 & 0.404 & 0.09 & 2.89\\
\hline
\end{tabular}
}
\caption{Average errors (Err$_\lambda$), coverage probabilities (CP), and credible (or confidence) interval lengths (IL) for estimated eigenvalues under varying $n$.}
\label{tbl:eigenval_case2}
\end{center}
\end{table}

Table \ref{tbl:eigenval_case2} represents the average errors (Err$_\lambda$), coverage probabilities (CP), and credible (or confidence) interval length (IL) for the estimated spiked eigenvalues, across different methods and values of $n$. The gSIW, SIW, and S-POET methods consistently outperform the sample covariance. In contrast, the IW and gIW methods exhibit significantly poorer performance across all settings.
For small sample size $n=20$, the SIW method shows superior performance in terms of both error and coverage probability, while for larger sample sizes, gSIW achieves the best overall performance. In particular, the coverage probability of gSIW increases with $n$, indicating that the method becomes more reliable as sample size grows.
Furthermore, across all values of $n$, the credible intervals produced by gSIW are shorter than those of S-POET, highlighting the efficiency of gSIW in quantifying uncertainty. These results demonstrate the advantages of the gSIW prior, particularly for moderate to large sample sizes.

\begin{table}[t!]
\begin{center}
{\scriptsize
\setlength{\tabcolsep}{4pt}
\begin{tabular}{ccc ccccc}
\hline
             Eigenvector & n & Sample & gSIW & gIW & SIW & IW & S-POET \\
\hline
\multirow{4}{*}{$\xi_1$} & 20 & 0.747 & 0.788 & 0.784 & 0.791 & 0.828 & \textbf{0.737} \\
                         & 40 & \textbf{0.614} & 0.643 & 0.671 & 0.671 & 0.726 & 0.634 \\
                         & 60 & \textbf{0.543} & 0.560 & 0.607 & 0.565 & 0.638 & 0.578 \\
                         & 80 & \textbf{0.492} & 0.521 & 0.528 & 0.490 & 0.543 & 0.536 \\
\hline
\multirow{4}{*}{$\xi_2$} & 20 & 0.897 & 0.893 & 0.913 & 0.930 & 0.938 & \textbf{0.885} \\
                         & 40 & \textbf{0.794} & 0.837 & 0.835 & 0.886 & 0.896 & 0.798 \\
                         & 60 & \textbf{0.709} & 0.732 & 0.752 & 0.788 & 0.811 & 0.738 \\
                         & 80 & \textbf{0.650} & 0.680 & 0.703 & 0.657 & 0.786 & 0.698 \\
\hline
\multirow{4}{*}{$\xi_3$} & 20 & 0.960 & 0.961 & 0.948 & 0.958 & 0.961 & \textbf{0.952} \\
                         & 40 & \textbf{0.889} & 0.908 & 0.905 & 0.945 & 0.938 & 0.879 \\
                         & 60 & \textbf{0.802} & 0.813 & 0.831 & 0.917 & 0.902 & 0.806 \\
                         & 80 & \textbf{0.740} & 0.745 & 0.792 & 0.841 & 0.915 & 0.756 \\
\hline
\end{tabular}
}
\caption{Average errors (Err$_\xi$) for estimated eigenvectors under varying $n$.}
\label{tbl:eigenvec_case2}
\end{center}
\end{table}

Table \ref{tbl:eigenvec_case2} presents the average error (Err$_\xi$) for the estimated spiked eigenvectors across different methods and varying sample sizes $n$. Since the gIW and IW methods fail to accurately estimate the eigenvalues, as shown in Table \ref{tbl:eigenval_case2}, we exclude them from further consideration. At $n=20$, S-POET shows the best performance, whereas sample covariance performs the best for all other cases. Overall, sample covariance, gSIW, and S-POET exhibit comparable performance in estimating eigenvectors.

Combining the results from Case 1 and Case 2, we observe that IW and gIW generally exhibit poor performance across all settings. In contrast, gSIW shows improved performance in estimating spiked eigenvalues when the quantity $\dfrac{p}{n\lambda_i}$ is small. For spiked eigenvectors, sample covariance consistently demonstrates strong performance, and both gSIW and S-POET perform comparably to sample covariance.


\section{Real Data}\label{sec:real_data}
In this section, we apply the proposed gSIW prior to the MNIST dataset. To estimate the number of spiked eigenvalues in the covariance matrix, we select 50 images from the MNIST dataset, all labeled as the digit $7$. Based on the selected number of spiked eigenvalues $k$, we perform dimensionality reduction.

We assume that the MNIST dataset follows a multivariate normal distribution with a gSIW prior on the covariance matrix. The number of spiked eigenvalues, denoted as $k$, is estimated 
using two model selection criteria: the Watanabe-Akaike Information Criterion (WAIC) and $\mathrm{IC}_{p3}$. Dimensionality reduction is then performed based on the selected $k$. 
For comparison, we also include results based on the Growth Ratio (GR) method.

\begin{figure}[!h]
    \centering
    \includegraphics[scale=0.15]{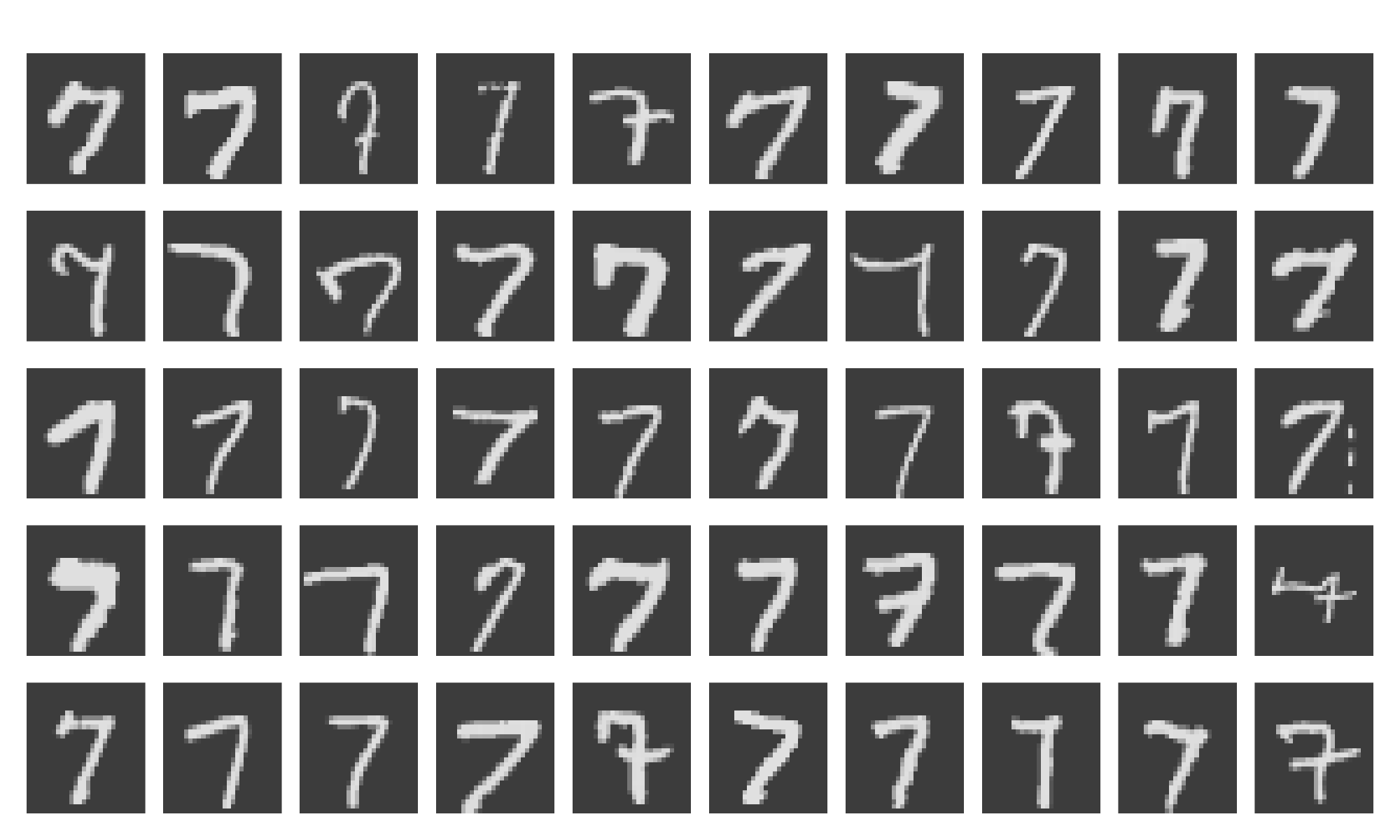}
    \caption{The images of selected $50$ MNIST samples labeled as $7$.\label{fig:real_data_origin}}
\end{figure} 

Figure \ref{fig:real_data_origin} displays the $50$ selected MNIST samples labeled as $7$. We flatten the $28 \times 28$ images into $784$-dimensional vectors. For each candidate $k$, the gSIW prior is applied, and the WAIC and $\mathrm{IC}_{p3}$ are computed. WAIC and $\mathrm{IC}_{p3}$ select $k = 7$ and $k = 9$, respectively, 
while the GR method selects $k = 1$ as the optimal number of spiked eigenvalues.

\begin{figure}[!h]
    \centering
    \begin{subfigure}[b]{0.45\textwidth}
        \centering
        \includegraphics[scale=0.12]{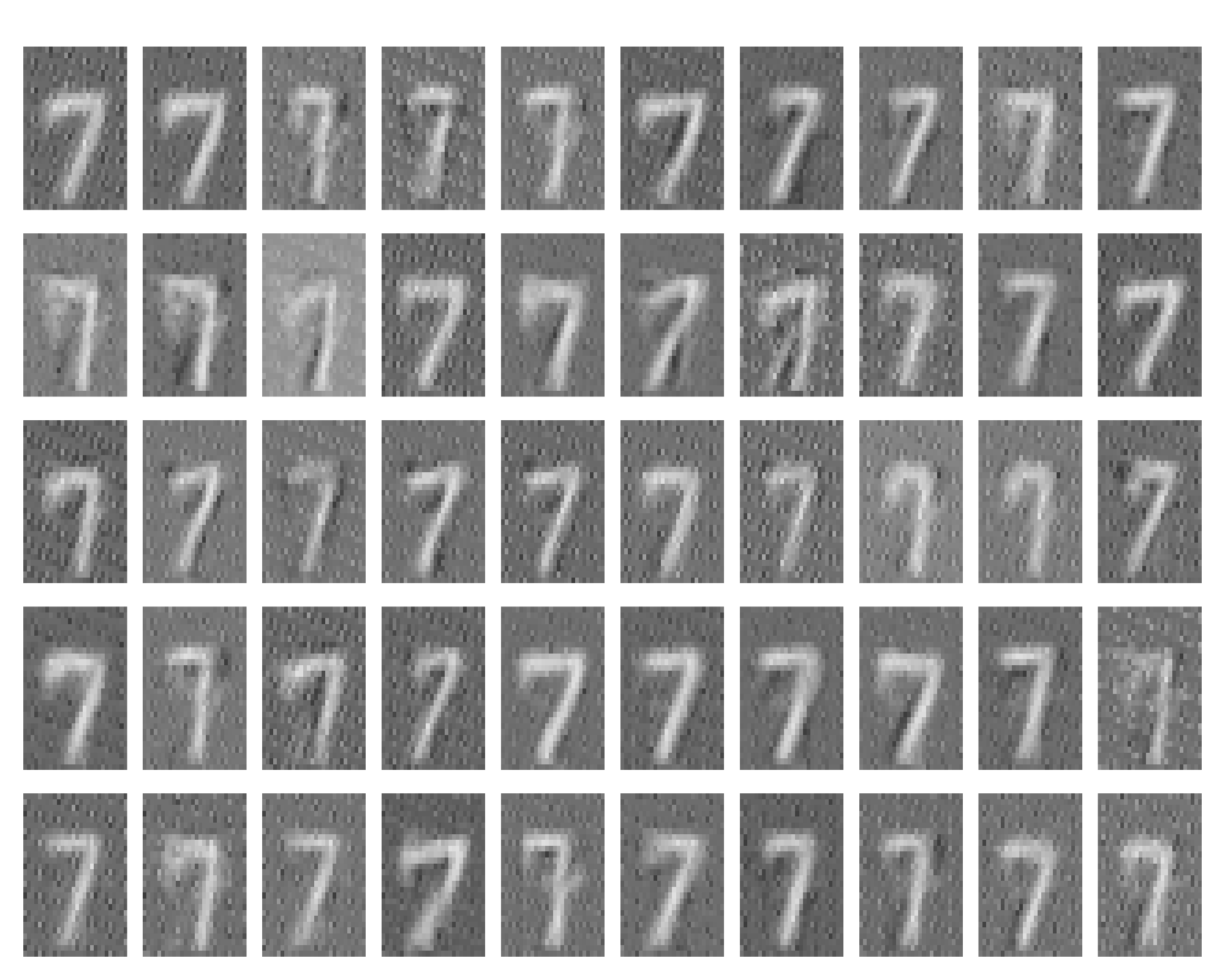}
        \caption{Dimensionality reduction with $k=7$ (WAIC).}
        \label{fig:mnist_waic}
    \end{subfigure}
    \hfill
    \begin{subfigure}[b]{0.45\textwidth}
        \centering
        \includegraphics[scale=0.12]{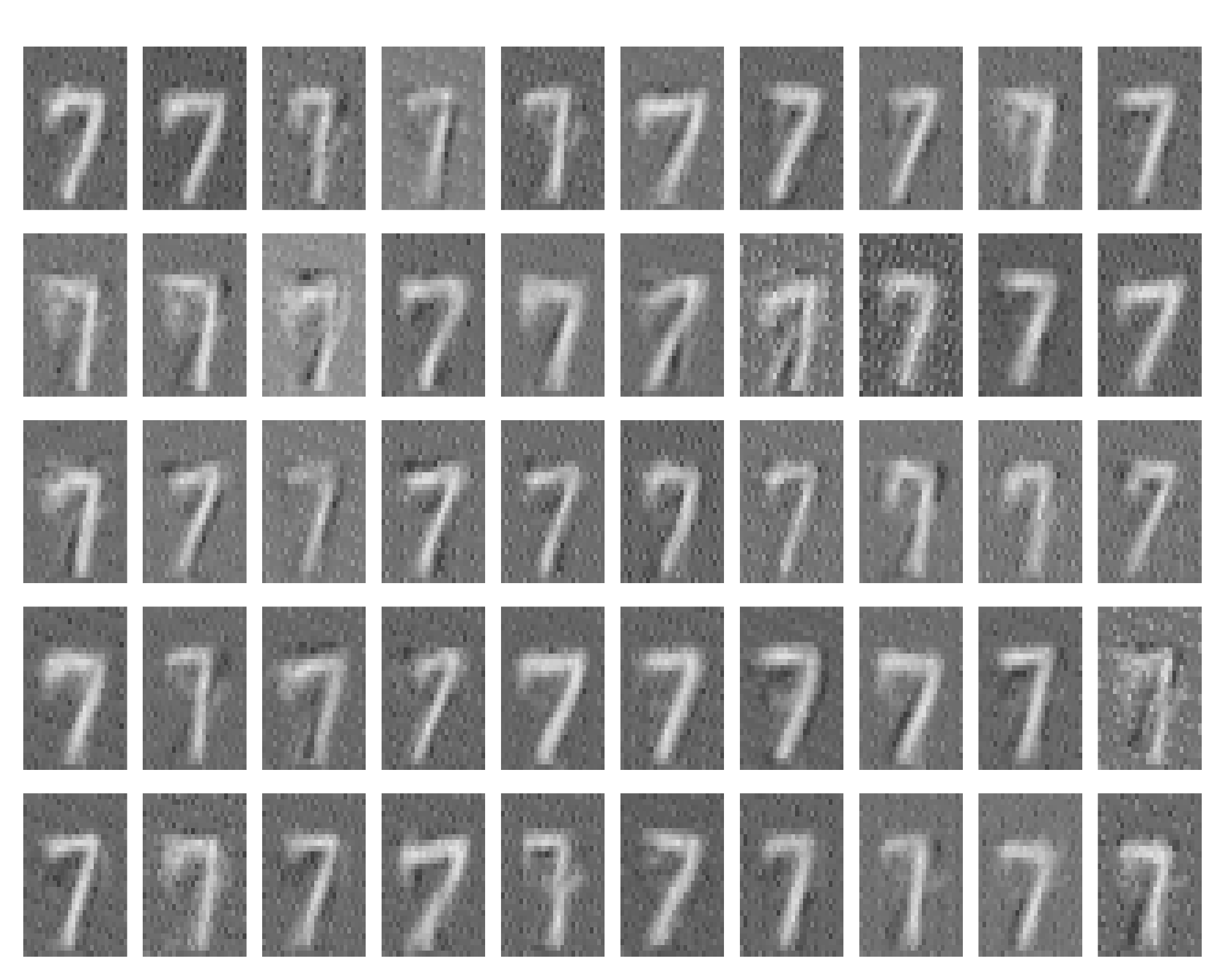}
        \caption{Dimensionality reduction with $k=9$ ($\mathrm{IC}_{p3}$).}
        \label{fig:mnist_ic}
    \end{subfigure}
    \caption{The images after dimensionality reduction of the selected $50$ MNIST samples labeled as $7$ 
    using different model selection criteria.}
    \label{fig:real_data_pca}
\end{figure}

Figure \ref{fig:real_data_pca} shows the images reconstructed using dimensionality reduction with $k=7$ and $k=9$, as determined by WAIC and $\mathrm{IC}_{p3}$, respectively. The reconstructed images successfully retain the shape of the digit $7$, indicating that both WAIC and $\mathrm{IC}_{p3}$ effectively capture the relevant subspace structure.

\begin{figure}[!ht]
    \centering
    \includegraphics[scale=0.15]{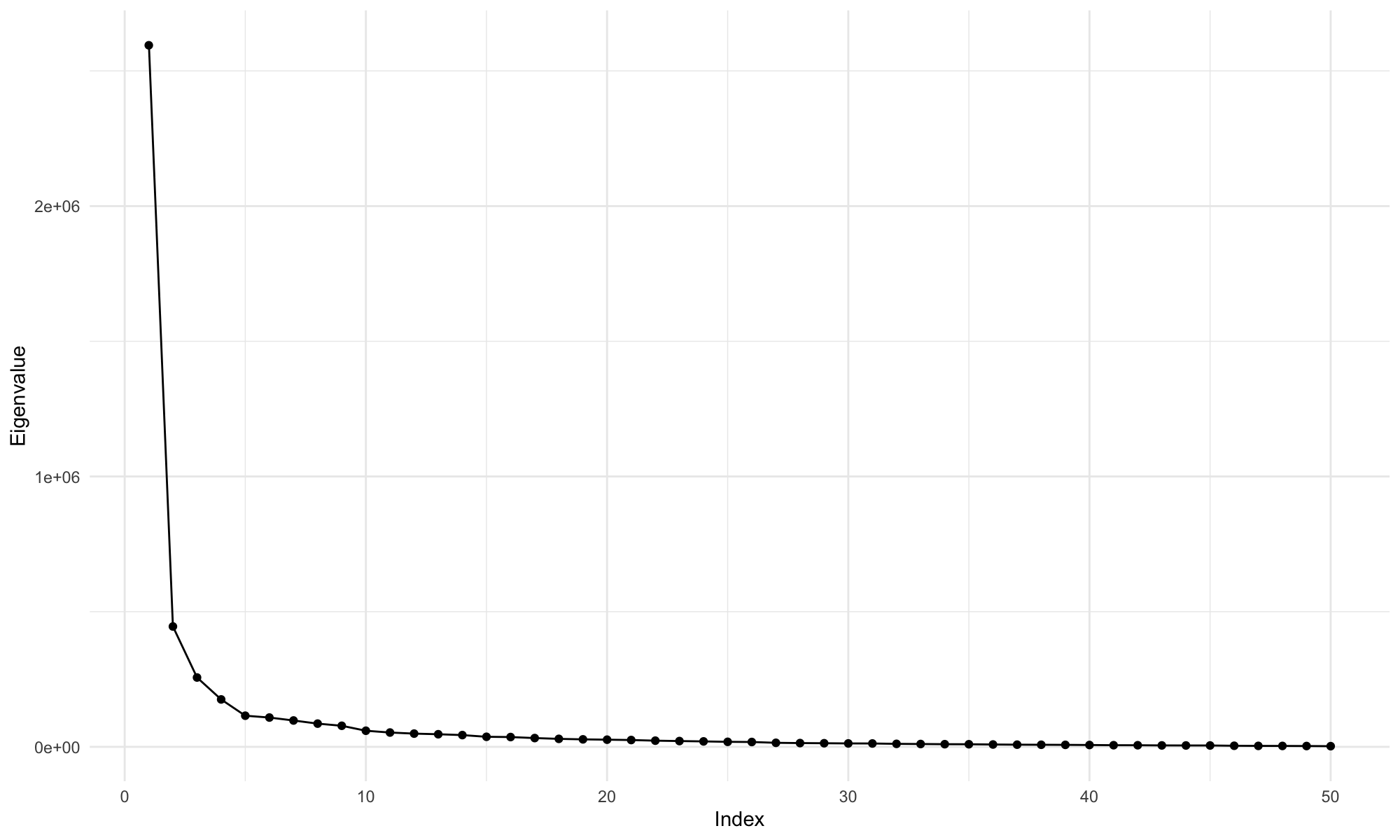}
    \caption{The first 50 eigenvalues of the sample covariance matrix.\label{fig:eigenval_mnist}}
\end{figure} 

In this case, where $(n, p) = (50, 784)$, the sample covariance eigenvalues are known to be heavily distorted due to high dimensionality. To illustrate this distortion, Figure \ref{fig:eigenval_mnist} shows the first $50$ eigenvalues of the sample covariance matrix. A few leading eigenvalues dominate the eigenvalue structure, while the remaining eigenvalues are much smaller. As a result, the GR method, which is based directly on sample eigenvalues, selects $k = 1$. By contrast, both WAIC and $\mathrm{IC}_{p3}$ select larger values of $k$ that lead to substantially better reconstruction.

To further assess the quality of the reduced representation, we compute the normalized mean squared reconstruction error (NMSE), defined as the mean squared error divided by the square of the range of the true values.  
The NMSE for $k = 7$ (WAIC) is 0.030, and for $k = 9$ ($\mathrm{IC}_{p3}$), it is 0.024, both substantially lower than the NMSE for $k = 1$ (GR), which is 0.069. This suggests that both WAIC and $\mathrm{IC}_{p3}$ identify subspaces that preserve the original variation well, while the GR-based choice performs poorly.

We also evaluate the cumulative explained variance ratio (CVE), defined as the proportion of total variance captured by the top $k$ eigenvalues of the sample covariance matrix. The CVE for $k = 7$ is 74.3\%, 
and for $k = 9$, it is 79.0\%, both substantially higher than the value for $k = 1$, which is 40.0\%.

These results demonstrate that, in this experiment, both WAIC and $\mathrm{IC}_{p3}$ provide reasonable estimates of $k$ that lead to effective dimensionality reduction and data reconstruction, while the GR method continues to underestimate $k$ due to the distortion of sample eigenvalues in high-dimensional settings.

\section{Conclusion}\label{sec:conclusion}

We proposed the generalized shrinkage inverse-Wishart (gSIW) prior,  
which generalizes the shrinkage inverse-Wishart (SIW) prior of \cite{berger2020bayesian}  
by allowing componentwise shrinkage and distinguishing between spiked and non-spiked eigenvalue structures.  
We established posterior expectation for both eigenvalues and eigenvectors under the gSIW prior, and showed that it outperforms existing priors both theoretically and empirically when Assumptions A1 through A5 are satisfied.  

Our theoretical analysis builds on the asymptotic behavior of sample covariance matrices developed in \cite{wang2017asymptotics}, particularly under the spiked covariance model. Future work may extend our results beyond the current assumptions, such as relaxing the spiked structure or improving asymptotic bounds for sample eigenstructures in more general settings.

\section*{Acknowledgements}
This work was supported by the National Research Foundation of Korea (NRF) grant funded by the Korea government(MSIT) (No. NRF-2023R1A2C1003050).

\bibliographystyle{dcu}
\bibliography{SIW}

\end{document}


\maketitle

\newpage
\section{Appendix}\label{sec:appendix}
Consider the transformed data $Y_i=\Gamma^T X_i$ which follows $N(0,\Lambda)$ where $\Lambda=diag(\lambda_{0,1},\cdots,\lambda_{0,p})$ and $\lambda_{0,i}$ is $i$th eigenvalue of true covariance. Since the eigenvalues of sample covariance are invariant under the orthogonal transform, we consider the sample covariance $S=\dfrac{1}{n}\boldsymbol{Y}^T\boldsymbol{Y}$ where $\boldsymbol{Y}=(Y_1,\ldots, Y_n)^T$.
Define $\tilde{S}= \dfrac{1}{n}\boldsymbol{Y}\boldsymbol{Y}^T $, which preserves the first $n$ eigenvalues of $S$. Let $\tilde{Y}_i$ is $i$th column of $\boldsymbol{Y}$, so that $\tilde{Y_i}\overset{ind}{\sim} N(0,\lambda_{0,i}I_n)$. Therefore, $\tilde{S}=\dfrac{1}{n} \sum\limits_{i=1}^p \tilde{Y}_i \tilde{Y}_i^T$. It can be expressed as 
\begin{align*}
    \tilde{S}=\dfrac{1}{n}\sum\limits_{i=1}^p\lambda_{0,i} Z_i Z_i^T,
\end{align*}
where $Z_1,\cdots,Z_p$ are independent $n$-dimensional vectors, and each elements of $Z_i$ follows an independent standard normal distribution.

\subsection*{Asymptotic properties of eigenstructure of sample covariance}
\begin{lemma}[Lemma A.1 of \citealp{wang2017asymptotics}]\label{lem:non-spiked_eig}
Suppose that $A_1,\ldots,A_p$ are independent $n$-dimensional Gaussian random vectors with mean $0$ and variance $I_n$. For all $t\geq 0$, the inequality holds at least $1-2\exp(ct^2)$ probability for some positive constant $c$
    $$\bar{w}-\max(\delta,\delta^2)\leq \lambda_n\Big(\dfrac{1}{p}\sum\limits_{i=1}^p w_iA_iA_i^T\Big)\leq \lambda_1\Big(\dfrac{1}{p}\sum\limits_{i=1}^p w_iA_iA_i^T\Big)\leq \bar{w}+\max(\delta,\delta^2),$$
where $\delta = C\sqrt{n/p}+t/\sqrt{p}$, $C$ is positive constant, $\abs{w_i}$ is bounded for all $i$, and $\bar{w}= p^{-1}\sum\limits_{i=1}^p w_i$.
\end{lemma}

The first $n$ eigenvalues of the sample covariance is equivalence with that of $\tilde{S}$. The $\tilde{S}$ can be decomposed into the sum of two matrices $A$ and $B$, where $A=\dfrac{1}{n}\sum\limits_{i=1}^k\lambda_{0,i} Z_i Z_i^T$ and $B=\dfrac{1}{n}\sum\limits_{i=k+1}^p\lambda_{0,i} Z_i Z_i^T$.

\begin{lemma}[Asymptotic properties of eigenvalues of sample covariance]\label{lem:eig}
Under model \eqref{main-eq:model}, the eigenvalues of sample covariance satisfy the following properties for all sufficiently large $n$,
$$\dfrac{\hat{\lambda}_j}{\lambda_{0,j}} =
\begin{cases}
    1+\bar{d}d_j +\alpha_j{\lambda_{0,j}}^{-1}\sqrt{\dfrac{p}{n}}+\beta_j,\quad & j=1,\ldots,k\\
    \bar{d}d_j +\alpha_j{\lambda_{0,j}}^{-1}\sqrt{\dfrac{p}{n}},\quad & j=k+1,\ldots,n,
\end{cases}$$
where $\alpha_j$ is a constant in the interval $[-C,C]$ for some positive constant $C$, and $\beta_j\lesssim n^{-1/2+\delta}$ for all small $\delta>0$.

\end{lemma}
\begin{proof}
According to Lemma \ref{lem:non-spiked_eig}, for \( t = \sqrt{n} \), the inequality holds with probability at least $1 - 2\exp(-cn)$:
\begin{align}\label{ineq:bound_lambda_B}
     \dfrac{1}{p-k} \sum\limits_{i=k+1}^p \lambda_{0,i} - C \sqrt{\dfrac{n}{p}}\leq \lambda_l\left( \dfrac{n}{p-k} B \right)\leq \dfrac{1}{p-k} \sum\limits_{i=k+1}^p \lambda_{0,i} + C \sqrt{\dfrac{n}{p}} ,
\end{align}
for some positive constant $C$ and for all sufficiently large $n$.

The inequality follows from Wely's theorem:
$$ \dfrac{\lambda_j(A)}{\lambda_{0,j}}+\dfrac{\lambda_p(B)}{\lambda_{0,j}}\leq\dfrac{\hat{\lambda}_j}{\lambda_{0,j}} = \dfrac{{\lambda}_j(A+B)}{\lambda_{0,j}}\leq \dfrac{\lambda_j(A)}{\lambda_{0,j}}+\dfrac{\lambda_1(B)}{\lambda_{0,j}},$$
for $j=1,\cdots,n$. Given the bound in \eqref{ineq:bound_lambda_B}, we obtain the following inequality
$$\dfrac{\lambda_j(A)}{\lambda_{0,j}} + \bar{d} d_j - C {\lambda_{0,j}}^{-1} \sqrt{\dfrac{p}{n}} \leq \dfrac{\hat{\lambda}_j}{\lambda_{0,j}} \leq \dfrac{\lambda_j(A)}{\lambda_{0,j}} + \bar{d} d_j + C {\lambda_{0,j}}^{-1} \sqrt{\dfrac{p}{n}},$$
where $\bar{d} = \dfrac{1}{p-k} \sum\limits_{i=k+1}^p \lambda_{0,i}$, and $d_j = \dfrac{p}{n\lambda_{0,i}}$.  Therefore, we obtain the following equality
$$ \dfrac{\hat{\lambda}_j}{\lambda_{0,j}} = \dfrac{\lambda_j(A)}{\lambda_{0,j}} + \bar{d} d_j + \alpha_j {\lambda_{0,j}}^{-1} \sqrt{\dfrac{p}{n}}, \quad j = 1, \ldots, n,$$
for some constant $\alpha_j\in[-C,C]$.

According to Lemma A.2 of \cite{wang2017asymptotics}, the following asymptotic normality and independence hold
$$\sqrt{n} \left( \dfrac{\lambda_j(A)}{\lambda_{0,j}} - 1 \right) \overset{d}{\rightarrow} N(0,2),\quad \forall i=1,\ldots,k.$$

We obtain the following asymptotic normality
$$ \sqrt{n} \left( \dfrac{\hat{\lambda}_j}{\lambda_{0,j}} - \left( 1 + \bar{d} d_j + \alpha_j {\lambda_{0,j}}^{-1} \sqrt{\dfrac{p}{n}} \right) \right) \overset{d}{\rightarrow} N(0,2), \quad j = 1, \ldots, k. $$

Since $k$ is constant, for $\delta_n \prec \sqrt{n}$, the following inequality holds for sufficiently large $n$ and $p$
$$ \delta_n \left| \dfrac{\hat{\lambda}_j}{\lambda_{0,j}} - \left( 1 + \bar{d} d_j + \alpha_j {\lambda_{0,j}}^{-1} \sqrt{\dfrac{p}{n}} \right) \right| < \epsilon, \quad \forall j = 1, \ldots, k. $$

Therefore, we obtain the following equality
\begin{align*}
    \dfrac{\hat{\lambda}_j}{\lambda_{0,j}} &= 1 + \beta_j + \bar{d} d_j + \alpha_j {\lambda_{0,j}}^{-1} \sqrt{\dfrac{p}{n}}, \quad  j=1,\ldots,k,
\end{align*}
for $\beta_j\lesssim n^{-1/2+\delta}$ for all small $\delta>0$.

For $j>k$, $\lambda_j(A) = 0$, and it implies
$$\dfrac{\hat{\lambda}_j}{\lambda_{0,j}} = \ \bar{d} d_j + \alpha_j {\lambda_{0,j}}^{-1} \sqrt{\dfrac{p}{n}}.$$

\end{proof}

Let $(S,d)$ be a metric space and $K\subseteq S$. For $\epsilon>0$, the $\epsilon$-covering number of $K$ with respect to the metric $d$, denoted by $N(K,d,\epsilon)$, is the minimum number of closed balls of radius $\epsilon$ (with respect to $d$) required to cover the set $K$. Similarly, the $\epsilon$-packing number of $K$ with respect to the metric $d$, denoted by $M(K,d,\epsilon)$, is the maximum number of points in $K$ such that any two distinct points are at least distance $\epsilon$ apart.

\begin{lemma}[Covering number of Orthogonal group ; Proposition 7 of \citealt{szarek1982nets}]\label{lem:Orthogonal_covering}
Let $\nu$ be a unitary ideal norm, for which $\nu(PAQ)=\nu(A)$ holds for any unitary matrices $P$ and $Q$, and any matrix $A$. Then, there exist universal positive constants $c$ and $C$ that satisfy the following inequality for $\epsilon\in(0,2\nu(I)]$
    $$(c\nu(I)/\epsilon)^d\leq N(O(m),\nu,\epsilon)\leq (C\nu(I)/\epsilon)^d,$$
where $d=m(m-1)/2$ and $N(O(m),\nu,\epsilon)$ is $\epsilon$-covering number of $O(m)$ with respect to the metric $\rho(x,y) = \nu(x-y)$, where $x,y\in O(m)$.    
\end{lemma}

Grassmann manifold $G_{n,p}$ denotes the $n$-dimensional subspaces of $\bbR^p$. It can be expressed as a quotient space $G_{n,p} = O(p)/(O(n)\times O(p-n))$. \cite{edelman1998geometry} states that a point in the Grassmann manifold can be represented as an equivalence class
$$[Q] = \{Q \begin{pmatrix}
        Q_1 & 0 \\ 
        0& Q_2
    \end{pmatrix} : Q_1\in O(n), Q_2\in O(p-n)\},$$
which consists of all orthogonal matrices that span the same subspace as the first $n$ columns of $Q$. \cite{szarek1982nets} defines the quotient metric on $G_{n,p}$ by the unitary ideal norm $\nu$ on $O(p)$ as follows
$$\rho_{\nu}(H_1,H_2)= \inf\{\nu(I-V) : V\in O(p),VH_1=H_2\},$$
for $H_1,H_2\in G_{n,p}$. Therefore, $H_1$ and $H_2$ can be represented as follows
\begin{align*}
    H_1 &= \{P_1 \begin{pmatrix}
        Q_1 & 0 \\ 
        0& Q_2
    \end{pmatrix} : Q_1\in O(n), Q_2\in O(p-n)\},\\
    H_2 &= \{P_2 \begin{pmatrix}
        Q_1 & 0 \\ 
        0& Q_2
    \end{pmatrix} : Q_1\in O(n), Q_2\in O(p-n)\},
\end{align*}
for some $P_1,P_2\in O(p)$. The quotient metric is then given by
\begin{align*}
    &\rho_{\nu}(H_1,H_2)\\
    &=\inf\Bigg\{\nu(I-V): V = P_2\begin{pmatrix}
        Q_1 & 0 \\ 
        0& Q_2
    \end{pmatrix}\begin{pmatrix}
        Q_3 & 0 \\ 
        0& Q_4
    \end{pmatrix} P_1^T , \;Q_1,Q_3\in O(n),\;Q_2,Q_4\in O(p-n)\Bigg\}\\
    &=\inf_{Q_1\in O(n),Q_2\in O(p-n)}\nu\Big(I-P_2\begin{pmatrix}
        Q_1 & 0 \\ 
        0& Q_2
    \end{pmatrix} P_1^T\Big)\\
    &= \inf_{Q_1\in O(n),Q_2\in O(p-n)}\nu\Big(P_1-P_2\begin{pmatrix}
        Q_1 & 0 \\ 
        0& Q_2
    \end{pmatrix}\Big).
\end{align*}

\begin{lemma}[Covering number of Grassmann manifold ; Proposition 8 of \citealt{szarek1982nets}]\label{lem:Grassmann_covering}
For all unitary ideal norms $\nu$, there exist universal positive constants $m$ and $M$ that satisfy the following inequality for $\epsilon\in(0,D_\nu]$
    $$(m D_\nu/\epsilon)^d\leq N(G_{n,p},\rho_\nu,\epsilon)\leq (M D_\nu/\epsilon)^d,$$
where $d=n(p-n)=dim G_{n,p}$, and $D_\nu$ is the diameter of $G_{n,p}$ with respect to the metric $\rho_\nu$. Specifically, for the Frobenius norm $\nu$, the diameter is given by $D_\nu = 2\sqrt{\min(n,p-n)}$.

\end{lemma}

\begin{lemma}[Probability of subset on orthogonal group]\label{lem:subset_prob1}
Consider the following subset of the orthogonal group $O(p)$, defined as
\begin{align*}
    B_{\epsilon} = \bigg\{\Gamma\in O(p):\displaystyle\inf_{Q_1\in O(n),\;Q_2\in O(p-n)} \abs{\abs{\begin{pmatrix}
       Q_1 & 0 \\ 
        0& Q_2
    \end{pmatrix}-\Gamma}}_F <\epsilon\bigg\}.
\end{align*}
Then, the probability measure of $B_\epsilon$ satisfies the following lower bound: 
\begin{align*}
    \bbP(B_\epsilon) := \displaystyle\int_{B_\epsilon}[d\Gamma] &\geq \Big(\dfrac{2M\sqrt{\min(n,p-n)}}{\epsilon}\Big)^{-n(p-n)},
\end{align*} 
for some positive constant $M$.
\end{lemma}
\begin{proof}
Let $S_1,\ldots,S_{N(G_{n,p},\rho_{\abs{\abs{\cdot}}_F},\epsilon)}$ form an $\epsilon$-covering of $G_{n,p}$ with respect to the distance $\rho_{\abs{\abs{\cdot}}_F}$. Then, the inequality holds
\begin{align*}
    1 &= \bbP\Big(\bigcup_{i=1}^{N(G_{n,p},\rho_{\abs{\abs{\cdot}}_F},\epsilon)}\{\Gamma \in O(p): \rho_{\abs{\abs{\cdot}}_F}([\Gamma],S_i)<\epsilon \}\Big)\\
    &\leq \sum\limits_{i=1}^{N(G_{n,p},\rho_{\abs{\abs{\cdot}}_F},\epsilon)}\bbP(\{\Gamma \in O(p): \rho_{\abs{\abs{\cdot}}_F}([\Gamma],S_i)<\epsilon \})\\
    &= N(G_{n,p},\rho_{\abs{\abs{\cdot}}_F},\epsilon) \cdot \bbP(\{\Gamma \in O(p): \rho_{\abs{\abs{\cdot}}_F}([\Gamma],S_i)<\epsilon \})\\
    &= N(G_{n,p},\rho_{\abs{\abs{\cdot}}_F},\epsilon) \cdot \bbP(\{\Gamma \in O(p): \rho_{\abs{\abs{\cdot}}_F}([\Gamma],[I_p])<\epsilon \})
\end{align*}
where $[\Gamma]:=\{\Gamma \begin{pmatrix}
        Q_1 & 0\\
        0 & Q_2
\end{pmatrix} : Q_1\in O(n),\;Q_2\in O(p-n) \}$ for $ \Gamma\in O(p)$.

We obtain the following lower bound on the probability measure of $B_\epsilon$
\begin{align*}
    \bbP(B_{\epsilon}) &= \bbP\bigg(\bigg\{\Gamma\in O(p):\displaystyle\inf_{Q_1\in O(n),\;Q_2\in O(p-n)} \abs{\abs{\begin{pmatrix}
       Q_1 & 0 \\ 
        0& Q_2
    \end{pmatrix}-\Gamma}}_F <\epsilon\bigg\}\bigg)\nonumber\\
     &= \bbP(\{\Gamma \in O(p): \rho_{\abs{\abs{\cdot}}_F}([\Gamma],[I_p])<\epsilon
    \})\\
    &\geq \dfrac{1}{N(G_{n,p},\rho_{\abs{\abs{\cdot}}_F},\epsilon)}\\
    &\geq \Big(\dfrac{2M\sqrt{\min(n,p-n)}}{\epsilon}\Big)^{-n(p-n)},
\end{align*}
where $M$ is a positive constant. The last inequality follows from Lemma \ref{lem:Grassmann_covering}.
\end{proof}

\begin{lemma}[Probability of subset on orthogonal group]\label{lem:subset_prob2}
Consider the following subset of the orthogonal group $O(p)$, defined as
\begin{align*}
    C_{\epsilon}  &= \bigg\{\Gamma\in O(p):\displaystyle\inf_{Q_2\in O(p-n)} \abs{\abs{\begin{pmatrix}
        I_n & 0 \\ 
        0& Q_2
    \end{pmatrix}-\Gamma}}_F <\epsilon
    \bigg\}.
\end{align*}
Then, the probability measure of $C_\epsilon$ satisfies the following lower bound: 
\begin{align*}
    \bbP(C_\epsilon) := \displaystyle\int_{C_\epsilon}[d\Gamma] &\geq \Big(\dfrac{c\sqrt{n}}{\epsilon}\Big)^{-n(p-\frac{n}{2}-\frac{1}{2})},\quad \epsilon\leq 2\sqrt{n},
\end{align*} 
for some positive constant $c$.

\end{lemma}
\begin{proof}

Let $S_1,\cdots,S_{N(O(n),\abs{\abs{\cdot}}_F,\epsilon_1)}$ form an $\epsilon_1$-covering of $O_{n}$ with respect to the distance $\abs{\abs{\cdot}}_F$ for $\epsilon_1\in(0,\epsilon)$. For all $Q_1\in O(n)$, there exists some $S_i$ such that $\abs{\abs{S_i - Q_1}}_F<\epsilon_1$. By triangular inequality, we obtain 
\begin{align}\label{ineq:tri}
     \abs{\abs{\begin{pmatrix}
        S_i & 0 \\ 
       0& Q_2
    \end{pmatrix}-\Gamma}}_F \leq \abs{\abs{\begin{pmatrix}
        S_i & 0 \\ 
       0& Q_2
    \end{pmatrix}-\begin{pmatrix}
       Q_1 & 0 \\ 
        0& Q_2
    \end{pmatrix}}}_F + \abs{\abs{\begin{pmatrix}
       Q_1 & 0 \\ 
        0& Q_2
    \end{pmatrix}-\Gamma}}_F.
\end{align}
Taking the infimum over $Q_2\in O(p-n)$ on both sides of \eqref{ineq:tri}, then the inequality holds
\begin{align}\label{ineq:inf_tri}
    \inf_{Q_2\in O(p-n)}\abs{\abs{\begin{pmatrix}
        S_i & 0 \\ 
       0& Q_2
    \end{pmatrix}-\Gamma}}_F \leq \inf_{Q_2\in O(p-n)}\abs{\abs{\begin{pmatrix}
        S_i & 0 \\ 
       0& Q_2
    \end{pmatrix}-\begin{pmatrix}
       Q_1 & 0 \\ 
        0& Q_2
    \end{pmatrix}}}_F + \inf_{Q_2\in O(p-n)}\abs{\abs{\begin{pmatrix}
       Q_1 & 0 \\ 
        0& Q_2
    \end{pmatrix}-\Gamma}}_F.
\end{align}

We obtain the following lower bound on the probability measure of $C_\epsilon$:
\begin{align*}
    \bbP(C_{\epsilon}) &=  \bbP\bigg(\bigg\{\Gamma \in O(p):\displaystyle\inf_{Q_2\in O(p-n)} \abs{\abs{\begin{pmatrix}
        I_n & 0 \\ 
       0& Q_2
    \end{pmatrix}-\Gamma}}_F <\epsilon
    \bigg\}\bigg)\nonumber\\
    &\geq \dfrac{1}{N(O(n),\abs{\abs{\cdot}}_F,\epsilon_1)} \bbP\bigg(\bigcup_{i=1}^{N(O(n),\abs{\abs{\cdot}}_F,\epsilon_1)}\bigg\{\Gamma \in O(p):\displaystyle\inf_{Q_2\in O(p-n)} \abs{\abs{\begin{pmatrix}
        S_i & 0 \\ 
       0& Q_2
    \end{pmatrix}-\Gamma}}_F <\epsilon
    \bigg\}\bigg)\nonumber\\
    &\geq\dfrac{1}{N(O(n),\abs{\abs{\cdot}}_F,\epsilon_1)} \bbP\bigg(\bigg\{\Gamma\in O(p):\displaystyle\inf_{Q_1\in O(n),\;Q_2\in O(p-n)} \abs{\abs{\begin{pmatrix}
       Q_1 & 0 \\ 
        0& Q_2
    \end{pmatrix}-\Gamma}}_F <\epsilon-\epsilon_1\bigg\}\bigg)\\
    &= \dfrac{ \bbP(A_{\epsilon-\epsilon_1})}{N(O(n),\abs{\abs{\cdot}}_F,\epsilon_1)}.
\end{align*} 
The first inequality holds by the following fact
$$\bbP(C_{\epsilon}) = \bbP\bigg(\bigg\{\Gamma \in O(p):\displaystyle\inf_{Q_2\in O(p-n)} \abs{\abs{\begin{pmatrix}
        S_i & 0 \\ 
       0& Q_2
    \end{pmatrix}-\Gamma}}_F <\epsilon
    \bigg\}\bigg),$$
for $i = 1,\ldots, N(O(n),\abs{\abs{\cdot}}_F,\epsilon_1)$. The second inequality holds by \eqref{ineq:inf_tri}.

Therefore, by applying Lemma \ref{lem:Orthogonal_covering} and Lemma \ref{lem:Grassmann_covering}, we obtain the following inequality:
\begin{align*}
    \bbP(C_{\epsilon}) &\geq \Big(\dfrac{2M\sqrt{\min(n,p-n)}}{\epsilon-\epsilon_1}\Big)^{-n(p-n)} \cdot \Big(\dfrac{C\sqrt{n}}{\epsilon_1}\Big)^{-n(n-1)/2}\\
    &\geq \Big(\dfrac{2M\sqrt{n}}{\epsilon-\epsilon_1}\Big)^{-n(p-n)} \cdot \Big(\dfrac{C\sqrt{n}}{\epsilon_1}\Big)^{-n(n-1)/2}\\
    &= \Big(\dfrac{4M\sqrt{n}}{\epsilon}\Big)^{-n(p-n)} \cdot \Big(\dfrac{2C\sqrt{n}}{\epsilon}\Big)^{-n(n-1)/2}\\
    &\geq \Big(\dfrac{\max(4M,2C)\sqrt{n}}{\epsilon}\Big)^{-n(p-\frac{n}{2}-\frac{1}{2})},
\end{align*}
where we set $\epsilon_1 =\dfrac{\epsilon}{2}$, and $M$ and $C$ are positive constants. 

As a result, for all sufficiently large $n$, we obtain the lower bound
\begin{align*}
    \bbP(C_\epsilon)\geq \Big(\dfrac{c\sqrt{n}}{\epsilon}\Big)^{-n(p-\frac{n}{2}-\frac{1}{2})},
\end{align*}
for some positive constant $c$.

\end{proof}


We define $c_i$ as $(i,i)$ elements of $\Gamma^T(hI_p + W)\Gamma$, where 
$$W = diag(n\hat{\lambda_1},\ldots,n\hat{\lambda_n},0,\ldots,0),$$
and $\hat{\lambda}_i$ is the $i$th eigenvalue of sample covariance $S$.
By Lemma \ref{lem:eig}, we express $c_i$ as follows:
\begin{align*}
    \dfrac{c_i}{n} &= \dfrac{h}{n} + \sum\limits_{j=1}^n \Gamma_{ji}^2\hat{\lambda}_j\\
    &= \dfrac{h}{n} + \sum\limits_{j=1}^k\Gamma_{ji}^2\Big\{(1+\beta_j)\lambda_{0,j}+(\bar{d}\dfrac{p}{n}+\alpha_j\sqrt{\dfrac{p}{n}})\Big\}+ \sum\limits_{j=k+1}^n\Gamma_{ji}^2\Big(\bar{d}\dfrac{p}{n}+\alpha_j\sqrt{\dfrac{p}{n}}\Big).
\end{align*}
Next, we establish two lemmas that provide upper and lower bounds for $c_i$.

\begin{lemma}[Upper bound of $c_i$]\label{lem:upper_bound_c}
Under the conditions $A1-A4$ in Section \ref{main-sec:main_results}, the quantity $c_i$ satisfies the following upper bound on $C_\epsilon$
$$\dfrac{c_i}{n} \leq
\begin{cases}
\bigg((1+\beta_i)\lambda_{0,i}  +\Big( \bar{d}\dfrac{p}{n}+\alpha_i\sqrt{\dfrac{p}{n}}\Big)+\dfrac{h}{n}\bigg) (1 + 4\epsilon^2\dfrac{\lambda_{0,1}}{\lambda_{0,k}}) \quad& i=1,\ldots,k,\\ 
\Big(\bar{d}\dfrac{p}{n}+\alpha_i\sqrt{\dfrac{p}{n}}+\dfrac{h}{n}\Big)\Big(1+4\epsilon^2\dfrac{n\lambda_{0,1}}{\bar{d}p}\Big) \quad& i=k+1,\ldots,n,\end{cases}$$
for all sufficiently large $n$. Furthermore, the following inequality also holds on $C_\epsilon$
$$\prod\limits_{i=n+1}^p\dfrac{c_i}{n} \leq \Big( \dfrac{h}{n}\Big)^{p-n} \bigg[1+\dfrac{2\epsilon^2n}{(p-n)h}\lambda_{0,1}  \bigg]^{p-n},$$
where $\beta_{\max} = \max_{1\leq j\leq k}\beta_j$, and for all sufficiently large $n$.
\end{lemma}

\begin{proof}
For $\Gamma\in C_\epsilon$, the following inequality holds:
$$(\Gamma_{ii}-1)^2 + \sum\limits_{\substack{j=1 \\ j \neq i}}^p\Gamma_{ji}^2<\epsilon^2,$$
for $i=1,\ldots,n$. Since $\sum\limits_{j=1}^p\Gamma_{ji}^2 =1$, it follows that
$$(1-\Gamma_{ii}^2)<\epsilon^2,$$
which implies $\Gamma_{ii}^2>1-\epsilon^2,$ for $i=1,\ldots,n$. 

For $i=1,\ldots,k$, we have the following bound:
\begin{align*}
    \dfrac{c_i}{n} &= \dfrac{h}{n}+\sum\limits_{j=1}^k (1+\beta_j)\Gamma_{ji}^2 \lambda_{0,j} + \sum\limits_{j=1}^n\Big(\bar{d}\dfrac{p}{n}+\alpha_j\sqrt{\dfrac{p}{n}}\Big)\Gamma_{ji}^2\\
    &\leq \dfrac{h}{n} + (1+\beta_i)\lambda_{0,i}\Gamma_{ii}^2   + \Big((1+\beta_{\max})\lambda_{0,1}\Big)\sum\limits_{\substack{j=1 \\ j \neq i}}^p\Gamma_{ji}^2+\Big( \bar{d}\dfrac{p}{n}+\alpha_i\sqrt{\dfrac{p}{n}}\Big)\Gamma_{ii}^2+\Big( \bar{d}\dfrac{p}{n}+C\sqrt{\dfrac{p}{n}}\Big)\sum\limits_{\substack{j=1 \\ j \neq i}}^p\Gamma_{ji}^2\\
    &\leq \dfrac{h}{n} + (1+\beta_i)\lambda_{0,i}  +\Big( \bar{d}\dfrac{p}{n}+\alpha_i\sqrt{\dfrac{p}{n}}\Big) + \Big((1+\beta_{\max})\lambda_{0,1} + \bar{d}\dfrac{p}{n}+C\sqrt{\dfrac{p}{n}} \Big)(1-\Gamma_{ii}^2)\\
    &\leq \dfrac{h}{n} + (1+\beta_i)\lambda_{0,i}  +\Big( \bar{d}\dfrac{p}{n}+\alpha_i\sqrt{\dfrac{p}{n}}\Big) + 2\epsilon^2\lambda_{0,1}\\
    &\leq \bigg((1+\beta_i)\lambda_{0,i}  +\Big( \bar{d}\dfrac{p}{n}+\alpha_i\sqrt{\dfrac{p}{n}}\Big)+\dfrac{h}{n}\bigg) (1 + 4\epsilon^2\dfrac{\lambda_{0,1}}{\lambda_{0,k}}),
\end{align*}
for all sufficiently large $n$. The third inequality follows from the fact that:
$$\lambda_{0,1}\succ \beta_{\max}\lambda_{0,1} +  \bar{d}\dfrac{p}{n}+C\sqrt{\dfrac{p}{n}}.$$
The fourth inequality follows from
$$(1+\beta_i)\lambda_{0,i}  +\Big( \bar{d}\dfrac{p}{n}+\alpha_i\sqrt{\dfrac{p}{n}}\Big)+\dfrac{h}{n}\geq\dfrac{\lambda_{0,k}}{2},$$
for all sufficiently large $n$.

For $i=k+1,\ldots,n$, we have the following bound:
\begin{align*}
    \dfrac{c_i}{n} &=\dfrac{h}{n}+\sum\limits_{j=1}^k (1+\beta_j)\Gamma_{ji}^2 \lambda_{0,j} + \sum\limits_{j=1}^n\Big(\bar{d}\dfrac{p}{n}+\alpha_j\sqrt{\dfrac{p}{n}}\Big)\Gamma_{ji}^2\\
    &\leq \dfrac{h}{n} + (1+\beta_{\max})\lambda_{0,1}\sum\limits_{\substack{j=1 \\ j \neq i}}^p\Gamma_{ji}^2 +\Big(\bar{d}\dfrac{p}{n}+\alpha_i\sqrt{\dfrac{p}{n}}\Big)\Gamma_{ii}^2 +\Big(\bar{d}\dfrac{p}{n}+C\sqrt{\dfrac{p}{n}}\Big)\sum\limits_{\substack{j=1 \\ j \neq i}}^p\Gamma_{ji}^2\\
    &\leq \dfrac{h}{n}+\Big(\bar{d}\dfrac{p}{n}+\alpha_i\sqrt{\dfrac{p}{n}}\Big) + \Big( (1+\beta_{\max})\lambda_{0,1} +\bar{d}\dfrac{p}{n}+C\sqrt{\dfrac{p}{n}}\Big)(1-\Gamma_{ii}^2)\\
    &\leq \dfrac{h}{n}+\Big(\bar{d}\dfrac{p}{n}+\alpha_i\sqrt{\dfrac{p}{n}}\Big) + 2\lambda_{0,1}\epsilon^2\\
    &\leq \Big(\bar{d}\dfrac{p}{n}+\alpha_i\sqrt{\dfrac{p}{n}}+\dfrac{h}{n}\Big)\Big(1+4\epsilon^2\dfrac{n\lambda_{0,1}}{\bar{d}p}\Big),
\end{align*}
for all sufficiently large $n$. The third inequality follows from the fact that
$$\lambda_{0,1}\succ \beta_{\max}\lambda_{0,1} +  \bar{d}\dfrac{p}{n}+C\sqrt{\dfrac{p}{n}}.$$
The fourth inequality follows from
$$\bar{d}\dfrac{p}{n}+\alpha_i\sqrt{\dfrac{p}{n}}+\dfrac{h}{n}\geq \dfrac{1}{2}\dfrac{\bar{d}p}{n},$$
for all sufficiently large $n$.

For $\Gamma\in C_\epsilon$, the following bound holds for $j=1,\ldots,n$:
$$\sum\limits_{i=n+1}^p\Gamma_{ji}^2<\epsilon^2.$$
Therefore, for $i=n+1,\ldots,p$, we obtain the following inequality:
\begin{align*}
    \sum\limits_{i=n+1}^p\dfrac{c_i}{n} &=\dfrac{h(p-n)}{n}+\sum\limits_{i=n+1}^p\sum\limits_{j=1}^k (1+\beta_j)\Gamma_{ji}^2 \lambda_{0,j} + \sum\limits_{i=n+1}^p\sum\limits_{j=1}^n\Big(\bar{d}\dfrac{p}{n}+\alpha_j\sqrt{\dfrac{p}{n}}\Big)\Gamma_{ji}^2\\
    &\leq \dfrac{h(p-n)}{n}+(1+\beta_{\max})\lambda_{0,1}\sum\limits_{i=n+1}^p\sum\limits_{j=1}^k \Gamma_{ji}^2  + \Big(\bar{d}\dfrac{p}{n}+C\sqrt{\dfrac{p}{n}}\Big) \sum\limits_{i=n+1}^p\sum\limits_{j=1}^n\Gamma_{ji}^2\\
    &\leq \dfrac{h(p-n)}{n}+\epsilon^2(1+\beta_{\max})\lambda_{0,1} + \epsilon^2\Big(\bar{d}\dfrac{p}{n}+C\sqrt{\dfrac{p}{n}}\Big)\\
    &\leq \dfrac{h(p-n)}{n}+2\epsilon^2\lambda_{0,1},
\end{align*}
where $\beta_{\max} = \max_{1\leq j\leq k}\beta_j$, and for all sufficiently large $n$. The last inequality follows from
$$\lambda_{0,1}\succ \beta_{\max}\lambda_{0,1} +  \bar{d}\dfrac{p}{n}+C\sqrt{\dfrac{p}{n}}.$$

By arithmetic geometric mean inequality, we obtain the following inequality:
\begin{align*}
    \prod\limits_{i=n+1}^p \dfrac{c_i}{n} &\leq \Big(\dfrac{1}{p-n}\sum\limits_{i=n+1}^p\dfrac{c_i}{n}\Big)^{p-n}\\
    &\leq \bigg[\dfrac{h}{n}+\dfrac{2\epsilon^2}{(p-n)}\lambda_{0,1}\Big)  \bigg]^{p-n}\\
    &= \Big( \dfrac{h}{n}\Big)^{p-n} \bigg[1+\dfrac{2\epsilon^2n}{(p-n)h}\lambda_{0,1}  \bigg]^{p-n}.
\end{align*}

\end{proof}

\begin{lemma}[Lower bound of $c_i$]\label{lem:lower_bound_c}
Under the conditions $A1-A4$ on Section \ref{main-sec:main_results}, $c_i$ has the following lower bound:
    $$\dfrac{c_i}{n}\geq \prod\limits_{j=1}^k \Big((1+\beta_j)\lambda_{0,j}+(\bar{d}\dfrac{p}{n}+\alpha_j\sqrt{\dfrac{p}{n}})+\dfrac{h}{n}\Big)^{\Gamma_{ji}^2}\cdot
    \prod\limits_{j=k+1}^n \Big(\dfrac{\bar{d}p}{n}+\alpha_j\sqrt{\dfrac{p}{n}}+\dfrac{h}{n}\Big)^{\Gamma_{ji}^2}\cdot
    \prod\limits_{j=n+1}^p \Big(\dfrac{h}{n}\Big)^{\Gamma_{ji}^2}.
$$
\end{lemma}
\begin{proof}
By arithmetic–geometric mean inequality, we obtain the following inequality
\begin{align*}
    \dfrac{c_i}{n} 
    &= \sum\limits_{j=1}^k\Gamma_{ji}^2\Big\{(1+\beta_j)\lambda_{0,j}+(\bar{d}\dfrac{p}{n}+\alpha_j\sqrt{\dfrac{p}{n}})+\dfrac{h}{n}\Big\}+ \sum\limits_{j=k+1}^n\Gamma_{ji}^2\Big(\bar{d}\dfrac{p}{n}+\alpha_j\sqrt{\dfrac{p}{n}}+\dfrac{h}{n} \Big) + \sum\limits_{j=n+1}^p\Gamma_{ji}^2\dfrac{h}{n}\\
    &\geq \prod\limits_{j=1}^k \Big((1+\beta_j)\lambda_{0,j}+(\bar{d}\dfrac{p}{n}+\alpha_j\sqrt{\dfrac{p}{n}})+\dfrac{h}{n}\Big)^{\Gamma_{ji}^2}\cdot
    \prod\limits_{j=k+1}^n \Big(\dfrac{\bar{d}p}{n}+\alpha_j\sqrt{\dfrac{p}{n}}+\dfrac{h}{n}\Big)^{\Gamma_{ji}^2}\cdot
    \prod\limits_{j=n+1}^p \Big(\dfrac{h}{n}\Big)^{\Gamma_{ji}^2}.
\end{align*}
\end{proof}

\begin{lemma}[Hoffman-Wielandt theorem ; \citealt{stewart1990matrix}]\label{lem:perturb_eigen}
    For $n\times n$ normal matrices $A,\tilde{A}$, the following inequality holds:
    $$\min_\pi\sqrt{\sum\limits_i (\tilde{\lambda}_{\pi(i)}-\lambda_i)^2}\leq \abs{\abs{\tilde{A}-A}}_F,$$
    where $\pi$ is permutation of the set $[n]=\{1,2,\ldots,n\}$, and $\min_\pi$ denotes the minimum among all permutation of the set $[n]$.
\end{lemma}

\begin{lemma}[Block matrix approximation]\label{lem:block_app}
    For all $\Gamma\in O(p)$, if $\abs{\abs{\Gamma_{1:n,n+1:p}}}_F<\eta$, then it follows that
    $$\displaystyle\inf_{Q_1\in O(n),\;Q_2\in O(p-n)} \abs{\abs{\begin{pmatrix}
Q_1 & 0 \\ 0 & Q_2
\end{pmatrix}-\Gamma}}_F <2\eta,$$
for all $\eta\in(0,1)$.
\end{lemma}
\begin{proof}
    Consider the following block matrix of $\Gamma$:
$$\Gamma = \begin{pmatrix}
    \Gamma_{11} & \Gamma_{12}\\
    \Gamma_{21} & \Gamma_{22}
\end{pmatrix}, $$
where $\Gamma_{11}\in\bbR^{n\times n}$ and $\Gamma_{22}\in\bbR^{(p-n)\times(p-n)}$.

Let $A=I_{p-n}$, $\tilde{A}= I_{p-n} - \Gamma_{12}^T\Gamma_{12}$, and $B_{\eta}=\{\Gamma\in O(p) :\abs{\abs{\Gamma_{12}}}_F<\eta \}$. For all $\Gamma\in B_\eta$, the following inequality holds by Lemma \ref{lem:perturb_eigen}:
$$\sum\limits_{i=1}^{p-n}(\tilde{\lambda}_i-1)^2\leq \abs{\abs{\Gamma_{12}^T\Gamma_{12}}}_F^2<\eta^4,$$
where $\tilde{\lambda}_i$ is $i$th eigenvalue of $\tilde{A}$.    

Consider the singular value decomposition $\Gamma_{22}=UDV^T$, where $U,V\in O(p-n)$ and $D = diag(\sqrt{\tilde{\lambda}_1},\cdots,\sqrt{\tilde{\lambda}_{p-n}})$. Then, the equality holds for all $\Gamma\in B_\eta$:
\begin{align*}
    \inf_{Q\in O(p-n)}\abs{\abs{\Gamma_{22}-Q}}_F^2 &= \inf_{Q\in O(p-n)}\abs{\abs{D-Q}}_F^2\\
    &=\inf_{Q\in O(p-n)}\sum\limits_{i=1}^{p-n}[(\sqrt{\tilde{\lambda}_i}-q_i)^2+(1-q_i^2)]\\
    &=\sum\limits_{i=1}^{p-n}(\sqrt{\tilde{\lambda}_i}-1)^2\\
    &\leq \max_{i\leq p-n} \dfrac{1}{(\sqrt{\tilde{\lambda}_i}+1)^2} \cdot \sum\limits_{i=1}^{p-n}(\tilde{\lambda}_i-1)^2\\
    &\leq \eta^4,
\end{align*}
where $q_i$ is $(i,i)$ element of $Q$. Likewise we attain 
$$\inf_{Q\in O(n)}\abs{\abs{\Gamma_{11}-Q}}_F^2\leq \eta^4.$$

For $\Gamma\in B_\eta$, we obtain the inequality:
\begin{align*}
    \inf_{Q_1\in O(n),Q_2\in O(p-n)}\abs{\abs{\Gamma-\begin{pmatrix}
        Q_1 & 0 \\ 0 & Q_2
    \end{pmatrix}}}_F^2 &= \inf_{Q_1\in O(n)}\abs{\abs{Q_1-\Gamma_{11}}}_F^2 +  \inf_{Q_2\in O(p-n)}\abs{\abs{Q_2-\Gamma_{22}}}_F^2 +\abs{\abs{\Gamma_{12}}}_F^2 + \abs{\abs{\Gamma_{21}}}_F^2\\
    &\leq 2\eta^4 + 2\eta^2\\
    &\leq 4\eta^2.
\end{align*}
\end{proof}

\begin{lemma}\label{lem:post_contract}
    Consider the following subset of $O(p)$:
\begin{align*}
    A_{\epsilon}  = \bigg\{\Gamma\in O(p):\displaystyle\inf_{Q_1\in O(k),\;Q_2\in O(p-k)} \abs{\abs{\Gamma-\begin{pmatrix}
        Q_1 & 0 \\ 
        0& Q_2
    \end{pmatrix}}}_F <\epsilon
    \bigg\}.
\end{align*}
Under the conditions $A1-A5$, suppose that $\epsilon^2a_n\succ np$. Then, the following inequality holds:
$$\dfrac{\displaystyle\int_{A_\epsilon^c} \prod\limits_{i=1}^p c_{i}^{-a_i-\frac{n}{2}+1}(d{\Gamma})}{\displaystyle\int\prod\limits_{i=1}^p c_{i}^{-a_i-\frac{n}{2}+1} (d{\Gamma})}\preccurlyeq\Big(\dfrac{n\lambda_{0,k}+p}{p}\Big)^{-\epsilon^2 a_n}.$$

\end{lemma}
\begin{proof}
Consider the following subset of the orthogonal group $O(p)$: 
\begin{align*}
    C_{\epsilon}  = \bigg\{\Gamma\in O(p):\displaystyle\inf_{Q_2\in O(p-n)} \abs{\abs{\Gamma-\begin{pmatrix}
        I_n & 0 \\ 
        0& Q_2
    \end{pmatrix}}}_F <\epsilon
    \bigg\}.
\end{align*}

Define the following quantities:
\begin{align*}
    K_i &= \Big((1+\beta_i)\lambda_{0,i}  +\Big( \bar{d}\dfrac{p}{n}+\alpha_i\sqrt{\dfrac{p}{n}}\Big)+\dfrac{h}{n}\Big),\quad i=1,\ldots,k,\\
    L_i &= \Big(\bar{d}\dfrac{p}{n}+\alpha_i\sqrt{\dfrac{p}{n}}+\dfrac{h}{n}\Big),\quad i=k+1,\ldots,n,\\
    M &= \Big( \dfrac{h}{n}\Big),\\
    r_n &= 4\epsilon_2^2\dfrac{\lambda_{0,1}}{\lambda_{0,k}},\\
    s_n &= 4\epsilon_2^2\dfrac{n\lambda_{0,1}}{\bar{d}p},\\
    t_n &= \dfrac{2\epsilon_2^2n}{(p-n)h}\lambda_{0,1} .
\end{align*}

Using Lemma \ref{lem:upper_bound_c} and Lemma \ref{lem:lower_bound_c}, we obtain the following inequality:
\begin{align*}
    &\dfrac{\displaystyle\int_{A_\epsilon^c} \prod c_{i}^{-a_i-\frac{n}{2}+1}(d{\Gamma})}{\displaystyle\int_{C_{\epsilon_2}}\prod c_{i}^{-a_i-\frac{n}{2}+1} (d{\Gamma})}\\
    &\leq \dfrac{\bbP(A_\epsilon^c)\cdot\sup_{A_\epsilon^c}\prod c_{i}^{-a_i-\frac{n}{2}+1}}{\bbP(C_{\epsilon_2})\cdot\inf_{C_{\epsilon_2}}\prod c_{i}^{-a_i-\frac{n}{2}+1}}\\
    &\leq \dfrac{1}{\bbP(C_{\epsilon_2})}\cdot \dfrac{\sup_{C_{\epsilon_2}}\prod c_{i}^{a_i+\frac{n}{2}-1}}{\inf_{A_\epsilon^c}\prod c_{i}^{a_i+\frac{n}{2}-1}}\\
    &\leq \dfrac{1}{\bbP(C_{\epsilon_2})}\cdot \dfrac{\prod\limits_{j=1}^k\bigg[K_j  (1+r_n) \bigg]^{(a_j+\frac{n}{2}-1)}  \cdot \prod\limits_{j=k+1}^n\bigg[L_j (1+s_n)\bigg]^{(a_n+\frac{n}{2}-1)}\cdot\bigg[M(1+t_n)\bigg]^{(p-n)(a_p+\frac{n}{2}-1)}}{\inf_{A_\epsilon^c}\prod\limits_{i=1}^p \bigg[\prod\limits_{j=1}^k K_j^{\Gamma_{ji}^2}\cdot
    \prod\limits_{j=k+1}^n L_j^{\Gamma_{ji}^2}\cdot
    \prod\limits_{j=n+1}^p M^{\Gamma_{ji}^2}\bigg]^{a_i+\frac{n}{2}-1}}\\
    &= \dfrac{1}{\bbP(C_{\epsilon_2})}\cdot (1+r_n)^{\sum\limits_{j=1}^k (a_j+\frac{n}{2}-1)}(1+s_n)^{(n-k)(a_n+\frac{n}{2}-1)}(1+t_n)^{(p-n)(a_p+\frac{n}{2}-1)}\\
    &\times \sup_{A_\epsilon^c} \dfrac{\prod\limits_{j=1}^k K_j ^{(a_j+\frac{n}{2}-1)}  \cdot\prod\limits_{j=k+1}^n L_j^{(a_n+\frac{n}{2}-1)}\cdot M^{(p-n)(a_p+\frac{n}{2}-1)}}{\prod\limits_{i=1}^p \bigg[\prod\limits_{j=1}^k K_j^{\Gamma_{ji}^2}\cdot
    \prod\limits_{j=k+1}^n L_j^{\Gamma_{ji}^2}\cdot
    \prod\limits_{j=n+1}^p M^{\Gamma_{ji}^2}\bigg]^{a_i+\frac{n}{2}-1}}.
\end{align*}

Now, we focus on bounding the term
$$\sup_{A_\epsilon^c} \dfrac{\prod\limits_{j=1}^k K_j ^{(a_j+\frac{n}{2}-1)}  \cdot\prod\limits_{j=k+1}^n L_j^{(a_n+\frac{n}{2}-1)}\cdot M^{(p-n)(a_p+\frac{n}{2}-1)}}{\prod\limits_{i=1}^p \bigg[\prod\limits_{j=1}^k K_j^{\Gamma_{ji}^2}\cdot
    \prod\limits_{j=k+1}^n L_j^{\Gamma_{ji}^2}\cdot
    \prod\limits_{j=n+1}^p M^{\Gamma_{ji}^2}\bigg]^{a_i+\frac{n}{2}-1}}.$$

We obtain the following upper bound:
\begin{adjustwidth}{-1.2cm}{-1.2cm}
{\small
\begin{align*}
    &\sup_{A_\epsilon^c} \dfrac{\prod\limits_{j=1}^k K_j ^{(a_j+\frac{n}{2}-1)}  \cdot \prod\limits_{j=k+1}^n L_j^{(a_n+\frac{n}{2}-1)}\cdot M^{(p-n)(a_p+\frac{n}{2}-1)}}{\prod\limits_{i=1}^p \bigg[\prod\limits_{j=1}^k K_j^{\Gamma_{ji}^2}\cdot
    \prod\limits_{j=k+1}^n L_j^{\Gamma_{ji}^2}\cdot
    \prod\limits_{j=n+1}^p M^{\Gamma_{ji}^2}\bigg]^{a_i+\frac{n}{2}-1}}\\
    &= \sup_{A_\epsilon^c}\Bigg[ \dfrac{\prod\limits_{j=1}^k K_j ^{(a_j+\frac{n}{2}-1)}}{\prod\limits_{j=1}^k K_j^{\sum\limits_{i=1}^k\Gamma_{ji}^2(a_i+\frac{n}{2}-1)+\sum\limits_{i=k+1}^n\Gamma_{ji}^2(a_n+\frac{n}{2}-1)+\sum\limits_{i=n+1}^p\Gamma_{ji}^2(a_p+\frac{n}{2}-1)}}\\
    &\times \dfrac{\prod\limits_{j=k+1}^n L_j^{(a_n+\frac{n}{2}-1)}}{\prod\limits_{j=k+1}^n L_j^{\sum\limits_{i=1}^k\Gamma_{ji}^2 (a_i+\frac{n}{2}-1)+\sum\limits_{i=k+1}^n\Gamma_{ji}^2 (a_n+\frac{n}{2}-1)+\sum\limits_{i=n+1}^p\Gamma_{ji}^2 (a_p+\frac{n}{2}-1)}}\times \dfrac{M^{(p-n)(a_p+\frac{n}{2}-1)}}{M^{\sum\limits_{i=1}^k\sum\limits_{j=n+1}^p\Gamma_{ji}^2 (a_i+\frac{n}{2}-1)+\sum\limits_{i=k+1}^n\sum\limits_{j=n+1}^p\Gamma_{ji}^2 (a_n+\frac{n}{2}-1)+\sum\limits_{i=n+1}^p\sum\limits_{j=n+1}^p\Gamma_{ji}^2 (a_p+\frac{n}{2}-1)}} \Bigg]\\
    &= \sup_{A_\epsilon^c}\Bigg[ \dfrac{\prod\limits_{j=1}^k K_j ^{(a_j+\frac{n}{2}-1)}}{\prod\limits_{j=1}^k K_j^{\sum\limits_{i=1}^k\Gamma_{ji}^2(a_i+\frac{n}{2}-1)+\sum\limits_{i=k+1}^p\Gamma_{ji}^2(a_k+\frac{n}{2}-1)}} \cdot \dfrac{1}{\prod\limits_{j=1}^k K_i^{\sum\limits_{i=k+1}^n\Gamma_{ji}^2(a_n-a_k)+\sum\limits_{i=n+1}^p\Gamma_{ji}^2(a_p-a_k)}}\\
    &\times \dfrac{\prod\limits_{j=k+1}^n L_j^{(1-\sum\limits_{i=k+1}^n\Gamma_{ji}^2)(a_n+\frac{n}{2}-1)}}{\prod\limits_{j=k+1}^n L_j^{\sum\limits_{i=1}^k\Gamma_{ji}^2 (a_i+\frac{n}{2}-1)+\sum\limits_{i=n+1}^p\Gamma_{ji}^2 (a_p+\frac{n}{2}-1)}}
    \times \dfrac{M^{(p-n-\sum\limits_{i=n+1}^p\sum\limits_{j=n+1}^p\Gamma_{ji}^2)(a_p+\frac{n}{2}-1)}}{M^{\sum\limits_{i=1}^k\sum\limits_{j=n+1}^p\Gamma_{ji}^2 (a_i+\frac{n}{2}-1)+\sum\limits_{i=k+1}^n\sum\limits_{j=n+1}^p\Gamma_{ji}^2 (a_n+\frac{n}{2}-1)}} \Bigg]\\
    &\leq \sup_{O(p)}\Bigg[\dfrac{\prod\limits_{j=1}^k K_j ^{(a_j+\frac{n}{2}-1)}}{\prod\limits_{j=1}^k K_j^{\sum\limits_{i=1}^k\Gamma_{ji}^2(a_i+\frac{n}{2}-1)+\sum\limits_{i=k+1}^p\Gamma_{ji}^2(a_k+\frac{n}{2}-1)}}\Bigg] \\
    &\times \sup_{A_\epsilon^c}\Bigg[ \dfrac{1}{\prod\limits_{j=1}^k K_i^{\sum\limits_{i=k+1}^n\Gamma_{ji}^2(a_n-a_k)+\sum\limits_{i=n+1}^p\Gamma_{ji}^2(a_p-a_k)}}\times \dfrac{\prod\limits_{j=k+1}^n L_j^{(1-\sum\limits_{i=k+1}^n\Gamma_{ji}^2)(a_n+\frac{n}{2}-1)}}{\prod\limits_{j=k+1}^n L_j^{\sum\limits_{i=1}^k\Gamma_{ji}^2 (a_i+\frac{n}{2}-1)+\sum\limits_{i=n+1}^p\Gamma_{ji}^2 (a_p+\frac{n}{2}-1)}}\times \dfrac{M^{(p-n-\sum\limits_{i=n+1}^p\sum\limits_{j=n+1}^p\Gamma_{ji}^2)(a_p+\frac{n}{2}-1)}}{M^{\sum\limits_{i=1}^k\sum\limits_{j=n+1}^p\Gamma_{ji}^2 (a_i+\frac{n}{2}-1)+\sum\limits_{i=k+1}^n\sum\limits_{j=n+1}^p\Gamma_{ji}^2 (a_n+\frac{n}{2}-1)}} \Bigg]\\
    &\leq  \sup_{A_\epsilon^c}\Bigg[ \dfrac{1}{\prod\limits_{j=1}^k K_i^{\sum\limits_{i=k+1}^n\Gamma_{ji}^2(a_n-a_k)+\sum\limits_{i=n+1}^p\Gamma_{ji}^2(a_p-a_k)}}\times \dfrac{\prod\limits_{j=k+1}^n L_j^{(1-\sum\limits_{i=k+1}^p\Gamma_{ji}^2)(a_n+\frac{n}{2}-1)}}{\prod\limits_{j=k+1}^n L_j^{\sum\limits_{i=1}^k\Gamma_{ji}^2 (a_i+\frac{n}{2}-1)+\sum\limits_{i=n+1}^p\Gamma_{ji}^2 (a_p-a_n)}}\times \dfrac{M^{(p-n-\sum\limits_{i=n+1}^p\sum\limits_{j=n+1}^p\Gamma_{ji}^2)(a_p+\frac{n}{2}-1)}}{M^{\sum\limits_{i=1}^k\sum\limits_{j=n+1}^p\Gamma_{ji}^2 (a_i+\frac{n}{2}-1)+\sum\limits_{i=k+1}^n\sum\limits_{j=n+1}^p\Gamma_{ji}^2 (a_n+\frac{n}{2}-1)}} \Bigg]\\
    &\leq  \sup_{A_\epsilon^c}\Bigg[ \dfrac{1}{ K_k^{\sum\limits_{j=1}^k\sum\limits_{i=k+1}^n\Gamma_{ji}^2 (a_n-a_k)+\sum\limits_{j=1}^k\sum\limits_{i=n+1}^p\Gamma_{ji}^2 (a_p-a_k)}}\times \dfrac{L_n^{(n-k-\sum\limits_{j=k+1}^n\sum\limits_{i=k+1}^p\Gamma_{ji}^2)(a_n+\frac{n}{2}-1)}}{L_{k+1}^{\sum\limits_{i=n+1}^p\sum\limits_{j=k+1}^n\Gamma_{ji}^2 (a_p-a_n)}}\times \dfrac{M^{(p-n-\sum\limits_{i=n+1}^p\sum\limits_{j=n+1}^p\Gamma_{ji}^2)(a_p+\frac{n}{2}-1)}}{M^{k (a_1+\frac{n}{2}-1) +\sum\limits_{i=k+1}^n\sum\limits_{j=n+1}^p\Gamma_{ji}^2(a_n+\frac{n}{2}-1)}}\Bigg]\\
    &=  \sup_{A_\epsilon^c}\Bigg[ \dfrac{1}{ K_k^{\sum\limits_{j=1}^k\sum\limits_{i=k+1}^p\Gamma_{ji}^2 (a_n-a_k)+\sum\limits_{j=1}^k\sum\limits_{i=n+1}^p\Gamma_{ji}^2 (a_p-a_k)}}\times \dfrac{L_n^{\sum\limits_{j=k+1}^n\sum\limits_{i=1}^k \Gamma_{ji}^2(a_n+\frac{n}{2}-1)}}{L_{k+1}^{\sum\limits_{i=n+1}^p\sum\limits_{j=k+1}^n\Gamma_{ji}^2 (a_p-a_n)}}\times \dfrac{M^{\sum\limits_{i=1}^n\sum\limits_{j=n+1}^p\Gamma_{ji}^2(a_p+\frac{n}{2}-1)}}{M^{k (a_1+\frac{n}{2}-1)+\sum\limits_{i=k+1}^n\sum\limits_{j=n+1}^p\Gamma_{ji}^2(a_n +\frac{n}{2}-1) }}\Bigg],
\end{align*}
}
\end{adjustwidth}

where the second inequality follows from the fact that the term
$$\prod\limits_{j=1}^k K_j^{\sum\limits_{i=1}^k\Gamma_{ji}^2(a_i+\frac{n}{2}-1)+\sum\limits_{i=k+1}^p\Gamma_{ji}^2(a_k+\frac{n}{2}-1)},$$
is minimized when $\Gamma_{jj}^2 = 1$ for $j=1,\ldots,k$. The third inequality is satisfied due to the ordering conditions $K_1\geq \cdots\geq K_k$ and $L_{k+1}\geq \cdots\geq  L_n$.  

Since Lemma $\ref{lem:block_app}$ implies that $\abs{\abs{\Gamma_{1:n,n+1:p}}}_F>\epsilon/2$ on the $A_\epsilon^c$, we obtain the inequality:
\begin{align*}
    &\dfrac{1}{ K_k^{\sum\limits_{j=1}^k\sum\limits_{i=k+1}^p\Gamma_{ji}^2 (a_n-a_k)+\sum\limits_{j=1}^k\sum\limits_{i=n+1}^p\Gamma_{ji}^2 (a_p-a_k)}}\times \dfrac{L_n^{\sum\limits_{j=k+1}^n\sum\limits_{i=1}^k \Gamma_{ji}^2(a_n+\frac{n}{2}-1)}}{L_{k+1}^{\sum\limits_{i=n+1}^p\sum\limits_{j=k+1}^n\Gamma_{ji}^2 (a_p-a_n)}}\times \dfrac{M^{\sum\limits_{i=1}^n\sum\limits_{j=n+1}^p\Gamma_{ji}^2(a_p+\frac{n}{2}-1)}}{M^{k (a_1+\frac{n}{2}-1)+\sum\limits_{i=k+1}^n\sum\limits_{j=n+1}^p\Gamma_{ji}^2(a_n +\frac{n}{2}-1) }}\\
    &\leq \dfrac{1}{ K_k^{\sum\limits_{j=1}^k\sum\limits_{i=k+1}^p\Gamma_{ji}^2 (a_n-a_k)}}\times L_n^{\sum\limits_{j=k+1}^n\sum\limits_{i=1}^k \Gamma_{ji}^2(a_n+\frac{n}{2}-1)}\times \dfrac{1}{M^{k (a_1+\frac{n}{2}-1)}}\\
    &= \Big(\dfrac{L_n}{K_k}\Big)^{\sum\limits_{j=1}^k\sum\limits_{i=k+1}^p\Gamma_{ji}^2 (a_n-a_k)}\times \dfrac{L_n^{\sum\limits_{j=1}^k\sum\limits_{i=k+1}^p\Gamma_{ji}^2 (a_k +\frac{n}{2}-1)}}{M^{k (a_1+\frac{n}{2}-1)}},
\end{align*}
where the first inequality is derived from $K_k,L_{k+1}>1$, $M<1$, and $a_p\geq a_n$.

Next, we derive an upper bound for the remaining term: 
\begin{align*}
    &(1+r_n)^{\sum\limits_{i=1}^k (a_i+\frac{n}{2}-1)}(1+s_n)^{(n-k)(a_n+\frac{n}{2}-1)}(1+t_n)^{(p-k)(a_p+\frac{n}{2}-1)}\\
    &\preccurlyeq\exp(\epsilon_2^2\dfrac{\lambda_{0,1}}{\lambda_{0,k}}a_k) \cdot\exp(4\epsilon_2^2\dfrac{n\lambda_{0,1}}{\bar{d}p} n a_n) \cdot \exp(\Big(\dfrac{2\epsilon_2^2n}{(p-n)h}\lambda_{0,1}\Big)p a_p ),
\end{align*}
using the bound $\log(1+x)\leq x$ for all small $x>0$.

Therefore, we obtain the final upper bound
\begin{align*}
    &\dfrac{\displaystyle\int_{A_\epsilon^c} \prod c_{i}^{-a_i-\frac{n}{2}+1}(d{\Gamma})}{\displaystyle\int_{C_{\epsilon_2}}\prod c_{i}^{-a_i-\frac{n}{2}+1} (d{\Gamma})}\\
    &\leq \dfrac{1}{\bbP(C_{\epsilon_2})}\cdot \exp(\epsilon_2^2\dfrac{\lambda_{0,1}}{\lambda_{0,k}}a_k) \cdot\exp(4\epsilon_2^2\dfrac{n\lambda_{0,1}}{\bar{d}p} n a_n) \cdot \exp(\Big(\dfrac{2\epsilon_2^2n}{(p-n)h}\lambda_{0,1}\Big)p a_p )\\
    &\times \sup_{A_\epsilon^c}\Bigg[\Big(\dfrac{L_n}{K_k}\Big)^{\sum\limits_{j=1}^k\sum\limits_{i=k+1}^p\Gamma_{ji}^2 (a_n-a_k)}\times \dfrac{L_n^{\sum\limits_{j=1}^k\sum\limits_{i=k+1}^p\Gamma_{ji}^2 (a_k +\frac{n}{2}-1)}}{M^{k (a_1+\frac{n}{2}-1)}}\Bigg]\\
    &\preccurlyeq \Big(\dfrac{c\sqrt{n}}{\epsilon_2}\Big)^{np-n^2}\cdot  
    \exp(\dfrac{n^2\epsilon_2^2\lambda_{0,1}a_n}{p})\cdot
    \exp(\dfrac{n\epsilon_2^2\lambda_{0,1} a_p}{h} )\cdot \Big(\dfrac{L_n}{K_k}\Big)^{\sum\limits_{i=1}^k\sum\limits_{j=k+1}^n \Gamma_{ji}^2 a_n}\\
    &\preccurlyeq \Big(\dfrac{K_k}{L_n}\Big)^{-\epsilon^2 a_n/4}\\
    &\asymp \Big(\dfrac{n\lambda_{0,k}+p}{p}\Big)^{-\epsilon^2 a_n},
\end{align*}
where $\epsilon_2^2\prec \dfrac{hp}{a_p\lambda_{0,1}}\wedge \dfrac{p^2}{n\lambda_{0,1}a_n}$, and $\epsilon^2 a_n\succ np$.
    
\end{proof}

\begin{lemma}\label{lem:rearrange}
Suppose that $a_1<\cdots < a_n$ and $b_1>\cdots>b_n$. Consider the permutation function $\pi$, which is a permutation of the set $[n]=\{1,\ldots,n\}$. If $\pi$ is not the identity permutation, then there exists an index $i\in\{1,\ldots,n\}$ such that $\pi(i)\neq i$. In this case, the following inequality holds
$$\sum\limits_{l=1}^n a_l b_{\pi(l)}\geq \sum\limits_{l=1}^n a_l b_l  + \min_{l<n}(a_{l+1}-a_l) \cdot \min_{l<n}(b_l - b_{l+1}) .$$
\end{lemma}
\begin{proof}
Let $i$ be the smallest integer which satisfies $\pi(i) \neq i$, then $\pi(i) = j > i$. Furthermore, there exists $k>i$ such that $\pi(k) = i$. Now, we define a new permutation function $\sigma$ as follows:
$$
    \sigma(l) = \begin{cases}
        \pi(l) \quad & \text{for } l\neq i,k,\\
        i \quad & \text{for } l=i,\\
        j \quad & \text{for } l=k.
    \end{cases}
$$
Using the newly defined permutation function $\sigma$, we have
\begin{align*}
    \sum\limits_{l=1}^n a_lb_{\pi(l)} &= \sum\limits_{l=1}^n a_l b_{\sigma(l)} + (a_ib_j + a_kb_i ) - (a_i b_i + a_kb_j)\\
    &= \sum\limits_{l=1}^n a_l b_{\sigma(l)} + (a_k - a_i)(b_i - b_j),
\end{align*}
Since the rearrangement inequality guarantees that:
$$\sum\limits_{l=1}^n a_l b_{\sigma(l)}\geq \sum\limits_{l=1}^n a_l b_l,$$
if follows that
$$\sum\limits_{l=1}^n a_lb_{\pi(l)}\geq \sum\limits_{l=1}^n a_l b_l+(a_k - a_i)(b_i - b_j).$$

For any non-identity permutation function $\pi$, we obtain the following inequality
\begin{align*}
    \sum\limits_{l=1}^n a_lb_{\pi(l)}&\geq \sum\limits_{l=1}^n a_l b_l +  \inf_{i,j,k:i<j,i<k}(a_k - a_i)(b_i - b_j)\\
    &\geq \sum\limits_{l=1}^n a_l b_l  + \min_{l<n}(a_{l+1}-a_l) \cdot \min_{l<n}(b_l - b_{l+1}).
\end{align*}
\end{proof}

\begin{lemma}\label{lem:upp_bound_D/E}
Assume that conditions $A1-A5$ hold. Define the sets as follows:
\begin{align*}
    D_\eta  &= \bigg\{\Gamma\in O(p):\displaystyle\inf_{Q_2\in O(p-k)} \abs{\abs{\begin{pmatrix}
       I_k & 0 \\ 
        0& Q_2
    \end{pmatrix}-\Gamma}}_F <\eta\bigg\}, \\
    E_{\eta} &= \bigg\{\Gamma\in O(p):\displaystyle\inf_{Q_2\in O(p-k)} \abs{\abs{\begin{pmatrix}
       S & 0 \\ 
        0& Q_2
    \end{pmatrix}-\Gamma}}_F <\eta\bigg\},
\end{align*}
where $S\in O(k)$ has non-negative diagonal elements, and satisfies $\abs{\abs{I_k-S}}_F\geq \epsilon$. If $\eta\prec \dfrac{\lambda_{0,k}}{\lambda_{0,1}}$, then the following inequality holds
$$\dfrac{\sup\limits_{D_{\eta}}\prod\limits_{i=1}^k \Big(\dfrac{c_i}{n}\Big)^{a_i+n/2-1} }{\inf\limits_{E_{\eta}}\prod\limits_{i=1}^k \Big(\dfrac{c_i}{n}\Big)^{a_i+n/2-1} }\preccurlyeq \dfrac{\exp(n\eta^2 \dfrac{\hat{\lambda}_1 }{\hat{\lambda}_k})}{\exp(\dfrac{\epsilon}{2\sqrt{k}}\min_{l<k}(a_{l+1} - a_l) \cdot \min_{l<k}\log(\dfrac{\hat{\lambda}_l}{\hat{\lambda}_{l+1}}))}.$$
\end{lemma}

\begin{proof}
\begin{align*}
    \dfrac{\sup\limits_{D_{\eta}}\prod\limits_{i=1}^k \Big(\dfrac{c_i}{n}\Big)^{a_i+n/2-1} }{\inf\limits_{E_{\eta}}\prod\limits_{i=1}^k \Big(\dfrac{c_i}{n}\Big)^{a_i+n/2-1} }
    &= \dfrac{\sup\limits_{D_{\eta}}\prod\limits_{i=1}^k \bigg[\dfrac{h}{n}+\sum\limits_{j=1}^n\Gamma_{ji}^2\hat{\lambda}_j \bigg]^{a_i+n/2-1}}{\inf\limits_{E_{\eta}} \prod\limits_{i=1}^k \bigg[\dfrac{h}{n}+\sum\limits_{j=1}^n\Gamma_{ji}^2\hat{\lambda}_j \bigg]^{a_i+ n/2-1}}.
\end{align*}

At first, we focus on the following term
\begin{equation}\label{eq:sup_D}
    \sup\limits_{D_{\eta}}\prod\limits_{i=1}^k \bigg[\dfrac{h}{n}+\sum\limits_{j=1}^n\Gamma_{ji}^2\hat{\lambda}_j \bigg]^{a_i+n/2-1}.
\end{equation}

For $\Gamma\in D_\eta$, $(1-\Gamma_{ii}^2) + \sum\limits_{\substack{j=1 \\ j \neq i}}^p\Gamma_{ji}^2 <\eta^2$ holds and it implies $\Gamma_{ii}^2\geq 1-\eta^2$.

Therefore, we obtain the upper bound of \eqref{eq:sup_D} as follows
\begin{align*}
    \eqref{eq:sup_D}&\leq \sup\limits_{D_{\eta}}\prod\limits_{i=1}^k \bigg[\dfrac{h}{n}+ \Gamma_{ii}^2\hat{\lambda}_i+\sum\limits_{\substack{j=1 \\ j \neq i}}\Gamma_{ji}^2\hat{\lambda}_j \bigg]^{a_i+n/2-1}\\
    &\leq \sup\limits_{D_{\eta}}\prod\limits_{i=1}^k \bigg[\dfrac{h}{n}+ \Gamma_{ii}^2\hat{\lambda}_i + (1-\Gamma_{ii}^2)\hat{\lambda}_1\bigg]^{a_i+n/2-1}\\
    &\leq \sup\limits_{D_{\eta}}\prod\limits_{i=1}^k \bigg[\dfrac{h}{n}+ \hat{\lambda}_i + \eta^2\hat{\lambda}_1\bigg]^{a_i+n/2-1}\\
    &\leq \prod\limits_{i=1}^k \bigg[\Big(\dfrac{h}{n}+ \hat{\lambda}_i \Big)(1+ \eta^2\dfrac{\hat{\lambda}_1}{\hat{\lambda}_i})\bigg]^{a_i+n/2-1}\\
    &\leq \prod\limits_{i=1}^k \bigg[\dfrac{h}{n}+ \hat{\lambda}_i \bigg]^{a_i+n/2-1} \cdot\Big(1+ \eta^2\dfrac{\hat{\lambda}_1}{\hat{\lambda}_k}\Big)^{\sum\limits_{i=1}^k (a_i+n/2-1)}.
\end{align*}

Now, we focus on the following term
\begin{equation}\label{eq:inf_E}
    \inf\limits_{E_{\eta}} \prod\limits_{i=1}^k \bigg[\dfrac{h}{n}+\sum\limits_{j=1}^n\Gamma_{ji}^2\hat{\lambda}_j \bigg]^{a_i+ n/2-1}.
\end{equation}

For $\Gamma\in E_{\eta}$, $\sum\limits_{j=1}^k(\Gamma_{ij}^2- S_{ij}^2)\geq -\eta^2$ for $i=1,\ldots,k$. Therefore, we obtain the lower bound of \eqref{eq:inf_E}  as follows
\begin{align*}
    \eqref{eq:inf_E} &\geq \inf\limits_{E_{\eta}}\prod\limits_{i=1}^k \bigg[\sum\limits_{j=1}^k\Gamma_{ji}^2\Big(\hat{\lambda}_j +\dfrac{h}{n} \Big)\bigg]^{a_i+ n/2-1}\\
     & = \inf\limits_{E_{\eta}}\prod\limits_{i=1}^k \bigg[\sum\limits_{j=1}^k S_{ji}^2\Big(\hat{\lambda}_j +\dfrac{h}{n} \Big) + \sum\limits_{j=1}^k(\Gamma_{ij}^2- S_{ij}^2) \Big(\hat{\lambda}_j +\dfrac{h}{n} \Big)\bigg]^{a_i+ n/2-1}\\
     &\geq \inf\limits_{E_{\eta}} \prod\limits_{i=1}^k \bigg[\sum\limits_{j=1}^k S_{ji}^2\Big(\hat{\lambda}_j +\dfrac{h}{n} \Big) -\eta^2 \sum\limits_{j=1}^k\Big(\hat{\lambda}_j +\dfrac{h}{n} \Big)\bigg]^{a_i+ n/2-1}\\
     &\geq  \prod\limits_{i=1}^k \bigg[\sum\limits_{j=1}^k S_{ji}^2\Big(\hat{\lambda}_j +\dfrac{h}{n} \Big) -k\eta^2\Big(\hat{\lambda}_1 +\dfrac{h}{n} \Big)\bigg]^{a_i+ n/2-1}\\
     &\geq  \prod\limits_{i=1}^k \bigg[\sum\limits_{j=1}^k S_{ji}^2\Big(\hat{\lambda}_j +\dfrac{h}{n} \Big)\Big(1 -k\eta^2 \dfrac{\hat{\lambda}_1 +\dfrac{h}{n} }{\hat{\lambda}_k +\dfrac{h}{n} }\Big)\bigg]^{a_i+ n/2-1}\\
     &\geq  \prod\limits_{i=1}^k \bigg[\sum\limits_{j=1}^k S_{ji}^2\Big(\hat{\lambda}_j +\dfrac{h}{n} \Big)\bigg]^{a_i+ n/2-1} \cdot \Big(1 -2k\eta^2 \dfrac{\hat{\lambda}_1 }{\hat{\lambda}_k}\Big)^{\sum\limits_{i=1}^k (a_i+n/2-1)}.
\end{align*}

Consider the $k\times k$ matrix $M$, whose $(j,i)$-th entry is given by $S_{ji}^2$. Since $S$ is an orthonormal matrix, the matrix $M$ becomes a doubly stochastic matrix, meaning that the sum of each row and the sum of each column both equal $1$. By the Birkhoff-von Neumann theorem, any doubly stochastic matrix can be expressed as a convex combination of permutation matrices. Therefore, there exist nonnegative weights $w_0,w_1,\ldots,w_l\geq 0$ satisfying $\sum w_i = 1$, such that:
$$M = w_0 I_k + \sum\limits_{i=1}^l w_i P_i,$$
where $I_k,P_1,\ldots,P_l$ are all $k\times k$ permutation matrices.

Since $\abs{\abs{S-I_k}}_F\geq \epsilon$, we obtain the following inequality
\begin{align*}
    \abs{\abs{S-I_k}}_F^2 &= \sum\limits_{j=1}^k\Big((S_{jj}-1)^2 + \sum\limits_{i\neq j}S_{ji}^2\Big)\\
    &=\sum\limits_{j=1}^k(2-2S_{jj})\\
    &\leq \sum\limits_{j=1}^k(2-2S_{jj}^2)\\
    &= \sum\limits_{j=1}^k\Big((S_{jj}^2-1)^2 + \sum\limits_{i\neq j}(S_{ji}^2)^2\Big)\\
    &= \abs{\abs{M-I_k}}_F^2,
\end{align*}
which implies $\abs{\abs{M-I_k}}_F^2\geq \epsilon^2$. Therefore, the inequality holds
\begin{align*}
    \abs{\abs{M-I_k}}_F &= \abs{\abs{\sum\limits_{i=1}^l w_l P_l - (1-w_0)I_k}}_F\\
    &= \abs{\abs{\sum\limits_{i=1}^l w_l (P_l-I_k)}}_F\\
    &\leq \sum\limits_{i=1}^l w_l \abs{\abs{P_l-I_k}}_F\\
    &\leq \sum\limits_{i=1}^l w_l (\abs{\abs{P_l}}_F+ \abs{\abs{I_k}}_F)\\
    & = 2\sqrt{k} (1-w_0),
\end{align*}
which implies $w_0< 1 - \dfrac{\epsilon}{2\sqrt{k}}$.

Let $g(A) = \sum\limits_{i=1}^k(a_i+ \dfrac{n}{2}-1)\log\Big(\sum\limits_{j=1}^k A_{ji}\Big(\hat{\lambda}_j +\dfrac{h}{n} \Big)\Big)$ for doubly stochastic matrix $A\in \bbR^{k\times k}$, then the inequality holds
\begin{align*}
    g(M) &= \sum\limits_{i=1}^k(a_i+ \dfrac{n}{2}-1)\log\Big(\sum\limits_{j=1}^k M_{ji}\Big(\hat{\lambda}_j +\dfrac{h}{n} \Big)\Big)\\
    &= \sum\limits_{i=1}^k(a_i+ \dfrac{n}{2}-1)\log\Big(\sum\limits_{j=1}^k \Big(w_0[I_{k}]_{ji}+\sum\limits_l w_l [P_l]_{ji}\Big)\Big(\hat{\lambda}_j +\dfrac{h}{n} \Big)\Big)\\
    &=\sum\limits_{i=1}^k(a_i+ \dfrac{n}{2}-1)\log\Big(  w_0 \sum\limits_{j=1}^k [I_{k}]_{ji}\Big(\hat{\lambda}_j +\dfrac{h}{n} \Big) +\sum\limits_l w_l  \sum\limits_{j=1}^k [P_l]_{ji}\Big(\hat{\lambda}_j +\dfrac{h}{n} \Big)\Big)\\
    &\geq \sum\limits_{i=1}^k(a_i+ \dfrac{n}{2}-1) \bigg[w_0 \log\Big(\sum\limits_{j=1}^k [I_{k}]_{ji}\Big(\hat{\lambda}_j +\dfrac{h}{n} \Big)\Big) + \sum\limits_l w_l \log\Big(\sum\limits_{j=1}^k [P_l]_{ji}\Big(\hat{\lambda}_j +\dfrac{h}{n} \Big)\Big) \bigg]\\
    & = w_0 g(I_k) + \sum\limits_l w_l g(P_l),
\end{align*}
where the inequality holds by Jensen inequality since $\log$ is concave.

Since $a_1<\cdots< a_k$ and $\hat{\lambda}_1>\cdots>\hat{\lambda}_k$, by applying Lemma \ref{lem:rearrange}, we obtain the following inequality
\begin{align*}
    g(P_l) &\geq g(I_k) + \min_{l<k}(a_{l+1} - a_l) \cdot \min_{l<k}(\log\hat{\lambda}_l - \log\hat{\lambda}_{l+1})\\
    &= g(I_k) + \min_{l<k}(a_{l+1} - a_l) \cdot \min_{l<k}\log(\dfrac{\hat{\lambda}_l}{\hat{\lambda}_{l+1}}).
\end{align*}

Therefore, we obtain the lower bound of $g(M)$ as follows:
\begin{align*}
    g(M) &\geq  w_0 g(I_k) + \sum\limits_l w_l g(P_l)\\
    &\geq w_0 g(I_k) + \sum\limits_l w_l (g(I_k) + \min_{l<k}(a_{l+1} - a_l) \cdot \min_{l<k}\log(\dfrac{\hat{\lambda}_l}{\hat{\lambda}_{l+1}}))\\
    & = g(I_k) + \sum\limits_l w_l \cdot \min_{l<k}(a_{l+1} - a_l) \cdot \min_{l<k}\log(\dfrac{\hat{\lambda}_l}{\hat{\lambda}_{l+1}})\\
    &\geq g(I_k) + \dfrac{\epsilon}{2\sqrt{k}}\min_{l<k}(a_{l+1} - a_l) \cdot \min_{l<k}\log(\dfrac{\hat{\lambda}_l}{\hat{\lambda}_{l+1}}).
\end{align*}

The lower bound of $\eqref{eq:inf_E}$ is given by
\begin{align*}
    \eqref{eq:inf_E} \geq \prod\limits_{i=1}^k \bigg[\hat{\lambda}_i +\dfrac{h}{n} \bigg]^{a_i+ n/2-1} \exp(\dfrac{\epsilon}{2\sqrt{k}}\min_{l<k}(a_{l+1} - a_l) \cdot \min_{l<k}\log(\dfrac{\hat{\lambda}_l}{\hat{\lambda}_{l+1}})) \cdot \Big(1 -2k\eta^2 \dfrac{\hat{\lambda}_1 }{\hat{\lambda}_k}\Big)^{\sum\limits_{i=1}^k (a_i+n/2-1)}.
\end{align*}

By appending the upper bound of \eqref{eq:sup_D} and the lower bound of \eqref{eq:inf_E}, the following inequality holds:
\begin{align*}
    &\dfrac{\sup\limits_{D_{\eta}}\prod\limits_{i=1}^k \Big(\dfrac{c_i}{n}\Big)^{a_i+n/2-1} }{\inf\limits_{E_{\eta}}\prod\limits_{i=1}^k \Big(\dfrac{c_i}{n}\Big)^{a_i+n/2-1} }\\
    &\leq \dfrac{\prod\limits_{i=1}^k \bigg[\dfrac{h}{n}+ \hat{\lambda}_i \bigg]^{a_i+n/2-1} \cdot\Big(1+ \eta^2\dfrac{\hat{\lambda}_1}{\hat{\lambda}_k}\Big)^{\sum\limits_{i=1}^k (a_i+n/2-1)}}{\prod\limits_{i=1}^k\bigg[\hat{\lambda}_i +\dfrac{h}{n} \bigg]^{a_i+ n/2-1} \exp(\dfrac{\epsilon}{2\sqrt{k}}\min_{l<k}(a_{l+1} - a_l) \cdot \min_{l<k}\log(\dfrac{\hat{\lambda}_l}{\hat{\lambda}_{l+1}})) \cdot \Big(1 -2k\eta^2 \dfrac{\hat{\lambda}_1 }{\hat{\lambda}_k}\Big)^{\sum\limits_{i=1}^k (a_i+n/2-1)}}\\
    & \asymp \dfrac{\exp(n\eta^2 \dfrac{\hat{\lambda}_1 }{\hat{\lambda}_k})}{\exp(\dfrac{\epsilon}{2\sqrt{k}}\min_{l<k}(a_{l+1} - a_l) \cdot \min_{l<k}\log(\dfrac{\hat{\lambda}_l}{\hat{\lambda}_{l+1}}))}.
\end{align*}

\end{proof}
\begin{lemma}\label{lem:upp_bound_D/E_2}
Assume that conditions $A1-A5$ hold. Define the sets as follows
\begin{align*}
    D_\eta  &= \bigg\{\Gamma\in O(p):\displaystyle\inf_{Q_2\in O(p-k)} \abs{\abs{\begin{pmatrix}
       I_k & 0 \\ 
        0& Q_2
    \end{pmatrix}-\Gamma}}_F <\eta\bigg\}, \\
    E_{\eta} &= \bigg\{\Gamma\in O(p):\displaystyle\inf_{Q_2\in O(p-k)} \abs{\abs{\begin{pmatrix}
       S & 0 \\ 
        0& Q_2
    \end{pmatrix}-\Gamma}}_F <\eta\bigg\},
\end{align*}
where $S\in O(k)$ with non-negative diagonal elements, and  $\abs{\abs{P-S}}_F\geq \epsilon$ for all permutation matrix $P$. If $a_1=\cdots = a_k$, then the inequality holds
$$\dfrac{\sup\limits_{D_{\eta}}\prod\limits_{i=1}^k \Big(\dfrac{c_i}{n}\Big)^{a_i+n/2-1} }{\inf\limits_{E_{\eta}}\prod\limits_{i=1}^k \Big(\dfrac{c_i}{n}\Big)^{a_i+n/2-1} }\preccurlyeq  \dfrac{\exp(n\eta^2\dfrac{\hat{\lambda}_1}{\hat{\lambda}_k})}{\exp(n\epsilon^2 \min\limits_{i<k}\log(\dfrac{\hat{\lambda}_i}{\hat{\lambda}_{i+1}}))}.$$
\end{lemma}

\begin{proof}
    By proof of Lemma \eqref{lem:upp_bound_D/E}, we obtain the following
    \begin{align}\label{eq:D/E_2}
    \dfrac{\sup\limits_{D_{\eta}}\prod\limits_{i=1}^k \Big(\dfrac{c_i}{n}\Big)^{a_i+n/2-1} }{\inf\limits_{E_{\eta}}\prod\limits_{i=1}^k \Big(\dfrac{c_i}{n}\Big)^{a_i+n/2-1} }
    = \dfrac{\sup\limits_{D_{\eta}}\prod\limits_{i=1}^k \bigg[\dfrac{h}{n}+\sum\limits_{j=1}^n\Gamma_{ji}^2\hat{\lambda}_j \bigg]^{a_1+n/2-1}}{\inf\limits_{E_{\eta}} \prod\limits_{i=1}^k \bigg[\dfrac{h}{n}+\sum\limits_{j=1}^n\Gamma_{ji}^2\hat{\lambda}_j \bigg]^{a_1+ n/2-1}}.
    \end{align}

    At first, the upper bound of numerator of \eqref{eq:D/E_2} is given by
    \begin{align*}
        \sup\limits_{D_{\eta}}\prod\limits_{i=1}^k \bigg[\dfrac{h}{n}+\sum\limits_{j=1}^n\Gamma_{ji}^2\hat{\lambda}_j \bigg]^{a_1+n/2-1}\leq \prod\limits_{i=1}^k \bigg[\dfrac{h}{n}+ \hat{\lambda}_i \bigg]^{a_1+n/2-1} \cdot\Big(1+ \eta^2\dfrac{\hat{\lambda}_1}{\hat{\lambda}_k}\Big)^{k(a_1+n/2-1)}.
    \end{align*}
    
    Next, the lower bound of denominator of \eqref{eq:D/E_2} is given by
    \begin{align*}
        \inf\limits_{E_{\eta}} \prod\limits_{i=1}^k \bigg[\dfrac{h}{n}+\sum\limits_{j=1}^n\Gamma_{ji}^2\hat{\lambda}_j \bigg]^{a_1+ n/2-1}\geq \prod\limits_{i=1}^k \bigg[\sum\limits_{j=1}^k S_{ji}^2\Big(\hat{\lambda}_j +\dfrac{h}{n} \Big)\bigg]^{a_1+ n/2-1} \cdot \Big(1 -2k\eta^2 \dfrac{\hat{\lambda}_1 }{\hat{\lambda}_k}\Big)^{k(a_1+n/2-1)}.
    \end{align*}

    Since $\abs{\abs{P-S}}_F\geq \epsilon$ for all permutation matrix $P$, there exist $S_{uv}$ which satisfies $S_{uv}^2\in[\dfrac{\epsilon^2}{k},1-\dfrac{\epsilon^2}{k}]$. Then, the inequality holds
    \begin{align}\label{eq:sum_three_term}
        &\sum\limits_{j=1}^k S_{jv}^2 (\hat{\lambda_j}+\dfrac{h}{n}) \nonumber\\
        &= \sum\limits_{j<u}S_{jv}^2 (\hat{\lambda_j}+\dfrac{h}{n}) + S_{uv}^2 (\hat{\lambda_u}+\dfrac{h}{n}) +\sum\limits_{j>u}S_{jv}^2 (\hat{\lambda_j}+\dfrac{h}{n})\nonumber\\
        &\geq(\sum\limits_{j<u}S_{jv}^2) \prod\limits_{j<u} (\hat{\lambda_j}+\dfrac{h}{n})^{\frac{S_{jv}^2}{\sum\limits_{j<u}S_{jv}^2}} +S_{uv}^2 (\hat{\lambda_u}+\dfrac{h}{n}) + (\sum\limits_{j>u}S_{jv}^2) \prod\limits_{j>u} (\hat{\lambda_j}+\dfrac{h}{n})^{\frac{S_{jv}^2}{\sum\limits_{j>u}S_{jv}^2}},
    \end{align}
    where the last inequality holds by arithmetic geometric mean inequality.

    Since $S_{uv}^2\in[\dfrac{\epsilon^2}{k},1-\dfrac{\epsilon^2}{k}]$, $\sum\limits_{j<u}S_{jv}^2$ or $\sum\limits_{j>u}S_{jv}^2$ is larger than $\dfrac{\epsilon^2}{2k}$. Without loss of generality, we assume that $\sum\limits_{j<u}S_{jv}^2$ is larger than $\dfrac{\epsilon^2}{2k}$. Then, the lower bound of \eqref{eq:sum_three_term} is given by
    \begin{align}\label{eq:sum_two_term}
        \eqref{eq:sum_three_term} &\geq (\sum\limits_{j<u}S_{jv}^2) \prod\limits_{j<u} (\hat{\lambda_j}+\dfrac{h}{n})^{\frac{S_{jv}^2}{\sum\limits_{j<u}S_{jv}^2}} +  (\sum\limits_{j\geq u}S_{jv}^2) \prod\limits_{j\geq u} (\hat{\lambda_j}+\dfrac{h}{n})^{\frac{S_{jv}^2}{\sum\limits_{j\geq u}S_{jv}^2}}\nonumber\\
        & = y_1^{t_1}y_2^{1-t_1}\cdot \Big[t_1\big(\dfrac{y_1}{y_2}\big)^{1-t_1}+(1-t_1)\big(\dfrac{y_2}{y_1}\big)^{t_1}\Big]\nonumber\\
        &= \prod\limits_{j=1}^k  (\hat{\lambda_j}+\dfrac{h}{n})^{S_{jv}^2}\cdot \Big[t_1\big(\dfrac{y_1}{y_2}\big)^{1-t_1}+(1-t_1)\big(\dfrac{y_2}{y_1}\big)^{t_1}\Big],
    \end{align}
    where $t_1 = \sum\limits_{j<u}S_{jv}^2$, $y_1= \prod\limits_{j<u} (\hat{\lambda_j}+\dfrac{h}{n})^{\frac{S_{jv}^2}{\sum\limits_{j<u}S_{jv}^2}}$, and $y_2 = \prod\limits_{j\geq u} (\hat{\lambda_j}+\dfrac{h}{n})^{\frac{S_{jv}^2}{\sum\limits_{j\geq u}S_{jv}^2}}$.

    Let $f(t) = t a ^{t - 1} + (1-t) a^{t}$ with $a>1$. Then, $f'(t) = a^{t-1}((\log a - a\log a)t + a\log a-a +1)$ holds. For $t\in (s,1-s)$ with $s>0$, the inequality holds
    $$f(t)\geq f(s) \wedge f(1-s).$$
    We apply talor expansion on $f(s)$ and $f(1-s)$, then we obtain the following equality
    \begin{align*}
        f(s) &= a^s (sa^{-1}+1-s)\\
        & = (1+(\log a) s + \dfrac{1}{2}(\log a)^2 s^2 + O(s^3))(sa^{-1}+1-s)\\
        & = 1 + s(a^{-1}+\log a -1) + O(s^2)\\
        f(1-s) &= a^{-s} (1-s+as)\\
        &=(1-(\log a) s + \dfrac{1}{2}(\log a)^2 s^2 + O(s^3)) (1-s+as)\\
        &= 1 + s (a -1 - \log a) + O(s^2),
    \end{align*}
    for sufficiently small $s$. It implies $f(s)\wedge f(1-s) = f(s) $.

    Therefore, the lower bound of \eqref{eq:sum_two_term} is given by
    \begin{align*}
        \eqref{eq:sum_two_term}&\geq \prod\limits_{j=1}^k(\hat{\lambda_j}+\dfrac{h}{n})^{S_{jv}^2}\cdot\Bigg[1+\dfrac{\epsilon^2}{2k}\Big(\big(\dfrac{y_2}{y_1}\big)^{-1} + \log(\dfrac{y_2}{y_1}) -1 \Big) + O\Big(\big(\dfrac{\epsilon^2}{2k}\big)^2\Big) \Big]\Bigg]\\
        &\geq \prod\limits_{j=1}^k(\hat{\lambda_j}+\dfrac{h}{n})^{S_{jv}^2} \cdot\Bigg[1+ \dfrac{c\epsilon^2}{2k}\log(\dfrac{y_2}{y_1}) +O(\epsilon^4)\Bigg]\\
        &\geq \prod\limits_{j=1}^k(\hat{\lambda_j}+\dfrac{h}{n})^{S_{jv}^2} \cdot\Big[1+ \dfrac{c\epsilon^2}{4k}\min_{i<k}\log(\dfrac{\hat{\lambda}_i}{\hat{\lambda}_{i+1}})\Big],
    \end{align*}
    for some positive constant $c$.

    Finally, the upper bound of \eqref{eq:D/E_2} is given by
    \begin{align*}
        \eqref{eq:D/E_2} &\leq \dfrac{\prod\limits_{i=1}^k \bigg[\dfrac{h}{n}+ \hat{\lambda}_i \bigg]^{a_1+n/2-1} \cdot\Big(1+ \eta^2\dfrac{\hat{\lambda}_1}{\hat{\lambda}_k}\Big)^{k(a_1+n/2-1)}}{\prod\limits_{i=1}^k \bigg[\sum\limits_{j=1}^k S_{ji}^2\Big(\hat{\lambda}_j +\dfrac{h}{n} \Big)\bigg]^{a_1+ n/2-1} \cdot \Big(1 -2k\eta^2 \dfrac{\hat{\lambda}_1 }{\hat{\lambda}_k}\Big)^{k(a_1+n/2-1)}}\\
        &\leq \dfrac{\prod\limits_{i=1}^k \bigg[\dfrac{h}{n}+ \hat{\lambda}_i \bigg]^{a_1+n/2-1} \cdot\Big(1+ \eta^2\dfrac{\hat{\lambda}_1}{\hat{\lambda}_k}\Big)^{k(a_1+n/2-1)}}{\prod\limits_{i=1}^k \bigg[\prod\limits_{j=1}^k \Big(\hat{\lambda}_j +\dfrac{h}{n} \Big)^{S_{ji}^2}\bigg]^{a_1+ n/2-1} \cdot \Big[1+ \dfrac{c\epsilon^2}{4k}\min\limits_{i<k}\log(\dfrac{\hat{\lambda}_i}{\hat{\lambda}_{i+1}})\Big]^{a_1+n/2-1}  \cdot \Big(1 -2k\eta^2 \dfrac{\hat{\lambda}_1 }{\hat{\lambda}_k}\Big)^{k(a_1+n/2-1)}}\\
        &\leq \dfrac{\exp(n\eta^2\dfrac{\hat{\lambda}_1}{\hat{\lambda}_k})}{\exp(n\epsilon^2 \min\limits_{i<k}\log(\dfrac{\hat{\lambda}_i}{\hat{\lambda}_{i+1}}))}.
    \end{align*}
\end{proof}
\hfill\break

\begin{theorem}\label{thm:post_exp_diff_a}
Under the assumptions of Lemma \ref{main-lem:shrink_post_diff_a}, for $i=1,\ldots,k$, we have
\begin{align*}
    \bbE\Big[\dfrac{\lambda_i-\lambda_{0,i}}{{\lambda_{0,i}}}\big\vert \bfX_n \Big]  &=\dfrac{n}{n+2a_i-4}\dfrac{1}{\lambda_{0,i}}\Big[(1+\beta_i)\lambda_{0,i}+\bar{d}\dfrac{p}{n}+\alpha_i\sqrt{\dfrac{p}{n}}+\dfrac{h}{n}\Big]-1\\
    &+O\Big( \epsilon^2\dfrac{\lambda_{0,1}}{\lambda_{0,i}} + \dfrac{\lambda_{0,1}}{\lambda_{0,i}}\exp(-\dfrac{\kappa}{2}\epsilon)+\dfrac{\lambda_{0,1}}{\lambda_{0,i}}\big(\dfrac{n\lambda_{0,k}+p}{p}\big)^{-\epsilon^2 a_n} \Big),
\end{align*}
where $\alpha_i\in[-C,C]$ for some positive constant $C$, $\beta_i\lesssim n^{-1/2+\delta}$ for all small $\delta>0$, $\bar{d} = \dfrac{1}{p-k}\sum\limits_{j=k+1}^p\lambda_{0,j}$, and $\kappa = \min\limits_{l<k}(a_{l+1} - a_l) \cdot \min\limits_{l<k}\log(\dfrac{\lambda_{0,l}}{\lambda_{0,l+1}})$.
\end{theorem}
\hfill\break

\begin{lemma}\label{lem:post_exp_diff_ordered}
Under the assumptions of Lemma \ref{main-lem:shrink_post_diff_a}, for $i=1,\ldots,k$, we have
\begin{align*}
 \dfrac{1}{\lambda_{0,i}}\abs{\bbE[\lambda_i\vert\bfX_n] - \bbE[\lambda_{(i)}\vert\bfX_n]} = \dfrac{\lambda_{0,1}}{\lambda_{0,i}} \cdot\bigg(O\Big(\dfrac{1}{n}\Big) + O\Big( \exp(-\frac{\kappa}{2} \epsilon) \Big) + O\Big( \Big( \frac{n\lambda_{0,k} + p}{p} \Big)^{-\epsilon^2 a_n} \Big)\bigg),
\end{align*}
where $\kappa = \min\limits_{l<k}(a_{l+1} - a_l) \cdot \min\limits_{l<k}\log(\dfrac{\lambda_{0,l}}{\lambda_{0,l+1}})$.
\end{lemma}

\begin{proof}

Let $D_{\epsilon} = \Big\{\Gamma \in O(p): \inf_{Q_2 \in O_{p-k}} \abs{\abs{ \begin{pmatrix} I_k & 0 \\ 0 & Q_2 \end{pmatrix} - \Gamma }}_F < \epsilon \Big\}.$
By Lemma \ref{main-lem:shrink_post_diff_a}, we have
$$\int_{D_\epsilon}\int\pi(\Lambda,\Gamma \vert \bfX_n) (d\Lambda) (d\Gamma) = 1 + O\Big(\exp(-\frac{\kappa}{2} \epsilon)\Big) + O\Big( \Big( \frac{n\lambda_{0,k} + p}{p} \Big)^{-\epsilon^2 a_n} \Big),$$
where $\kappa = \min_{l<k}(a_{l+1} - a_l) \cdot \min_{l<k} \log\Big( \frac{\lambda_{0,l}}{\lambda_{0,l+1}} \Big).$

On $D_\epsilon$, since $(\Gamma_{ii} - 1)^2 + (1 - \Gamma_{ii}^2) < \epsilon^2$, it follows that $\Gamma_{ii}^2 > 1 - \epsilon^2$. Given $\frac{c_i}{n} = \frac{h}{n} + \sum_{j=1}^n \Gamma_{ji}^2 \hat{\lambda}_j,$
we obtain the bounds:
\begin{align*}
\frac{c_i}{n} &\in \left[ (1 - \epsilon^2) \hat{\lambda}_i + \frac{h}{n},\ (1 - \epsilon^2) \hat{\lambda}_i + \epsilon^2 \hat{\lambda}_1 + \frac{h}{n} \right], \quad \text{for } i = 1, \ldots, n, \\
\frac{c_i}{n} &\in \left[ \frac{h}{n},\ \epsilon^2 \hat{\lambda}_1 + \frac{h}{n} \right], \quad \text{for } i = n+1, \ldots, p.
\end{align*}

Under $\Gamma \in D_\epsilon$, each $\lambda_i$ follows an inverse gamma distribution, $\lambda_i \overset{iNd}{\sim} InvGam\Big(a_i + \frac{n}{2} - 1, \frac{c_i}{2} \Big)$. By Chebyshev's inequality,
\begin{align*}
    \pi(\abs{\lambda_i - \bbE[\lambda_i\vert \Gamma,\bfX_n]} > \alpha_i \vert \Gamma,\bfX_n) &\leq \frac{Var(\lambda_i)}{\alpha_i^2} = \frac{2}{\alpha_i^2} \Big( \frac{c_i}{n} \Big)^2 \cdot \frac{n^2}{(2a_i + n - 4)^2 (2a_i + n - 6)}.
\end{align*}

Let $p_i$ be an upper bound for each of these probabilities
\begin{align*}
    p_i = \begin{cases}
        \frac{2}{\alpha_i^2} \Big( (1 - \epsilon^2) \hat{\lambda}_i + \epsilon^2 \hat{\lambda}_1 + \frac{h}{n} \Big)^2 \cdot \frac{n^2}{(2a_i + n - 4)^2 (2a_i + n - 6)} & \quad\text{for }i=1,\ldots,n,\\
        \frac{2}{\alpha_i^2} \Big( \epsilon^2 \hat{\lambda}_1 + \frac{h}{n} \Big)^2 \cdot \frac{n^2}{(2a_i + n - 4)^2 (2a_i + n - 6)} & \quad\text{for }i=n+1,\ldots,p.
    \end{cases}
\end{align*}

Choose $\alpha_i=\dfrac{1}{4}\delta_0 \hat{\lambda}_i$ for $i=1,\ldots,n$, and $\alpha_i = \dfrac{1}{4}\delta_0 \hat{\lambda}_n$ for $i=n+1,\ldots,p$. Then with probability at least $\prod_{i=1}^p (1-p_i)$, the inequalities  
$$\lambda_i \in [\bbE[\lambda_i\vert \Gamma,\bfX_n] - \alpha_i, \bbE[\lambda_i\vert \Gamma,\bfX_n] + \alpha_i],$$
hold for $i = 1,\ldots,p$. This implies:
\begin{align*}
\lambda_1 &> \cdots > \lambda_k, \\
\lambda_k &> \lambda_i \quad \text{for } i = k+1, \ldots, p,
\end{align*}
and therefore,
$$
\lambda_{(i)} = \lambda_i \quad \text{for all } i = 1, \ldots, k.
$$

We now bound the difference of expectation:
\begin{align}\label{eq:diff_exp}
\abs{ \bbE[\lambda_{(i)} - \lambda_i \vert \bfX_n] } 
&= \abs{ \bbE[(\lambda_{(i)} - \lambda_i) \cdot I(\lambda_{(i)} \neq \lambda_i) \vert \bfX_n] } \nonumber\\
&\leq 2  \bbE[\sum_{j=1}^p\lambda_j \cdot I(\lambda_{(i)} \neq \lambda_i) \vert \bfX_n] \nonumber\\
&= 2  \bbE[ \sum_{j=1}^p\bbE[\lambda_j \cdot I(\lambda_{(i)} \neq \lambda_i) \vert \Gamma,\bfX_n] \cdot (I(\Gamma \in D_\epsilon) + I(\Gamma \notin D_\epsilon) ) \vert \bfX_n ] \nonumber\\
&\leq 2  \bbE[ \sum_{j=1}^p\bbE[\lambda_j \cdot I(\lambda_{(i)} \neq \lambda_i) \vert \Gamma,\bfX_n] \cdot I(\Gamma \in D_\epsilon) \vert \bfX_n ]  + 2  \bbE[ \sum_{j=1}^p\bbE[\lambda_j \vert \Gamma,\bfX_n] \cdot I(\Gamma \notin D_\epsilon) \vert \bfX_n ].
\end{align}

There exists $\beta_j(\Gamma) > 0$ such that
$$\pi(\lambda_{(i)} \neq \lambda_i \vert \Gamma, \bfX_n) = \pi(\lambda_j \geq \beta_j(\Gamma) \vert\Gamma, \bfX_n),$$
for each $j=1,\ldots,p$. Therefore, the inequality holds
$$\bbE[ \sum_{j=1}^p\bbE[\lambda_j \cdot I(\lambda_{(i)} \neq \lambda_i) \vert \Gamma, \bfX_n] \cdot I(\Gamma \in D_\epsilon) \vert \bfX_n ]\leq \bbE[ \sum_{j=1}^p\bbE[\lambda_j \cdot I(\lambda_j \geq \beta_j(\Gamma)) \vert \Gamma, \bfX_n] \cdot I(\Gamma \in D_\epsilon) \vert \bfX_n ].$$

By applying Hölder's inequality, we obtain the following bound on the first term on the right-hand side of \eqref{eq:diff_exp} as follows
\begin{align*}
&\bbE[ \sum_{j=1}^p\bbE[\lambda_j \cdot I(\lambda_j \geq \beta_j(\Gamma)) \vert \Gamma, \bfX_n] \cdot I(\Gamma \in D_\epsilon) \vert \bfX_n ] \\
&\leq  \bbE[ \sum_{j=1}^p\sqrt{\bbE[\lambda_j^2\vert \Gamma]}\cdot\sqrt{\bbE[I(\lambda_j>\beta_j(\Gamma))\vert \Gamma, \bfX_n]} I(\Gamma\in D_\epsilon)\vert \bfX_n]\\
&\leq  \bbE[\sum_{j=1}^p\bbE[\lambda_j\vert \Gamma]\cdot\sqrt{\pi(\lambda_j>\beta_j(\Gamma)\vert \Gamma, \bfX_n)} I(\Gamma\in D_\epsilon)\vert \bfX_n]\\
&\leq  \sup_{\Gamma \in D_\epsilon} \Big(\sum_{j=1}^p\dfrac{c_j}{n+2a_j-4} \Big) \cdot \sqrt{1-\prod_{i=1}^p(1-p_i)}\\
&\leq  \sup_{\Gamma \in D_\epsilon} \Big(\sum_{j=1}^p\dfrac{c_j}{n+2a_j-4}  \Big)\cdot \sqrt{\sum_{i=1}^p p_i}.
\end{align*}

We obtain the following bound on the second term on the right-hand side of \eqref{eq:diff_exp} as follows
\begin{align*}
     &\bbE[ \sum_{j=1}^p\bbE[\lambda_j \vert \Gamma,\bfX_n] \cdot I(\Gamma \notin D_\epsilon) \vert \bfX_n ]\\
     &=\bbE[\Big(\sum_{j=1}^p\dfrac{c_j}{n+2a_j-4} \Big) \cdot I(\Gamma \notin D_\epsilon) \vert \bfX_n]\\
    &\leq \sup_{\Gamma \notin D_\epsilon}\Big(\sum_{j=1}^p  \frac{c_j}{n + 2a_j - 4}\Big)
\cdot \Big( O\Big( \exp(-\frac{\kappa}{2} \epsilon) \Big) + O\Big( \Big( \frac{n\lambda_{0,k} + p}{p} \Big)^{-\epsilon^2 a_n} \Big) \Big).
\end{align*}

Therefore, we obtain
\begin{align*}
     \abs{\bbE[\lambda_i\vert\bfX_n] - \bbE[\lambda_{(i)}\vert\bfX_n]} &\leq   \sup_{\Gamma \in D_\epsilon} \Big(\sum_{j=1}^p\dfrac{c_j}{n+2a_j-4}\Big)  \cdot \sum_{i=1}^p p_i \\
     &+ \sup_{\Gamma \notin D_\epsilon} \Big(\sum_{j=1}^p\dfrac{c_j}{n+2a_j-4}\Big)
\cdot \Big( O\Big( \exp(-\frac{\kappa}{2} \epsilon) \Big) + O\Big( \Big( \frac{n\lambda_{0,k} + p}{p} \Big)^{-\epsilon^2 a_n} \Big) \Big)\\
    &\leq   \sup_{\Gamma } \Big(\sum_{j=1}^p\dfrac{c_j}{n+2a_j-4}\Big)  \cdot \sum_{i=1}^p p_i \\
     &+ \sup_{\Gamma} \Big(\sum_{j=1}^p\dfrac{c_j}{n+2a_j-4}\Big)
\cdot \Big( O\Big( \exp(-\frac{\kappa}{2} \epsilon) \Big) + O\Big( \Big( \frac{n\lambda_{0,k} + p}{p} \Big)^{-\epsilon^2 a_n} \Big) \Big).
\end{align*}

Since $a_1\leq\cdots\leq a_p$, the term $\sum_{j=1}^p\dfrac{c_j}{n+2a_j-4}$ is maximized when $\Gamma = I_p$, and hence
\begin{align*}
    \sup_{\Gamma} \Big(\sum_{j=1}^p\dfrac{c_j}{n+2a_j-4}\Big) &\leq \sum_{j=1}^n\dfrac{n\hat{\lambda}_j + h}{n+2a_j-4} + \sum_{j=n+1}^p\dfrac{h}{n+2a_j-4}\\
    &= O(  \sum_{j=1}^k\lambda_{0,j}  +  \dfrac{np}{a_n} + \dfrac{hp}{a_p})\\
    &= O( \sum_{j=1}^k\lambda_{0,j} ),\\
    &= O(\lambda_{0,1}),
\end{align*}
where we used the condition $a_n \succ n^{3/2}p$.

Next, we evaluate the upper bound of $\sum_{i=1}^p p_i$
\begin{align*}
    \sum_{i=1}^p p_i & = \sum_{i=1}^n \frac{2}{\alpha_i^2} \Big( (1 - \epsilon^2) \hat{\lambda}_i + \epsilon^2 \hat{\lambda}_1 + \frac{h}{n} \Big)^2 \cdot \frac{n^2}{(2a_i + n - 4)^2 (2a_i + n - 6)}\\
    & + \sum_{i=n+1}^p \frac{2}{\alpha_i^2} \Big(  \epsilon^2 \hat{\lambda}_1 + \frac{h}{n} \Big)^2 \cdot \frac{n^2}{(2a_i + n - 4)^2 (2a_i + n - 6)} \\
    &\asymp \sum_{i=1}^k \dfrac{n^2}{(2a_i + n - 4)^2 (2a_i + n - 6)} +\sum_{i=k+1}^n \dfrac{n^2}{(2a_i + n - 4)^2 (2a_i + n - 6)}+ \dfrac{h^2}{n^2\hat{\lambda}_n^2}\sum_{i=n+1}^p \dfrac{n^2}{(2a_i + n - 4)^2 (2a_i + n - 6)}\\
    & =  O( \dfrac{k}{n}) + O(\dfrac{n^2(n-k)}{a_n^3}) + O(\dfrac{h^2}{p^2}\cdot\dfrac{n^2(p-n)}{a_p^3})\\
    &= O( \dfrac{1}{n}),
\end{align*}
again using $a_n \succ n^{3/2}p$.

Combining these, we conclude
\begin{align*}
    &\dfrac{1}{\lambda_{0,i}}\abs{\bbE[\lambda_i\vert\bfX_n] - \bbE[\lambda_{(i)}\vert\bfX_n]}\\
    &\leq \dfrac{1}{\lambda_{0,i}}\cdot\Big(  \sup_{\Gamma} \Big(\sum_{j=1}^p\dfrac{c_j}{n+2a_j-4}\Big) \cdot \Big(\sum_{i=1}^p p_i  + O\Big( \exp(-\frac{\kappa}{2} \epsilon) \Big) + O\Big( \Big( \frac{n\lambda_{0,k} + p}{p} \Big)^{-\epsilon^2 a_n} \Big)\Big)\\
    &= \dfrac{\lambda_{0,1}}{\lambda_{0,i}} \cdot\bigg(O\Big(\dfrac{1}{n}\Big) + O\Big( \exp(-\frac{\kappa}{2} \epsilon) \Big) + O\Big( \Big( \frac{n\lambda_{0,k} + p}{p} \Big)^{-\epsilon^2 a_n} \Big)\bigg).
\end{align*}

\end{proof}


\hfill\break

\begin{theorem}\label{thm:post_exp_equal_a}
 Under the assumptions of Lemma \ref{main-lem:shrink_post_equal_a}, for $i=1,\ldots,k$, we have
\begin{align*}
    \bbE\Big[\dfrac{\lambda_i-\lambda_{0,i}}{{\lambda_{0,i}}}\big\vert \bfX_n \Big]  &=\dfrac{n}{n+2a_i-4}\dfrac{1}{\lambda_{0,i}}\Big[(1+\beta_i)\lambda_{0,i}+\bar{d}\dfrac{p}{n}+\alpha_i\sqrt{\dfrac{p}{n}}+\dfrac{h}{n}\Big]-1\\
    &+O\Big( \epsilon^2\dfrac{\lambda_{0,1}}{\lambda_{0,i}} + \dfrac{\lambda_{0,1}}{\lambda_{0,i}}\exp(-\dfrac{\tau}{2}n\epsilon^2)\cdot\Big(\dfrac{2\lambda_{0,1}}{\lambda_{0,k}}\Big)^{kq}+\dfrac{\lambda_{0,1}}{\lambda_{0,i}}\big(\dfrac{n\lambda_{0,k}+p}{p}\big)^{-\epsilon^2 a_n} \Big),
\end{align*}
where $\alpha_i\in[-C,C]$ for some positive constant $C$, $\beta_i\lesssim n^{-1/2+\delta}$ for all small $\delta>0$, $\bar{d} = \dfrac{1}{p-k}\sum\limits_{j=k+1}^p\lambda_{0,j}$, and $\tau = \min\limits_{l<k}\log(\dfrac{\lambda_{0,l}}{\lambda_{0,l+1}})$.    
\end{theorem}
\hfill\break

\begin{lemma}\label{lem:post_exp_equal_ordered}
Under the assumptions of Lemma \ref{main-lem:shrink_post_equal_a}, for $i=1,\ldots,k$, we have
\begin{align*}
 \dfrac{1}{\lambda_{0,i}}\abs{\bbE[\lambda_i\vert\bfX_n] - \bbE[\lambda_{(i)}\vert\bfX_n]} = \dfrac{\lambda_{0,1}}{\lambda_{0,i}} \cdot\bigg(O\Big(\dfrac{1}{n}\Big) + O\Big( \exp(-\frac{\tau}{2}n \epsilon^2)\cdot\Big(\dfrac{2\lambda_{0,1}}{\lambda_{0,k}}\Big)^{kq} \Big) + O\Big( \Big( \frac{n\lambda_{0,k} + p}{p} \Big)^{-\epsilon^2 a_n} \Big)\bigg),
\end{align*}
where $\tau = \min\limits_{l<k}\log(\dfrac{\lambda_{0,l}}{\lambda_{0,l+1}})$.
\end{lemma}

\begin{proof}

The subset $D_{\epsilon,l}$ of the orthogonal group $O(p)$ is defined as:
    \begin{align*}
        D_{\epsilon,l} = \bigg\{\Gamma\in O(p):\displaystyle\inf_{Q_2\in O_{p-k}} \abs{\abs{\begin{pmatrix}
           P_l & 0 \\ 
        0& Q_2
        \end{pmatrix}-\Gamma}}_F <\epsilon\bigg\},
    \end{align*}
    where $P_1,\ldots,P_L$ denotes all $k\times k$ permutation matrices.

The proof proceeds by applying the same argument as in Lemma \ref{lem:post_exp_diff_ordered},  
with $D_{\epsilon}$ replaced by the union 
$$\bigcup_{l=1}^L D_{\epsilon,l} .$$

By Lemma \ref{main-lem:shrink_post_diff_a}, we have
$$\int_{\bigcup_{l=1}^L D_{\epsilon,l}}\int\pi(\Lambda,\Gamma \vert \bfX_n) (d\Lambda) (d\Gamma) = 1 + O\Big(\exp(-\frac{\tau}{2}n \epsilon^2)\cdot\Big(\dfrac{2\lambda_{0,1}}{\lambda_{0,k}}\Big)^{kq}\Big) + O\Big( \Big( \frac{n\lambda_{0,k} + p}{p} \Big)^{-\epsilon^2 a_n} \Big),$$
where $\tau = \min\limits_{l<k}\log(\dfrac{\lambda_{0,l}}{\lambda_{0,l+1}}).$

Following the same bounding steps as in Lemma \ref{lem:post_exp_diff_ordered}, we obtain
\begin{align*}
 \dfrac{1}{\lambda_{0,i}}\abs{\bbE[\lambda_i\vert\bfX_n] - \bbE[\lambda_{(i)}\vert\bfX_n]} = \dfrac{\lambda_{0,1}}{\lambda_{0,i}} \cdot\bigg(O\Big(\dfrac{1}{n}\Big) + O\Big( \exp(-\frac{\tau}{2}n \epsilon^2)\cdot\Big(\dfrac{2\lambda_{0,1}}{\lambda_{0,k}}\Big)^{kq} \Big) + O\Big( \Big( \frac{n\lambda_{0,k} + p}{p} \Big)^{-\epsilon^2 a_n} \Big)\bigg).
\end{align*}

\end{proof}

\subsection*{Proof of Theorem}

\subsubsection*{Proof of Lemma \ref{main-lem:shrink_post_diff_a}}
\begin{proof}

We consider the following inequality:
\begin{align}
    &\pi\bigg(\bigg\{\Gamma\in O(p):\displaystyle\inf_{Q_1\in O(k),\;Q_2\in O(p-k)} \abs{\abs{\begin{pmatrix}
       Q_1 & 0 \\ 
        0& Q_2
    \end{pmatrix}-\Gamma}}_F <\epsilon\bigg\}\vert \bfX_n\bigg)\label{ineq:proof_lem3.1_0}\\
    &\leq\pi\bigg(\bigg\{\Gamma\in O(p):\displaystyle\inf_{Q_2\in O(p-k)} \abs{\abs{\begin{pmatrix}
       I_k & 0 \\ 
        0& Q_2
    \end{pmatrix}-\Gamma}}_F <\epsilon_4\bigg\}\vert \bfX_n\bigg)\label{ineq:proof_lem3.1_1}\\
    &+ \pi\bigg(\bigg\{\Gamma\in O(p):\displaystyle\inf_{Q_1\in O(k),\;Q_2\in O(p-k)} \abs{\abs{\begin{pmatrix}
       Q_1 & 0 \\ 
        0& Q_2
    \end{pmatrix}-\Gamma}}_F <\epsilon,\quad\displaystyle\inf_{Q_2\in O(p-k)} \abs{\abs{\begin{pmatrix}
       I_k & 0 \\ 
        0& Q_2
    \end{pmatrix}-\Gamma}}_F \geq\epsilon_4\bigg\} \vert \bfX_n\bigg),\label{ineq:proof_lem3.1_2}
\end{align}
for some positive $\epsilon_4$. By Lemma \ref{lem:post_contract}, we have derived the following lower bound on the posterior probability:
\begin{align*}
    \eqref{ineq:proof_lem3.1_0} \geq 1-\big(\dfrac{n\lambda_{0,k}+p}{p}\big)^{-\epsilon^2 a_n}.
\end{align*}
Thus, in order to establish the convergence of \eqref{ineq:proof_lem3.1_1} to $1$, it suffices to show that \eqref{ineq:proof_lem3.1_2} vanishes as $n\rightarrow\infty$.

Now, we focus on the term:
$$\pi\bigg(\bigg\{\Gamma\in O(p):\displaystyle\inf_{Q_1\in O(k),\;Q_2\in O(p-k)} \abs{\abs{\begin{pmatrix}
       Q_1 & 0 \\ 
        0& Q_2
    \end{pmatrix}-\Gamma}}_F <\epsilon,\quad\displaystyle\inf_{Q_2\in O(p-k)} \abs{\abs{\begin{pmatrix}
       I_k & 0 \\ 
        0& Q_2
    \end{pmatrix}-\Gamma}}_F \geq\epsilon_4\bigg\} \vert \bfX_n\bigg).$$

Suppose that $\{R_1,\ldots, R_{M(O(k),\abs{\abs{\cdot}}_F,\epsilon_5)}\}$ forms a maximal $\epsilon_5$-packing of $O(k)$, where $R_{M(O(k),\abs{\abs{\cdot}}_F,\epsilon_5)}=I_k$, which also serves as an $\epsilon_5$-covering. For all $Q_1\in O(k)$, there exists some $R_i$ such that $\abs{\abs{R_i - Q_1}}_F<\epsilon_5$. By triangular inequality, we obtain 
\begin{align}\label{ineq:tri_lem3.1}
     \abs{\abs{\begin{pmatrix}
        R_i & 0 \\ 
       0& Q_2
    \end{pmatrix}-\Gamma}}_F \leq \abs{\abs{\begin{pmatrix}
        R_i & 0 \\ 
       0& Q_2
    \end{pmatrix}-\begin{pmatrix}
       Q_1 & 0 \\ 
        0& Q_2
    \end{pmatrix}}}_F + \abs{\abs{\begin{pmatrix}
       Q_1 & 0 \\ 
        0& Q_2
    \end{pmatrix}-\Gamma}}_F.
\end{align}

Taking the infimum over $Q_2\in O(p-k)$ on both sides of \eqref{ineq:tri_lem3.1}, then the inequality holds
\begin{align}\label{ineq:inf_tri_lem3.1}
    \inf_{Q_2\in O(p-k)}\abs{\abs{\begin{pmatrix}
        R_i & 0 \\ 
       0& Q_2
    \end{pmatrix}-\Gamma}}_F &\leq \inf_{Q_2\in O(p-k)}\abs{\abs{\begin{pmatrix}
        R_i & 0 \\ 
       0& Q_2
    \end{pmatrix}-\begin{pmatrix}
       Q_1 & 0 \\ 
        0& Q_2
    \end{pmatrix}}}_F + \inf_{Q_2\in O(p-k)}\abs{\abs{\begin{pmatrix}
       Q_1 & 0 \\ 
        0& Q_2
    \end{pmatrix}-\Gamma}}_F\nonumber\\
    &\leq \epsilon_5 + \inf_{Q_2\in O(p-k)}\abs{\abs{\begin{pmatrix}
       Q_1 & 0 \\ 
        0& Q_2
    \end{pmatrix}-\Gamma}}_F.
\end{align}

Let $\{S_1,\ldots,S_m\}$ be a subset of $\{R_1,\ldots, R_{M(O(k),\abs{\abs{\cdot}}_F,\epsilon_5)}\}$ which satisfies the following condition
\begin{equation}\label{eq:cond_S}
    \bigg\{\Gamma\in O(p):\displaystyle\inf_{Q_2\in O(p-k)} \abs{\abs{\begin{pmatrix}
       S_i & 0 \\ 
        0& Q_2
    \end{pmatrix}-\Gamma}}_F <\epsilon+\epsilon_5,\quad \displaystyle\inf_{Q_2\in O(p-k)} \abs{\abs{\begin{pmatrix}
       I_k & 0 \\ 
        0& Q_2
    \end{pmatrix}-\Gamma}}_F \geq\epsilon_4\bigg\} \neq \phi.
\end{equation}
Moreover, it follows that $I_k\notin \{S_1,\ldots,S_m\}$, where $\epsilon_4>\epsilon+\epsilon_5$.

Then, the following inequality holds:
\begin{align*}
    &\pi\bigg(\bigg\{\Gamma\in O(p):\displaystyle\inf_{Q_1\in O(k),\;Q_2\in O(p-k)} \abs{\abs{\begin{pmatrix}
       Q_1 & 0 \\ 
        0& Q_2
    \end{pmatrix}-\Gamma}}_F <\epsilon,\quad\displaystyle\inf_{Q_2\in O(p-k)} \abs{\abs{\begin{pmatrix}
       I_k & 0 \\ 
        0& Q_2
    \end{pmatrix}-\Gamma}}_F \geq\epsilon_4\bigg\} \vert \bfX_n\bigg)\\
    &\leq\pi\bigg(\bigcup_{i=1}^{M(O(k),\abs{\abs{\cdot}}_F,\epsilon_5)}\bigg\{\Gamma\in O(p):\displaystyle\inf_{Q_2\in O(p-k)} \abs{\abs{\begin{pmatrix}
       R_i & 0 \\ 
        0& Q_2
    \end{pmatrix}-\Gamma}}_F <\epsilon+\epsilon_5,\quad \displaystyle\inf_{Q_2\in O(p-k)} \abs{\abs{\begin{pmatrix}
       I_k & 0 \\ 
        0& Q_2
    \end{pmatrix}-\Gamma}}_F \geq\epsilon_4\bigg\} \vert \bfX_n\bigg)\\
    &\leq \pi\bigg(\bigcup_{i=1}^{m}\bigg\{\Gamma\in O(p):\displaystyle\inf_{Q_2\in O(p-k)} \abs{\abs{\begin{pmatrix}
       S_i & 0 \\ 
        0& Q_2
    \end{pmatrix}-\Gamma}}_F <\epsilon+\epsilon_5\bigg\} \vert \bfX_n\bigg),
\end{align*}
for some positive $\epsilon_5$ which satisfies $\epsilon_4>\epsilon+\epsilon_5$. The first inequality follows from \eqref{ineq:inf_tri_lem3.1}, and the second inequality follows from \eqref{eq:cond_S}.

By the triangular inequality, the following inequality holds:
\begin{align*}
    \abs{\abs{\begin{pmatrix}
       S_i & 0 \\ 
        0& Q_2
    \end{pmatrix}-\begin{pmatrix}
       I_k & 0 \\ 
        0& Q_2
    \end{pmatrix}}}_F + \abs{\abs{\begin{pmatrix}
       I_k & 0 \\ 
        0& Q_2
    \end{pmatrix}-\Gamma}}_F \geq  \abs{\abs{\begin{pmatrix}
       S_i & 0 \\ 
        0& Q_2
    \end{pmatrix}-\Gamma}}_F ,
\end{align*}
for $i=1,\ldots,m$. Taking the infimum over $Q_2\in O(p-k)$ on both sides yields
\begin{align*}
    \abs{\abs{S_i-I_k}}_F + \displaystyle\inf_{Q_2\in O(p-k)}\abs{\abs{\begin{pmatrix}
       I_k & 0 \\ 
        0& Q_2
    \end{pmatrix}-\Gamma}}_F \geq  \displaystyle\inf_{Q_2\in O(p-k)}\abs{\abs{\begin{pmatrix}
       S_i & 0 \\ 
        0& Q_2
    \end{pmatrix}-\Gamma}}_F ,
\end{align*}
for $i=1,\ldots,m$. By the definition of $S_i$, there exists $\Gamma_i\in O(p)$ satisfying condition \eqref{eq:cond_S}. Therefore, we obtain the inequality
$$\abs{\abs{S_i- I_k}}_F>\epsilon_4-(\epsilon+\epsilon_5),$$
for $i=1,\ldots,m$.

Next, we define the sets:
\begin{align*}
    D_\eta  &= \bigg\{\Gamma\in O(p):\displaystyle\inf_{Q_2\in O(p-k)} \abs{\abs{\begin{pmatrix}
       I_k & 0 \\ 
        0& Q_2
    \end{pmatrix}-\Gamma}}_F <\eta\bigg\}, \\
    E_{\eta,l} &= \bigg\{\Gamma\in O(p):\displaystyle\inf_{Q_2\in O(p-k)} \abs{\abs{\begin{pmatrix}
       S_l & 0 \\ 
        0& Q_2
    \end{pmatrix}-\Gamma}}_F <\eta\bigg\} ,
\end{align*}
where $S_l\in O(k)$ and $\abs{\abs{S_l-I_k}}_F>\epsilon_4-(\epsilon+\epsilon_5)$.

By considering the variable transformation $$\Gamma^{*} = \Gamma\begin{pmatrix}
       S_l^T & 0 \\ 
        0& I_{p-k}
    \end{pmatrix},$$
it follows that $(d\Gamma) = (d\Gamma^*)$, since $(d\Gamma)$ is an invariant measure.
Therefore, we obtain the following equality:
\begin{align*}
    \int_{E_{\eta,l}} \prod\limits_{i=k+1}^p\Big(\dfrac{c_i}{n}\Big)^{-a_i-n/2+1}(d\Gamma) &= \int_{E_{\eta,l}} \prod\limits_{i=k+1}^p \Big(\dfrac{h}{n}+\sum\limits_{j=1}^n\Gamma_{ji}^2\hat{\lambda}_j\Big)^{-a_i-n/2+1}(d\Gamma)\\
    &= \int_{D_\eta} \prod\limits_{i=k+1}^p \Big(\dfrac{h}{n}+\sum\limits_{j=1}^n(\Gamma^{*}_{ji})^2\hat{\lambda}_j\Big)^{-a_i-n/2+1}(d\Gamma^*)\\
    &= \int_{D_\eta} \prod\limits_{i=k+1}^p\Big(\dfrac{c_i}{n}\Big)^{-a_i-n/2+1}(d\Gamma).
\end{align*}
The last equality holds due to 
$$\abs{\abs{\begin{pmatrix}
       S_l & 0 \\ 
        0& Q_2
    \end{pmatrix}-\Gamma}}_F = \abs{\abs{\begin{pmatrix}
       I_k & 0 \\ 
        0& Q_2
    \end{pmatrix}-\Gamma^*}}_F,$$
and $$\Gamma_{ji}^2 = (\Gamma^{*}_{ji})^2 \quad \text{for }i>k.$$

Therefore, we obtain the following inequality:
\begin{align}\label{eq:D/E}
    \dfrac{\displaystyle\int_{E_{\epsilon+\epsilon_5,l}}\prod\limits_{i=1}^p c_i^{-a_i-n/2+1}(d\Gamma)}{\displaystyle\int_{D_{\epsilon+\epsilon_5}}\prod\limits_{i=1}^p c_i^{-a_i-n/2+1}(d\Gamma)}
    &\leq \dfrac{\sup\limits_{E_{\epsilon+\epsilon_5,l}}\prod\limits_{i=1}^k c_i^{-a_i-n/2+1} \cdot \displaystyle\int_{E_{\epsilon+\epsilon_5,l}}\prod\limits_{i=k+1}^p c_i^{-a_i-n/2+1}(d\Gamma) }{\inf\limits_{D_{\epsilon+\epsilon_5}}\prod\limits_{i=1}^k c_i^{-a_i-n/2+1} \cdot \displaystyle\int_{D_{\epsilon+\epsilon_5}}\prod\limits_{i=k+1}^p c_i^{-a_i-n/2+1}(d\Gamma)}\nonumber\\
    &= \dfrac{\sup\limits_{D_{\epsilon+\epsilon_5}}\prod\limits_{i=1}^k \Big(\dfrac{c_i}{n}\Big)^{a_i+n/2-1} }{\inf\limits_{E_{\epsilon+\epsilon_5,l}}\prod\limits_{i=1}^k \Big(\dfrac{c_i}{n}\Big)^{a_i+n/2-1} }.
\end{align}

By applying Lemma \eqref{lem:upp_bound_D/E}, the upper bound of \eqref{eq:D/E} is given by
\begin{align*}
    \eqref{eq:D/E} &\preccurlyeq \dfrac{\exp(n(\epsilon+\epsilon_5)^2 \dfrac{\hat{\lambda}_1 }{\hat{\lambda}_k})}{\exp(\dfrac{\epsilon_4-(\epsilon+\epsilon_5)}{2\sqrt{k}}\min\limits_{l<k}(a_{l+1} - a_l) \cdot \min\limits_{l<k}\log(\dfrac{\hat{\lambda}_l}{\hat{\lambda}_{l+1}}))}\\
    &\asymp \dfrac{\exp(n\epsilon^2 \dfrac{\lambda_{0,1} }{\lambda_{0,k}})}{\exp(\epsilon_4\min\limits_{l<k}(a_{l+1} - a_l) \cdot \min\limits_{l<k}\log(\dfrac{\lambda_{0,l}}{\lambda_{0,l+1}}))},
\end{align*}
where $\epsilon_4 \succ \epsilon+\epsilon_5$ and $\epsilon\asymp \epsilon_5$.

Let $\kappa = \min\limits_{l<k}(a_{l+1} - a_l) \cdot \min\limits_{l<k}\log(\dfrac{\lambda_{0,l}}{\lambda_{0,l+1}})$. Then, we get the following inequality
\begin{align*}
    &\pi\bigg(\bigg\{\Gamma\in O(p):\displaystyle\inf_{Q_1\in O_{k},\;Q_2\in O_{p-k}} \abs{\abs{\begin{pmatrix}
       Q_1 & 0 \\ 
        0& Q_2
    \end{pmatrix}-\Gamma}}_F <\epsilon\bigg\}\vert \bfX_n\bigg)\\
    &\leq \pi(D_{\epsilon_4}\vert \bfX_n) + m \pi(E_{\epsilon+\epsilon_5,1}\vert \bfX_n)\\
    &\leq \pi(D_{\epsilon_4}\vert \bfX_n) + M (O(k),\abs{\abs{\cdot}}_F ,\epsilon_5) \cdot \dfrac{\exp(n\epsilon^2 \dfrac{\lambda_{0,1} }{\lambda_{0,k}})}{\exp(\kappa\epsilon_4)} \cdot \pi(D_{\epsilon+\epsilon_5}\vert \bfX_n) \\
    &\leq \Bigg(1+ M (O(k),\abs{\abs{\cdot}}_F ,\epsilon_5) \cdot \dfrac{\exp(n\epsilon^2 \dfrac{\lambda_{0,1} }{\lambda_{0,k}})}{\exp(\kappa\epsilon_4)} \Bigg)\pi(D_{\epsilon_4}\vert \bfX_n).
\end{align*}

We obtain the following probability 
\begin{align*}
    \pi(D_{\epsilon_4}\vert \bfX_n) &\geq \Bigg(1+ M (O(k),\abs{\abs{\cdot}}_F ,\epsilon_5) \cdot \dfrac{\exp(n\epsilon^2 \dfrac{\lambda_{0,1} }{\lambda_{0,k}})}{\exp(\kappa\epsilon_4)} \Bigg) ^{-1} (1-\big(\dfrac{n\lambda_{0,k}+p}{p}\big)^{-\epsilon_4^2 a_n})\\
    &\geq \Bigg(1+  \exp( n\epsilon^2 \dfrac{\lambda_{0,1} }{\lambda_{0,k}} -\kappa \epsilon_4  )\Bigg) ^{-1} (1-\big(\dfrac{n\lambda_{0,k}+p}{p}\big)^{-\epsilon_4^2 a_n})\\
    &\geq \Bigg(1+  \exp(-\dfrac{\kappa}{2}\epsilon_4)\Bigg) ^{-1} (1-\big(\dfrac{n\lambda_{0,k}+p}{p}\big)^{-\epsilon_4^2 a_n}),
\end{align*}
where $ n\epsilon^2 \dfrac{\lambda_{0,1} }{\lambda_{0,k}}\leq \dfrac{\kappa}{2}\epsilon_4 $.

Therefore, we obtain:
\begin{align*}
    \pi(D_{\epsilon_4}^c\vert \bfX_n) &\leq 1 - \Bigg(1+  \exp(-\dfrac{\kappa}{2}\epsilon_4)\Bigg) ^{-1} (1-\big(\dfrac{n\lambda_{0,k}+p}{p}\big)^{-\epsilon_4^2 a_n})\\
    & = \dfrac{\exp(-\dfrac{\kappa}{2}\epsilon_4)+\big(\dfrac{n\lambda_{0,k}+p}{p}\big)^{-\epsilon_4^2 a_n} }{1+  \exp(-\dfrac{\kappa}{2}\epsilon_4)}\\
    &= O\Big(\exp(-\dfrac{\kappa}{2}\epsilon_4) \Big)  + O\Big(\big(\dfrac{n\lambda_{0,k}+p}{p}\big)^{-\epsilon_4^2 a_n}\Big),
\end{align*}
where $\epsilon_4 \succ \kappa^{-1}$.

In summary, if there exists $\epsilon>0$ satisfying
$$\dfrac{np}{a_n}\prec \epsilon^2 \leq \dfrac{\lambda_{0,k}}{n\lambda_{0,1}}\dfrac{\kappa}{2}\epsilon_4,$$ 
and 
$$\epsilon_4\succ \epsilon\vee \kappa^{-1},$$
then the following holds
$$\pi(D_{\epsilon_4}^c\vert \bfX_n)= O\Big(\exp(-\dfrac{\kappa}{2}\epsilon_4) \Big)  + O\Big(\big(\dfrac{n\lambda_{0,k}+p}{p}\big)^{-\epsilon_4^2 a_n}\Big).$$

In particular, the same posterior bound holds under the simplified condition 
$\dfrac{np}{a_n} \prec \dfrac{\lambda_{0,k}}{n \lambda_{0,1}} \cdot \dfrac{\kappa}{2} \epsilon_4$ 
and $\epsilon_4 \succ \sqrt{\dfrac{np}{a_n}} \vee \kappa^{-1}$.

\end{proof}

\subsubsection*{Proof of Lemma \ref{main-lem:shrink_post_equal_a}}
\begin{proof}
The proof follows the same procedure as the proof of Lemma \ref{main-lem:shrink_post_diff_a}. By replacing $$\displaystyle\inf_{Q_2\in O(p-k)} \abs{\abs{\begin{pmatrix}
       I_k & 0 \\ 
        0& Q_2
    \end{pmatrix}-\Gamma}}_F \geq\epsilon_4,$$
with 
$$\displaystyle\inf_{Q_2\in O(p-k)} \abs{\abs{\begin{pmatrix}
       P & 0 \\ 
        0& Q_2
    \end{pmatrix}-\Gamma}}_F \geq\epsilon_4,\quad \text{for all permutation matrix }P,$$
we obtain the following.

Suppose that $\{R_1,\ldots, R_{M(O(k),\abs{\abs{\cdot}}_F,\epsilon_5)}\}$ forms a maximal $\epsilon_5$-packing of $O(k)$, where $R_{M(O(k),\abs{\abs{\cdot}}_F,\epsilon_5)}=I_k$, which also serves as an $\epsilon_5$-covering. Let $\{S_1,\ldots,S_m\}$ be a subset of $\{R_1,\ldots, R_{M(O(k),\abs{\abs{\cdot}}_F,\epsilon_5)}\}$ which satisfies the following condition
\begin{equation}\label{eq:cond_S_2}
    \bigg\{\Gamma\in O(p):\displaystyle\inf_{Q_2\in O(p-k)} \abs{\abs{\begin{pmatrix}
       S_i & 0 \\ 
        0& Q_2
    \end{pmatrix}-\Gamma}}_F <\epsilon+\epsilon_5,\; \displaystyle\inf_{Q_2\in O(p-k)} \abs{\abs{\begin{pmatrix}
       P_l & 0 \\ 
        0& Q_2
    \end{pmatrix}-\Gamma}}_F \geq\epsilon_4,\; l=1,\ldots,L\bigg\} \neq \phi.
\end{equation}
Moreover, it follows that $P_l\notin \{S_1,\ldots,S_m\}$ for $l=1,\ldots,L$, where $\epsilon_4>\epsilon+\epsilon_5$.

Next, we define the sets:
\begin{align*}
    D_{\eta,l}  &= \bigg\{\Gamma\in O(p):\displaystyle\inf_{Q_2\in O(p-k)} \abs{\abs{\begin{pmatrix}
       P_l & 0 \\ 
        0& Q_2
    \end{pmatrix}-\Gamma}}_F <\eta\bigg\} ,\\
    E_{\eta,r} &= \bigg\{\Gamma\in O(p):\displaystyle\inf_{Q_2\in O(p-k)} \abs{\abs{\begin{pmatrix}
       S_r & 0 \\ 
        0& Q_2
    \end{pmatrix}-\Gamma}}_F <\eta\bigg\} ,
\end{align*}
where $S_r\in O(k)$ and $\abs{\abs{S_r-P_l}}_F>\epsilon_4-(\epsilon+\epsilon_5)$.

By considering the variable transformation $$\Gamma^{*} = \Gamma\begin{pmatrix}
       S_r^T P_l & 0 \\ 
        0& I_{p-k}
    \end{pmatrix},$$
it follows that $(d\Gamma) = (d\Gamma^*)$, since $(d\Gamma)$ is an invariant measure.
Therefore, we obtain the following equality:
\begin{align*}
    \int_{E_{\eta,r}} \prod\limits_{i=k+1}^p\Big(\dfrac{c_i}{n}\Big)^{-a_i-n/2+1}(d\Gamma) &= \int_{E_{\eta,r}} \prod\limits_{i=k+1}^p \Big(\dfrac{h}{n}+\sum\limits_{j=1}^n\Gamma_{ji}^2\hat{\lambda}_j\Big)^{-a_i-n/2+1}(d\Gamma)\\
    &= \int_{D_{\eta,l}} \prod\limits_{i=k+1}^p \Big(\dfrac{h}{n}+\sum\limits_{j=1}^n(\Gamma^{*}_{ji})^2\hat{\lambda}_j\Big)^{-a_i-n/2+1}(d\Gamma^*)\\
    &= \int_{D_{\eta,l}} \prod\limits_{i=k+1}^p\Big(\dfrac{c_i}{n}\Big)^{-a_i-n/2+1}(d\Gamma).
\end{align*}
The last equality holds due to 
$$\abs{\abs{\begin{pmatrix}
       S_r & 0 \\ 
        0& Q_2
    \end{pmatrix}-\Gamma}}_F = \abs{\abs{\begin{pmatrix}
       P_l & 0 \\ 
        0& Q_2
    \end{pmatrix}-\Gamma^*}}_F,$$
and $$\Gamma_{ji}^2 = (\Gamma^{*}_{ji})^2 \quad \text{for }i>k.$$

Therefore, we obtain the following inequality:
\begin{align}\label{eq:D/E_new}
    \dfrac{\displaystyle\int_{E_{\epsilon+\epsilon_5,r}}\prod\limits_{i=1}^p c_i^{-a_i-n/2+1}(d\Gamma)}{\displaystyle\int_{D_{\epsilon+\epsilon_5,l}}\prod\limits_{i=1}^p c_i^{-a_i-n/2+1}(d\Gamma)}
    &\leq \dfrac{\sup\limits_{E_{\epsilon+\epsilon_5,r}}\prod\limits_{i=1}^k c_i^{-a_i-n/2+1} \cdot \displaystyle\int_{E_{\epsilon+\epsilon_5,r}}\prod\limits_{i=k+1}^p c_i^{-a_i-n/2+1}(d\Gamma) }{\inf\limits_{D_{\epsilon+\epsilon_5,l}}\prod\limits_{i=1}^k c_i^{-a_i-n/2+1} \cdot \displaystyle\int_{D_{\epsilon+\epsilon_5,l}}\prod\limits_{i=k+1}^p c_i^{-a_i-n/2+1}(d\Gamma)}\nonumber\\
    &= \dfrac{\sup\limits_{D_{\epsilon+\epsilon_5,l}}\prod\limits_{i=1}^k \Big(\dfrac{c_i}{n}\Big)^{a_i+n/2-1} }{\inf\limits_{E_{\epsilon+\epsilon_5,r}}\prod\limits_{i=1}^k \Big(\dfrac{c_i}{n}\Big)^{a_i+n/2-1} }\nonumber\\
    &\leq \Bigg[\dfrac{\sup\limits_{D_{\epsilon+\epsilon_5,l}}\prod\limits_{i=1}^k \Big(\dfrac{c_i}{n}\Big) }{\inf\limits_{E_{\epsilon+\epsilon_5,r}}\prod\limits_{i=1}^k \Big(\dfrac{c_i}{n}\Big) }\Bigg]^{a_1+n/2-1} \times \dfrac{\sup\limits_{D_{\epsilon+\epsilon_5,l}}\prod\limits_{i=1}^k \Big(\dfrac{c_i}{n}\Big)^{a_i-a_1} }{\inf\limits_{E_{\epsilon+\epsilon_5,r}}\prod\limits_{i=1}^k \Big(\dfrac{c_i}{n}\Big)^{a_i-a_1} }\nonumber\\
    &\leq \Bigg[\dfrac{\sup\limits_{D_{\epsilon+\epsilon_5,l}}\prod\limits_{i=1}^k \Big(\dfrac{c_i}{n}\Big) }{\inf\limits_{E_{\epsilon+\epsilon_5,r}}\prod\limits_{i=1}^k \Big(\dfrac{c_i}{n}\Big) }\Bigg]^{a_1+n/2-1} \times \Big(\dfrac{2\lambda_{0,1}}{\lambda_{0,k}}\Big)^{kq}
\end{align}

By applying Lemma \eqref{lem:upp_bound_D/E_2}, the upper bound of \eqref{eq:D/E_new} is given by
\begin{align*}
    \eqref{eq:D/E_new}&\preccurlyeq \dfrac{\exp(n(\epsilon+\epsilon_5)^2\dfrac{\hat{\lambda}_1}{\hat{\lambda}_k})}{\exp(n(\epsilon_4-(\epsilon+\epsilon_5))^2 \min\limits_{i<k}\log(\dfrac{\hat{\lambda}_i}{\hat{\lambda}_{i+1}}))}\times \Big(\dfrac{2\lambda_{0,1}}{\lambda_{0,k}}\Big)^{kq}\\
    &\preccurlyeq \dfrac{\exp(n\epsilon^2\dfrac{\lambda_{0,1}}{\lambda_{0,k}})}{\exp(n\epsilon_4^2 \min\limits_{i<k}\log(\dfrac{\lambda_{0,i}}{\lambda_{0,i+1}}))}\times \Big(\dfrac{2\lambda_{0,1}}{\lambda_{0,k}}\Big)^{kq},
\end{align*}
where $\epsilon\succ \epsilon+\epsilon_5$ and $\epsilon\asymp \epsilon_5$.

Let $\tau = \min\limits_{i<k}\log(\dfrac{\lambda_{0,i}}{\lambda_{0,i+1}})$. Then, we get the following inequality
\begin{align*}
    &\pi\bigg(\bigg\{\Gamma\in O(p):\displaystyle\inf_{Q_1\in O_{k},\;Q_2\in O_{p-k}} \abs{\abs{\begin{pmatrix}
       Q_1 & 0 \\ 
        0& Q_2
    \end{pmatrix}-\Gamma}}_F <\epsilon\bigg\}\vert \bfX_n\bigg)\\
    &\leq \pi(\bigcup_{l=1}^L D_{\epsilon_4,l}\vert \bfX_n) + m \pi(E_{\epsilon+\epsilon_5,1}\vert \bfX_n)\\
    &\leq \pi(\bigcup_{l=1}^L D_{\epsilon_4,l}\vert \bfX_n) + M (O(k),\abs{\abs{\cdot}}_F ,\epsilon_5) \cdot \dfrac{\exp(n\epsilon^2 \dfrac{\lambda_{0,1} }{\lambda_{0,k}})\cdot\Big(\dfrac{2\lambda_{0,1}}{\lambda_{0,k}}\Big)^{kq}}{\exp(\tau n\epsilon_4^2)} \cdot \pi(D_{\epsilon+\epsilon_5,1}\vert \bfX_n) \\
    &\leq \Bigg(1+ M (O(k),\abs{\abs{\cdot}}_F ,\epsilon_5) \cdot \dfrac{\exp(n\epsilon^2 \dfrac{\lambda_{0,1} }{\lambda_{0,k}})\cdot\Big(\dfrac{2\lambda_{0,1}}{\lambda_{0,k}}\Big)^{kq}}{\exp(\tau n\epsilon_4^2)} \Bigg)\pi(\bigcup_{l=1}^L D_{\epsilon_4,l}\vert \bfX_n).
\end{align*}

We obtain the following probability 
\begin{align*}
    \pi(\bigcup_{l=1}^L D_{\epsilon_4,l} \vert \bfX_n) &\geq \Bigg(1+ M (O(k),\abs{\abs{\cdot}}_F ,\epsilon_5) \cdot \dfrac{\exp(n\epsilon^2 \dfrac{\lambda_{0,1} }{\lambda_{0,k}})\cdot\Big(\dfrac{2\lambda_{0,1}}{\lambda_{0,k}}\Big)^{kq}}{\exp(\tau n\epsilon_4^2)} \Bigg) ^{-1} (1-\big(\dfrac{n\lambda_{0,k}+p}{p}\big)^{-\epsilon_4^2 a_n})\\
    &\geq \Bigg(1+  \exp( n\epsilon^2 \dfrac{\lambda_{0,1} }{\lambda_{0,k}} -\tau n\epsilon_4^2  )\cdot\Big(\dfrac{2\lambda_{0,1}}{\lambda_{0,k}}\Big)^{kq}\Bigg) ^{-1} (1-\big(\dfrac{n\lambda_{0,k}+p}{p}\big)^{-\epsilon_4^2 a_n})\\
    &\geq \Bigg(1+  \exp(-\dfrac{\tau}{2}n\epsilon_4^2)\cdot\Big(\dfrac{2\lambda_{0,1}}{\lambda_{0,k}}\Big)^{kq}\Bigg) ^{-1} (1-\big(\dfrac{n\lambda_{0,k}+p}{p}\big)^{-\epsilon_4^2 a_n}),
\end{align*}
where $ n\epsilon^2 \dfrac{\lambda_{0,1} }{\lambda_{0,k}}\leq \dfrac{\tau}{2}n\epsilon_4^2 $.

Therefore, we obtain:
\begin{align*}
    \pi((\bigcup_{l=1}^L D_{\epsilon_4,l})^c \vert \bfX_n) &\leq 1 - \Bigg(1+  \exp(-\dfrac{\tau}{2}n\epsilon_4^2)\cdot\Big(\dfrac{2\lambda_{0,1}}{\lambda_{0,k}}\Big)^{kq}\Bigg) ^{-1} (1-\big(\dfrac{n\lambda_{0,k}+p}{p}\big)^{-\epsilon_4^2 a_n})\\
    & = \dfrac{\exp(-\dfrac{\tau}{2}n\epsilon_4^2)\cdot\Big(\dfrac{2\lambda_{0,1}}{\lambda_{0,k}}\Big)^{kq}+\big(\dfrac{n\lambda_{0,k}+p}{p}\big)^{-\epsilon_4^2 a_n} }{1+  \exp(-\dfrac{\tau}{2}n\epsilon_4^2)}\\
    &= O\Big(\exp(-\dfrac{\tau}{2}n\epsilon_4^2)\cdot\Big(\dfrac{2\lambda_{0,1}}{\lambda_{0,k}}\Big)^{kq} \Big)  + O\Big(\big(\dfrac{n\lambda_{0,k}+p}{p}\big)^{-\epsilon_4^2 a_n}\Big),
\end{align*}
where $\epsilon_4\succ (n\tau)^{-1/2}$.

In summary, if there exists $\epsilon>0$ satisfying
$$\dfrac{np}{a_n}\prec \epsilon^2 \leq \dfrac{\lambda_{0,k}}{n\lambda_{0,1}}\dfrac{\tau}{2}n\epsilon_4^2,$$ 
and 
$$\epsilon_4\succ \epsilon\vee (n\tau)^{-1/2},$$
then the following holds
$$\pi((\bigcup_{l=1}^L D_{\epsilon_4,l})^c\vert \bfX_n) = O\Big(\exp(-\dfrac{\tau}{2}n\epsilon_4^2) \cdot\Big(\dfrac{2\lambda_{0,1}}{\lambda_{0,k}}\Big)^{kq}\Big)  + O\Big(\big(\dfrac{n\lambda_{0,k}+p}{p}\big)^{-\epsilon_4^2 a_n}\Big).$$

In particular, the same posterior bound holds under the simplified condition 
$\dfrac{np}{a_n}\prec \dfrac{\lambda_{0,k}}{n\lambda_{0,1}}\dfrac{\tau}{2}n\epsilon_4^2$ 
and $\epsilon_4 \succ \sqrt{\dfrac{np}{a_n}} \vee (n\tau)^{-1/2}$.

\end{proof}


\subsubsection*{Proof of Theorem \ref{thm:post_exp_equal_a}}
\begin{proof}

The proof follows the same procedure as the proof of Theorem \ref{thm:post_exp_diff_a}. By using the reparemetrization as we mentioned at \eqref{main-eq:spiked_repara}, we obtain the following inequality
    \begin{align*}
        \bbE\Big[\dfrac{\tilde{\lambda}_i-\lambda_{0,i}}{{\lambda_{0,i}}}\big\vert \bfX_n \Big]  &=\dfrac{\displaystyle\int\dfrac{\tilde{\lambda_i}}{{\lambda_{0,i}}} \prod {\tilde{\lambda}_i}^{-a_i-\frac{n}{2}}\exp(-\dfrac{\tilde{c}_i}{2\tilde{\lambda}_i})(d\tilde{\Lambda })(d\tilde{\Gamma})}{\displaystyle\int\prod \tilde{\lambda}_i^{-a_i-\frac{n}{2}}\exp(-\dfrac{\tilde{c}_i}{2\tilde{\lambda}_i})(d\tilde{\Lambda })(d\tilde{\Gamma})}-1\\
        &=\dfrac{\displaystyle\int f(\tilde{c})\Big(\prod \tilde{c}_{i}\Big)^{-a_i-\frac{n}{2}+1} (d\tilde{\Gamma}) }{\displaystyle\int \Big(\prod \tilde{c}_{i}\Big)^{-a_i-\frac{n}{2}+1}(d\tilde{\Gamma})}-1\\
        &\leq\dfrac{\displaystyle\int_{D_\epsilon} f(\tilde{c})\Big(\prod \tilde{c}_{i}\Big)^{-a_i-\frac{n}{2}+1} (d\tilde{\Gamma}) }{\displaystyle\int_{D_\epsilon} \Big(\prod \tilde{c}_{i}\Big)^{-a_i-\frac{n}{2}+1}(d\tilde{\Gamma})} + \dfrac{\displaystyle\int_{D_\epsilon^c} f(\tilde{c})\Big(\prod \tilde{c}_{i}\Big)^{-a_i-\frac{n}{2}+1} (d\tilde{\Gamma}) }{\displaystyle\int \Big(\prod \tilde{c}_{i}\Big)^{-a_i-\frac{n}{2}+1}(d\tilde{\Gamma})}-1\\
        &\leq \sup_{D_\epsilon}f(\tilde{c}) -1 + \sup_{D_\epsilon^c}f(\tilde{c}) \dfrac{\displaystyle\int_{D_\epsilon^c} \Big(\prod \tilde{c}_{i}\Big)^{-a_i-\frac{n}{2}+1} (d\tilde{\Gamma}) }{\displaystyle\int \Big(\prod \tilde{c}_{i}\Big)^{-a_i-\frac{n}{2}+1}(d\tilde{\Gamma})}\\
        &= \sup_{D_\epsilon}f(\tilde{c})-1 + \sup_{D_\epsilon^c}f(\tilde{c}) \dfrac{\displaystyle\int_{(\bigcup_{l=1}^L D_{\epsilon,l})^c} \Big(\prod c_{i}\Big)^{-a_i-\frac{n}{2}+1} (d{\Gamma}) }{\displaystyle\int \Big(\prod c_{i}\Big)^{-a_i-\frac{n}{2}+1}(d{\Gamma})},
    \end{align*}
    where 
    $$f(\tilde{c})=\bbE_{\lambda_i}[\dfrac{\lambda_i}{\lambda_{0,i}}]= \dfrac{\Gamma(a_i+n/2-2)}{2\Gamma(a_i+n/2-1)} \dfrac{\tilde{c}_{i}}{\lambda_{0,i}},$$
    with $\lambda_i\overset{iid}{\sim}InvGam(a_i+\dfrac{n}{2}-1,\dfrac{\tilde{c}_{i}}{2})$. 
    
    Since $\sup_{D_\epsilon^c} f(\tilde{c}) \preccurlyeq \dfrac{\lambda_{0,1}}{\lambda_{0,i}}$, applying Lemma \ref{main-lem:shrink_post_equal_a}, we obtain:
    \begin{align*}
        \bbE\Big[\dfrac{\tilde{\lambda}_i-\lambda_{0,i}}{{\lambda_{0,i}}}\big\vert \bfX_n \Big] 
        &\leq \dfrac{n}{n+2a_i-4}\dfrac{1}{\lambda_{0,i}}\cdot\sup_{D_\epsilon}\Big( \dfrac{\tilde{c}_i}{n}\Big) - 1 \\
        &+O\Big( \dfrac{\lambda_{0,1}}{\lambda_{0,i}}\exp(-\dfrac{\tau}{2}n\epsilon^2)\cdot\Big(\dfrac{2\lambda_{0,1}}{\lambda_{0,k}}\Big)^{kq} \Big)  + O\Big( \dfrac{\lambda_{0,1}}{\lambda_{0,i}}\big(\dfrac{n\lambda_{0,k}+p}{p}\big)^{-\epsilon^2 a_n}\Big)\\
        &= \dfrac{n}{n+2a_i-4}\dfrac{1}{\lambda_{0,i}}\cdot\sup_{D_\epsilon}\Big( \dfrac{c_i}{n}\Big) - 1 \\
        &+O\Big( \dfrac{\lambda_{0,1}}{\lambda_{0,i}}\exp(-\dfrac{\tau}{2}n\epsilon^2)\cdot\Big(\dfrac{2\lambda_{0,1}}{\lambda_{0,k}}\Big)^{kq} \Big)  + O\Big( \dfrac{\lambda_{0,1}}{\lambda_{0,i}}\big(\dfrac{n\lambda_{0,k}+p}{p}\big)^{-\epsilon^2 a_n}\Big).
    \end{align*}
    The last equality follows from $c_i = \tilde{c}_i$ on $D_\epsilon$.

    Since the $c_i$ is given by 
    $$\dfrac{c_i}{n} 
    = \sum\limits_{j=1}^k\Gamma_{ji}^2\Big\{(1+\beta_j)\lambda_{0,j}+(\bar{d}\dfrac{p}{n}+\alpha_j\sqrt{\dfrac{p}{n}})+\dfrac{h}{n}\Big\}+ \sum\limits_{j=k+1}^n\Gamma_{ji}^2\Big(\bar{d}\dfrac{p}{n}+\alpha_j\sqrt{\dfrac{p}{n}}+\dfrac{h}{n} \Big) + \sum\limits_{j=n+1}^p\Gamma_{ji}^2\dfrac{h}{n},$$
    we obtain the upper bound of $\sup_{D_\epsilon}\Big( \dfrac{c_i}{n}\Big)$ as follows:
    \begin{align*}
        \sup_{D_\epsilon}\Big( \dfrac{c_i}{n}\Big) 
        &\leq (1-\epsilon^2)\Big[(1+\beta_i)\lambda_{0,i}+(\bar{d}\dfrac{p}{n}+\alpha_i\sqrt{\dfrac{p}{n}})+\dfrac{h}{n}\Big]  + \epsilon^2\Big[(1+\beta_1)\lambda_{0,1}+(\bar{d}\dfrac{p}{n}+\alpha_1\sqrt{\dfrac{p}{n}})+\dfrac{h}{n}\Big] \\
        &\leq \Big[(1+\beta_i)\lambda_{0,i}+(\bar{d}\dfrac{p}{n}+\alpha_i\sqrt{\dfrac{p}{n}})+\dfrac{h}{n}\Big] + 2\lambda_{0,1}\epsilon^2.
    \end{align*}

    Therefore, the posterior expectation given by,
    \begin{align*}
        \bbE\Big[\dfrac{\tilde{\lambda}_i-\lambda_{0,i}}{{\lambda_{0,i}}}\big\vert \bfX_n \Big] 
        &\leq \dfrac{n}{n+2a_i-4}\dfrac{1}{\lambda_{0,i}}\Big[(1+\beta_i)\lambda_{0,i}+\bar{d}\dfrac{p}{n}+\alpha_i\sqrt{\dfrac{p}{n}}+\dfrac{h}{n}\Big]-1 \\
        & + O\Big( \epsilon^2\dfrac{\lambda_{0,1}}{\lambda_{0,i}} + \dfrac{\lambda_{0,1}}{\lambda_{0,i}}\exp(-\dfrac{\tau}{2}n\epsilon^2)\cdot\Big(\dfrac{2\lambda_{0,1}}{\lambda_{0,k}}\Big)^{kq}+\dfrac{\lambda_{0,1}}{\lambda_{0,i}}\big(\dfrac{n\lambda_{0,k}+p}{p}\big)^{-\epsilon^2 a_n} \Big).
    \end{align*}
    
    Next, we derive the lower bound of the posterior expectation. Applying Lemma \ref{main-lem:shrink_post_equal_a}, we obtain
    \begin{align*}
        \bbE\Big[\dfrac{\tilde{\lambda}_i}{{\lambda_{0,i}}}\big\vert \bfX_n \Big]  &=\dfrac{\displaystyle\int f(\tilde{c})\Big(\prod \tilde{c}_{i}\Big)^{-a_i-\frac{n}{2}+1} (d\tilde{\Gamma}) }{\displaystyle\int \Big(\prod \tilde{c}_{i}\Big)^{-a_i-\frac{n}{2}+1}(d\tilde{\Gamma})}\\
        &= \dfrac{\displaystyle\int_{D_\epsilon} f(\tilde{c})\Big(\prod \tilde{c}_{i}\Big)^{-a_i-\frac{n}{2}+1} (d\tilde{\Gamma}) }{\displaystyle\int_{D_\epsilon} \Big(\prod \tilde{c}_{i}\Big)^{-a_i-\frac{n}{2}+1}(d\tilde{\Gamma})} \cdot \dfrac{\displaystyle\int_{D_\epsilon} \Big(\prod \tilde{c}_{i}\Big)^{-a_i-\frac{n}{2}+1} (d\tilde{\Gamma}) }{\displaystyle\int \Big(\prod \tilde{c}_{i}\Big)^{-a_i-\frac{n}{2}+1}(d\tilde{\Gamma})} + \dfrac{\displaystyle\int_{D_\epsilon^c} f(\tilde{c})\Big(\prod \tilde{c}_{i}\Big)^{-a_i-\frac{n}{2}+1} (d\tilde{\Gamma}) }{\displaystyle\int \Big(\prod \tilde{c}_{i}\Big)^{-a_i-\frac{n}{2}+1}(d\tilde{\Gamma})}\\
        &\geq \inf_{D_\epsilon} f(\tilde{c}) \cdot \dfrac{\displaystyle\int_{D_\epsilon} \Big(\prod \tilde{c}_{i}\Big)^{-a_i-\frac{n}{2}+1} (d\tilde{\Gamma}) }{\displaystyle\int \Big(\prod \tilde{c}_{i}\Big)^{-a_i-\frac{n}{2}+1}(d\tilde{\Gamma})} + \inf_{O(p)} f(\tilde{c})  \dfrac{\displaystyle\int_{D_\epsilon^c} \Big(\prod \tilde{c}_{i}\Big)^{-a_i-\frac{n}{2}+1} (d\tilde{\Gamma}) }{\displaystyle\int \Big(\prod \tilde{c}_{i}\Big)^{-a_i-\frac{n}{2}+1}(d\tilde{\Gamma})}.
    \end{align*}

    By using the reparemetrization as we mentioned at \eqref{main-eq:spiked_repara}, we obtain
    \begin{align*}
    &\inf_{D_\epsilon} f(\tilde{c}) \cdot \dfrac{\displaystyle\int_{D_\epsilon} \Big(\prod \tilde{c}_{i}\Big)^{-a_i-\frac{n}{2}+1} (d\tilde{\Gamma}) }{\displaystyle\int \Big(\prod \tilde{c}_{i}\Big)^{-a_i-\frac{n}{2}+1}(d\tilde{\Gamma})} + \inf_{O(p)} f(\tilde{c})  \dfrac{\displaystyle\int_{D_\epsilon^c} \Big(\prod \tilde{c}_{i}\Big)^{-a_i-\frac{n}{2}+1} (d\tilde{\Gamma}) }{\displaystyle\int \Big(\prod \tilde{c}_{i}\Big)^{-a_i-\frac{n}{2}+1}(d\tilde{\Gamma})}\\
        &= \inf_{D_\epsilon} f(\tilde{c}) \cdot \dfrac{\displaystyle\int_{(\bigcup_{l=1}^L D_{\epsilon,l})} \Big(\prod c_{i}\Big)^{-a_i-\frac{n}{2}+1} (d{\Gamma}) }{\displaystyle\int \Big(\prod c_{i}\Big)^{-a_i-\frac{n}{2}+1}(d{\Gamma})} + \inf_{O(p)} f(\tilde{c}) \cdot \dfrac{\displaystyle\int_{(\bigcup_{l=1}^L D_{\epsilon,l})^c} \Big(\prod c_{i}\Big)^{-a_i-\frac{n}{2}+1} (d{\Gamma}) }{\displaystyle\int \Big(\prod c_{i}\Big)^{-a_i-\frac{n}{2}+1}(d{\Gamma})}\\
        &\geq \inf_{D_\epsilon} f(\tilde{c}) \times \Bigg(1-\dfrac{\displaystyle\int_{(\bigcup_{l=1}^L D_{\epsilon,l})^c} \Big(\prod c_{i}\Big)^{-a_i-\frac{n}{2}+1} (d{\Gamma}) }{\displaystyle\int \Big(\prod c_{i}\Big)^{-a_i-\frac{n}{2}+1}(d{\Gamma})}\Bigg)\\
        &= \inf_{D_\epsilon} f(c) -\inf_{D_\epsilon} f(c)\times \dfrac{\displaystyle\int_{(\bigcup_{l=1}^L D_{\epsilon,l})^c} \Big(\prod c_{i}\Big)^{-a_i-\frac{n}{2}+1} (d{\Gamma}) }{\displaystyle\int \Big(\prod c_{i}\Big)^{-a_i-\frac{n}{2}+1}(d{\Gamma})}\\
        &\geq \inf_{D_\epsilon} f(c) -O\Big(\dfrac{\lambda_{0,1}}{\lambda_{0,i}}\Big)\times O\Big( \exp(-\dfrac{\tau}{2}n\epsilon^2)\cdot\Big(\dfrac{2\lambda_{0,1}}{\lambda_{0,k}}\Big)^{kq} + \big(\dfrac{n\lambda_{0,k}+p}{p}\big)^{-\epsilon^2 a_n}\Big)
    \end{align*}
    The first inequality holds by $ f(\tilde{c})\geq 0$, and the last inequality holds by Lemma \ref{main-lem:shrink_post_equal_a} and $f(c) = O(\dfrac{\lambda_{0,1}}{\lambda_{0,i}})$.
    
    We obtain the lower bound of $\inf_{D_\epsilon}\Big( \dfrac{c_i}{n}\Big)$ as follows:
        \begin{align*}
            \inf_{D_\epsilon}\Big( \dfrac{c_i}{n}\Big) &\geq (1-\epsilon^2)\Big[(1+\beta_i)\lambda_{0,i}+(\bar{d}\dfrac{p}{n}+\alpha_i\sqrt{\dfrac{p}{n}})+\dfrac{h}{n}\Big]  + \epsilon^2\Big[\dfrac{h}{n}\Big] \\
            &\geq \Big[(1+\beta_i)\lambda_{0,i}+(\bar{d}\dfrac{p}{n}+\alpha_i\sqrt{\dfrac{p}{n}})+\dfrac{h}{n}\Big]  - 2\epsilon^2\lambda_{0,1}.
        \end{align*}
    
    Therefore, the lower bound of posterior expectation is given by 
    \begin{align*}
        \bbE\Big[\dfrac{\tilde{\lambda}_i-\lambda_{0,i}}{{\lambda_{0,i}}}\big\vert \bfX_n \Big] &\geq 
      \Bigg[\dfrac{n}{n+2a_i-4}\dfrac{1}{\lambda_{0,i}} \bigg(\Big[(1+\beta_i)\lambda_{0,i}+(\bar{d}\dfrac{p}{n}+\alpha_i\sqrt{\dfrac{p}{n}})+\dfrac{h}{n}\Big]   - 2\epsilon^2\lambda_{0,1}\bigg)  \Bigg]-1 \\
      &-O\Big(\dfrac{\lambda_{0,1}}{\lambda_{0,i}}\Big)\times O\Big( \exp(-\dfrac{\tau}{2}n\epsilon^2)\cdot\Big(\dfrac{2\lambda_{0,1}}{\lambda_{0,k}}\Big)^{kq} + \big(\dfrac{n\lambda_{0,k}+p}{p}\big)^{-\epsilon^2 a_n}\Big)\\
      &= \dfrac{n}{n+2a_i-4}\dfrac{1}{\lambda_{0,i}}\Big[(1+\beta_i)\lambda_{0,i}+\bar{d}\dfrac{p}{n}+\alpha_i\sqrt{\dfrac{p}{n}}+\dfrac{h}{n}\Big]-1  \\
      &+ O\Big(\epsilon^2\dfrac{\lambda_{0,1}}{\lambda_{0,i}} +\dfrac{\lambda_{0,1}}{\lambda_{0,i}}\exp(-\dfrac{\tau}{2}n\epsilon^2)\cdot\Big(\dfrac{2\lambda_{0,1}}{\lambda_{0,k}}\Big)^{kq}  +\dfrac{\lambda_{0,1}}{\lambda_{0,i}}\big(\dfrac{n\lambda_{0,k}+p}{p}\big)^{-\epsilon^2 a_n}\Big).
    \end{align*}
    
    Thus, we conclude that the posterior expectation satisfies:
    \begin{align*}
        \bbE\Big[\dfrac{\tilde{\lambda}_i-\lambda_{0,i}}{{\lambda_{0,i}}}\big\vert \bfX_n \Big]  &=\dfrac{n}{n+2a_i-4}\dfrac{1}{\lambda_{0,i}}\Big[(1+\beta_i)\lambda_{0,i}+\bar{d}\dfrac{p}{n}+\alpha_i\sqrt{\dfrac{p}{n}}+\dfrac{h}{n}\Big]-1\\
        &+O\Big( \epsilon^2\dfrac{\lambda_{0,1}}{\lambda_{0,i}} + \dfrac{\lambda_{0,1}}{\lambda_{0,i}}\exp(-\dfrac{\tau}{2}n\epsilon^2)\cdot\Big(\dfrac{2\lambda_{0,1}}{\lambda_{0,k}}\Big)^{kq}+\dfrac{\lambda_{0,1}}{\lambda_{0,i}}\big(\dfrac{n\lambda_{0,k}+p}{p}\big)^{-\epsilon^2 a_n} \Big).
    \end{align*}

\end{proof}    

\subsubsection*{Proof of Theorem \ref{main-cor:conv_rate_eigenvalue_equal_a}}
\begin{proof}

By Theorem \ref{thm:post_exp_equal_a}, we obtain the following posterior expectation,
\begin{align*}
    \bbE\Big[\dfrac{\lambda_i-\lambda_{0,i}}{{\lambda_{0,i}}}\big\vert \bfX_n \Big]  &=\dfrac{n}{n+2a_i-4}\dfrac{1}{\lambda_{0,i}}\Big[(1+\beta_i)\lambda_{0,i}+\bar{d}\dfrac{p}{n}+\alpha_i\sqrt{\dfrac{p}{n}}+\dfrac{h}{n}\Big]-1\\
    &+O\Big( \epsilon^2\dfrac{\lambda_{0,1}}{\lambda_{0,i}} + \dfrac{\lambda_{0,1}}{\lambda_{0,i}}\exp(-\dfrac{\tau}{2}n\epsilon^2)\cdot\Big(\dfrac{2\lambda_{0,1}}{\lambda_{0,k}}\Big)^{kq}+\dfrac{\lambda_{0,1}}{\lambda_{0,i}}\big(\dfrac{n\lambda_{0,k}+p}{p}\big)^{-\epsilon^2 a_n} \Big)\\
    &= \Big(\dfrac{(1+\beta_i)n}{n+2a_i-4}-1\Big) + \dfrac{n}{n+2a_i-4}\dfrac{1}{\lambda_{0,i}}\Big[\bar{d}\dfrac{p}{n}+\alpha_i\sqrt{\dfrac{p}{n}}+\dfrac{h}{n}\Big]\\
    &+O\Big( \epsilon^2\dfrac{\lambda_{0,1}}{\lambda_{0,i}} + \dfrac{\lambda_{0,1}}{\lambda_{0,i}}\exp(-\dfrac{\tau}{2}n\epsilon^2)\cdot\Big(\dfrac{2\lambda_{0,1}}{\lambda_{0,k}}\Big)^{kq}+\dfrac{\lambda_{0,1}}{\lambda_{0,i}}\big(\dfrac{n\lambda_{0,k}+p}{p}\big)^{-\epsilon^2 a_n} \Big)\\
    &= \dfrac{\beta_i n -2a_i+4}{n+2a_i-4} + O\Big(\dfrac{p}{n\lambda_{0,i}}\Big)\\
    &+O\Big( \epsilon^2\dfrac{\lambda_{0,1}}{\lambda_{0,i}} + \dfrac{\lambda_{0,1}}{\lambda_{0,i}}\exp(-\dfrac{\tau}{2}n\epsilon^2)\cdot\Big(\dfrac{2\lambda_{0,1}}{\lambda_{0,k}}\Big)^{kq}+\dfrac{\lambda_{0,1}}{\lambda_{0,i}}\big(\dfrac{n\lambda_{0,k}+p}{p}\big)^{-\epsilon^2 a_n} \Big).
\end{align*}

If $a_i\preccurlyeq n^{1/2}$, then the following inequality holds,
$$\dfrac{\beta_i n -2a_i+4}{n+2a_i-4}\preccurlyeq \beta_i = O(n^{-1/2+\delta}).$$

Therefore, the posterior expectation is given by
\begin{align*}
    \bbE\Big[\dfrac{\lambda_i-\lambda_{0,i}}{{\lambda_{0,i}}}\big\vert \bfX_n \Big]&= O(n^{-1/2+\delta}) + O\Big(\dfrac{p}{n\lambda_{0,i}}\Big)\\
    &+O\Big( \epsilon^2\dfrac{\lambda_{0,1}}{\lambda_{0,i}} + \dfrac{\lambda_{0,1}}{\lambda_{0,i}}\exp(-\dfrac{\tau}{2}n\epsilon^2)\cdot\Big(\dfrac{2\lambda_{0,1}}{\lambda_{0,k}}\Big)^{kq}+\dfrac{\lambda_{0,1}}{\lambda_{0,i}}\big(\dfrac{n\lambda_{0,k}+p}{p}\big)^{-\epsilon^2 a_n} \Big).
\end{align*}

By the assumptions of the theorem, in particular $\lambda_{0,1}/\lambda_{0,k} = O(1)$ and $\epsilon = n^{-1/4}$, we apply Lemma \ref{lem:post_exp_equal_ordered} to obtain
\begin{align*}
    \bbE\Big[\dfrac{\lambda_{(i)}-\lambda_{0,i}}{{\lambda_{0,i}}}\big\vert \bfX_n \Big] &=  \bbE\Big[\dfrac{\lambda_{i}-\lambda_{0,i}}{{\lambda_{0,i}}}\big\vert \bfX_n \Big] + \bbE\Big[\dfrac{\lambda_{(i)}-\lambda_{i}}{{\lambda_{0,i}}}\big\vert \bfX_n \Big] \\
    &= O(n^{-1/2+\delta}) + O\Big(\dfrac{p}{n\lambda_{0,i}}\Big)\\
    &+O\Big( \epsilon^2\dfrac{\lambda_{0,1}}{\lambda_{0,i}} + \dfrac{\lambda_{0,1}}{\lambda_{0,i}}\exp(-\dfrac{\tau}{2}n\epsilon^2)\cdot\Big(\dfrac{2\lambda_{0,1}}{\lambda_{0,k}}\Big)^{kq}+\dfrac{\lambda_{0,1}}{\lambda_{0,i}}\big(\dfrac{n\lambda_{0,k}+p}{p}\big)^{-\epsilon^2 a_n} \Big)\\
    &+\dfrac{\lambda_{0,1}}{\lambda_{0,i}}\cdot\bigg(O\Big( \dfrac{1}{n} + \exp(-\dfrac{\tau}{2}n\epsilon^2)\cdot\Big(\dfrac{2\lambda_{0,1}}{\lambda_{0,k}}\Big)^{kq}+\big(\dfrac{n\lambda_{0,k}+p}{p}\big)^{-\epsilon^2 a_n} \Big)\bigg)\\
    &=  O(n^{-1/2+\delta}) + O\Big(\dfrac{p}{n\lambda_{0,i}}\Big) \\
    &+ O\Big( (\epsilon^2+\dfrac{1}{n})\dfrac{\lambda_{0,1}}{\lambda_{0,i}} + \dfrac{\lambda_{0,1}}{\lambda_{0,i}}\exp(-\dfrac{\tau}{2}n\epsilon^2)\cdot\Big(\dfrac{2\lambda_{0,1}}{\lambda_{0,k}}\Big)^{kq}+\dfrac{\lambda_{0,1}}{\lambda_{0,i}}\big(\dfrac{n\lambda_{0,k}+p}{p}\big)^{-\epsilon^2 a_n} \Big),
\end{align*}
where $\lambda_{(i)}$ is the $i$-th largest eigenvalue of $\Sigma$.

Since $\lambda_{0,1}/\lambda_{0,k}$ is bounded by a positive constant, $a_1,\ldots,a_k \prec n^{1/2}$, and $\epsilon \asymp n^{-1/4}$, it follows that
\begin{equation*}
    \epsilon^2 + \exp(-\dfrac{\tau}{2}n\epsilon^2)\cdot\Big(\dfrac{2\lambda_{0,1}}{\lambda_{0,k}}\Big)^{kq}+\big(\dfrac{n\lambda_{0,k}+p}{p}\big)^{-\epsilon^2 a_n} \preccurlyeq n^{-1/2+\delta}+\dfrac{p}{n\lambda_{0,i}}.
\end{equation*}

Therefore, we conclude
\begin{align*}
    \bbE\Big[\dfrac{\lambda_{(i)}-\lambda_{0,i}}{{\lambda_{0,i}}}\big\vert \bfX_n \Big]=   O\Big(\dfrac{p}{n\lambda_{0,i}}\Big)+ O(n^{-1/2+\delta}).
\end{align*}

\end{proof}


\subsubsection*{Proof of Theorem \ref{thm:post_exp_diff_a}}
\begin{proof}
    Using the Lemma \ref{main-lem:shrink_post_diff_a}, we obtain the following inequality:
    \begin{align*}
        \bbE\Big[\dfrac{\lambda_i-\lambda_{0,i}}{{\lambda_{0,i}}}\big\vert \bfX_n \Big]  &=\dfrac{\displaystyle\int\dfrac{\lambda_i}{{\lambda_{0,i}}} \prod \lambda_i^{-a_i-\frac{n}{2}}\exp(-\dfrac{c_i}{2\lambda_i})(d\Lambda) (d\Gamma)}{\displaystyle\int\prod \lambda_i^{-a_i-\frac{n}{2}}\exp(-\dfrac{c_i}{2\lambda_i})(d\Lambda) (d\Gamma)}-1\nonumber \\
        &=\dfrac{\displaystyle\int f(c)\Big(\prod c_{i}\Big)^{-a_i-\frac{n}{2}+1} (d\Gamma) }{\displaystyle\int \Big(\prod c_{i}\Big)^{-a_i-\frac{n}{2}+1}(d\Gamma)}-1\\
        &\leq \dfrac{\displaystyle\int_{D_\epsilon} f(c)\Big(\prod c_{i}\Big)^{-a_i-\frac{n}{2}+1} (d\Gamma)}{\displaystyle\int_{D_\epsilon} \Big(\prod c_{i}\Big)^{-a_i-\frac{n}{2}+1}(d\Gamma)} + \dfrac{\displaystyle\int_{D_\epsilon^c} f(c)\Big(\prod c_{i}\Big)^{-a_i-\frac{n}{2}+1} (d\Gamma)}{\displaystyle\int \Big(\prod c_{i}\Big)^{-a_i-\frac{n}{2}+1}(d\Gamma)}-1\\
        &\leq \sup_{D_\epsilon} f(c)-1 + \sup_{D_\epsilon^c} f(c) \cdot \dfrac{\displaystyle\int_{D_\epsilon^c} f(c)\Big(\prod c_{i}\Big)^{-a_i-\frac{n}{2}+1} (d\Gamma) }{\displaystyle\int \Big(\prod c_{i}\Big)^{-a_i-\frac{n}{2}+1}(d\Gamma)}
    \end{align*}
    where 
    $$f(c)=\bbE_{\lambda_i}[\dfrac{\lambda_i}{\lambda_{0,i}}]= \dfrac{\Gamma(a_i+n/2-2)}{2\Gamma(a_i+n/2-1)} \dfrac{c_{i}}{\lambda_{0,i}} ,$$
    with $\lambda_i\overset{iid}{\sim}InvGam(a_i+\dfrac{n}{2}-1,\dfrac{c_{i}}{2})$. 

      By the same procedure of proof of Theorem \ref{thm:post_exp_equal_a}, we conclude that the posterior expectation satisfies
    \begin{align*}
        \bbE\Big[\dfrac{{\lambda}_i-\lambda_{0,i}}{{\lambda_{0,i}}}\big\vert \bfX_n \Big]  &=\dfrac{n}{n+2a_i-4}\dfrac{1}{\lambda_{0,i}}\Big[(1+\beta_i)\lambda_{0,i}+\bar{d}\dfrac{p}{n}+\alpha_i\sqrt{\dfrac{p}{n}}+\dfrac{h}{n}\Big]-1\\
        &+O\Big( \epsilon^2\dfrac{\lambda_{0,1}}{\lambda_{0,i}} + \dfrac{\lambda_{0,1}}{\lambda_{0,i}}\exp(-\dfrac{\kappa}{2}\epsilon)+\dfrac{\lambda_{0,1}}{\lambda_{0,i}}\big(\dfrac{n\lambda_{0,k}+p}{p}\big)^{-\epsilon^2 a_n} \Big).
    \end{align*}

\end{proof}


\subsubsection*{Proof of Theorem \ref{main-cor:conv_rate_eigenvalue_diff_a}}
\begin{proof}
    Let $2a_i-4 = \dfrac{nt}{\hat{\lambda}_i-t}$. Then, by Theorem \ref{thm:post_exp_diff_a}, the posterior expectation can be written as:
\begin{align*}
    \bbE\Big[\dfrac{\lambda_i-\lambda_{0,i}}{{\lambda_{0,i}}}\big\vert \bfX_n \Big] 
    &=\dfrac{1}{\lambda_{0,i}}\Big[\dfrac{n\hat{\lambda_{i}}}{n+2a_i-4}-\lambda_{0,i}\Big]\\
    &+ O\Big( \epsilon^2\dfrac{\lambda_{0,1}}{\lambda_{0,i}} + \dfrac{\lambda_{0,1}}{\lambda_{0,i}}\exp(-\dfrac{\kappa}{2}\epsilon)+\dfrac{\lambda_{0,1}}{\lambda_{0,i}}\big(\dfrac{n\lambda_{0,k}+p}{p}\big)^{-\epsilon^2 a_n} \Big)\\
    &=\dfrac{1}{\lambda_{0,i}}\Big[\hat{\lambda_{i}}-t-\lambda_{0,i}\Big]\\
    &+ O\Big( \epsilon^2\dfrac{\lambda_{0,1}}{\lambda_{0,i}} + \dfrac{\lambda_{0,1}}{\lambda_{0,i}}\exp(-\dfrac{\kappa}{2}\epsilon)+\dfrac{\lambda_{0,1}}{\lambda_{0,i}}\big(\dfrac{n\lambda_{0,k}+p}{p}\big)^{-\epsilon^2 a_n} \Big)\\
    &=\dfrac{1}{\lambda_{0,i}}\Big[\dfrac{h}{n}+\beta_i\lambda_{0,i}+\bar{d}\dfrac{p}{n}+\alpha_i\sqrt{\dfrac{p}{n}}-t\Big]\\
    &+ O\Big( \epsilon^2\dfrac{\lambda_{0,1}}{\lambda_{0,i}} + \dfrac{\lambda_{0,1}}{\lambda_{0,i}}\exp(-\dfrac{\kappa}{2}\epsilon)+\dfrac{\lambda_{0,1}}{\lambda_{0,i}}\big(\dfrac{n\lambda_{0,k}+p}{p}\big)^{-\epsilon^2 a_n} \Big).
\end{align*}

If $t\in [\hat{\lambda}_{k+1},\hat{\lambda}_{n}]$, then for some constant $\alpha_0\in[-C,C]$, we have
$$t= \bar{d}\dfrac{p}{n}+\alpha_0\sqrt{\dfrac{p}{n}}.$$

Therefore, the posterior expectation simplifies to:
\begin{align*}
    \bbE\Big[\dfrac{\lambda_i-\lambda_{0,i}}{{\lambda_{0,i}}}\big\vert \bfX_n \Big]  &= \dfrac{1}{\lambda_{0,i}}(\alpha_i-\alpha_0)\sqrt{\dfrac{p}{n}}+\beta_i + \dfrac{1}{\lambda_{0,i}}\dfrac{h}{n}\\
    &+ O\Big( \epsilon^2\dfrac{\lambda_{0,1}}{\lambda_{0,i}} + \dfrac{\lambda_{0,1}}{\lambda_{0,i}}\exp(-\dfrac{\kappa}{2}\epsilon)+\dfrac{\lambda_{0,1}}{\lambda_{0,i}}\big(\dfrac{n\lambda_{0,k}+p}{p}\big)^{-\epsilon^2 a_n} \Big)\\
    &=O({\lambda_{0,i}}^{-1}\sqrt{\dfrac{p}{n}})+ O(\beta_i)\\
    &+ O\Big( \epsilon^2\dfrac{\lambda_{0,1}}{\lambda_{0,i}} + \dfrac{\lambda_{0,1}}{\lambda_{0,i}}\exp(-\dfrac{\kappa}{2}\epsilon)+\dfrac{\lambda_{0,1}}{\lambda_{0,i}}\big(\dfrac{n\lambda_{0,k}+p}{p}\big)^{-\epsilon^2 a_n} \Big).
\end{align*}

By the assumptions of the theorem, in particular $\lambda_{0,1}/\lambda_{0,k} = O(1)$ and $\epsilon = n^{-1/4}$, we apply Lemma \ref{lem:post_exp_equal_ordered} to obtain
\begin{align*}
    \bbE\Big[\dfrac{\lambda_{(i)}-\lambda_{0,i}}{{\lambda_{0,i}}}\big\vert \bfX_n \Big] &=  \bbE\Big[\dfrac{\lambda_{i}-\lambda_{0,i}}{{\lambda_{0,i}}}\big\vert \bfX_n \Big] + \bbE\Big[\dfrac{\lambda_{(i)}-\lambda_{i}}{{\lambda_{0,i}}}\big\vert \bfX_n \Big] \\
    &= O({\lambda_{0,i}}^{-1}\sqrt{\dfrac{p}{n}})+ O(\beta_i)\\
    &+ O\Big( \epsilon^2\dfrac{\lambda_{0,1}}{\lambda_{0,i}} + \dfrac{\lambda_{0,1}}{\lambda_{0,i}}\exp(-\dfrac{\kappa}{2}\epsilon)+\dfrac{\lambda_{0,1}}{\lambda_{0,i}}\big(\dfrac{n\lambda_{0,k}+p}{p}\big)^{-\epsilon^2 a_n} \Big)\\
    &+\dfrac{\lambda_{0,1}}{\lambda_{0,i}}\cdot\bigg(O\Big( \dfrac{1}{n} + \exp(-\dfrac{\kappa}{2}\epsilon)+\big(\dfrac{n\lambda_{0,k}+p}{p}\big)^{-\epsilon^2 a_n} \Big)\bigg)\\
    &=  O(\dfrac{1}{\lambda_{0,i}}\sqrt{\dfrac{p}{n}})+ O(n^{-\frac{1}{2}+\delta}) \\
    &+O\Big( \dfrac{1}{n} + \exp(-\dfrac{\kappa}{2}\epsilon)+\big(\dfrac{n\lambda_{0,k}+p}{p}\big)^{-\epsilon^2 a_n} \Big),
\end{align*}
where $\lambda_{(i)}$ is the $i$-th largest eigenvalue of $\Sigma$.

Since $2a_i - 4 = \frac{nt}{\hat{\lambda}_i - t}, \quad \text{for } i = 1, \ldots, k$, $\epsilon\asymp n^{-1/4}$, and $\lambda_{0,k}\preccurlyeq p/n^{1/4}$, it follows that
\begin{equation*}
    \epsilon^2 + \exp(-\dfrac{\kappa}{2}\epsilon)+\big(\dfrac{n\lambda_{0,k}+p}{p}\big)^{-\epsilon^2 a_n} \preccurlyeq n^{-1/2+\delta}+\dfrac{1}{\lambda_{0,i}}\sqrt{\dfrac{p}{n}}.
\end{equation*}

Therefore, we conclude
\begin{align*}
    \bbE\Big[\dfrac{\lambda_{(i)}-\lambda_{0,i}}{{\lambda_{0,i}}}\big\vert \bfX_n \Big]=   O\Big(\dfrac{1}{\lambda_{0,i}}\sqrt{\dfrac{p}{n}}\Big)+ O(n^{-1/2+\delta}),
\end{align*}
where $\lambda_{0,k}\preccurlyeq p/n^{1/4}$.

On the other hand, if $\lambda_{0,k} \succ p / n^{1/4}$, then it follows that
$$\max_{i\leq k}a_i-\min_{i\leq k}a_i \preccurlyeq \dfrac{p}{\lambda_{0,k}}\prec n^{1/2}.$$

In this case, Theorem \ref{main-cor:conv_rate_eigenvalue_equal_a} applies and yields
\begin{align*}
    \bbE\Big[\dfrac{\lambda_{(i)}-\lambda_{0,i}}{{\lambda_{0,i}}}\big\vert \bfX_n \Big]=   O\Big(\dfrac{p}{n\lambda_{0,i}}\Big)+ O(n^{-1/2+\delta}).
\end{align*}

Since $\lambda_{0,i} \succ p / n^{1/4}$ in this regime, it follows that
$$O\left( \frac{p}{n\lambda_{0,i}} \right) + O\left( n^{-1/2 + \delta} \right) 
= O\left( \frac{1}{\lambda_{0,i}} \sqrt{\frac{p}{n}} \right) + O\left( n^{-1/2 + \delta} \right) 
= O(n^{-1/2 + \delta}),$$
holds. 

Therefore, in both regimes, we obtain the unified posterior bound:
\begin{align*}
    \bbE\Big[\dfrac{\lambda_{(i)}-\lambda_{0,i}}{{\lambda_{0,i}}}\big\vert \bfX_n \Big]=   O\Big(\dfrac{1}{\lambda_{0,i}}\sqrt{\dfrac{p}{n}}\Big)+ O(n^{-1/2+\delta}).
\end{align*}

\end{proof}


\subsubsection*{Proof of Corollary \ref{main-cor:conv_rate_eigenvector}}
\begin{proof}
Consider the spectral decomposition of $nS = QWQ^T$. Then, on $D_\epsilon$, we have the following inequality:
\begin{align*}
    (\xi_j^T\hat{\xi}_j)^2 &= ([Q\Gamma]_{\cdot j}^T Q_{\cdot j})^2\\
    &= [(Q\Gamma e_j^T)^T (Qe_j^T)]^2\\
    &= \Gamma_{jj}^2\\
    &\geq 1-\epsilon^2,
\end{align*}
for $j=1,\ldots,k$.

By Theorem 3.2 of \cite{wang2017asymptotics}, the following limit holds:
$$\abs{\xi_{0,j}^T\hat{\xi}_j}-(1+\bar{d}d_j)^{-\frac{1}{2}}= O_p(\zeta),$$
where $\zeta_j = \dfrac{1}{\lambda_{0,j}}\sqrt{\dfrac{p}{n}}+ \dfrac{p}{n^{3/2}\lambda_{0,j}}+\dfrac{1}{n}$.

Applying the triangular inequality of angles from \cite{castano2016angles}, where
$$\theta(u,v) = \arccos\dfrac{\abs{u^Tv}}{\abs{\abs{u}}\abs{\abs{v}}},$$
we obtain:
\begin{align*}
    \abs{\xi_{0,j}^T\xi_j}&\leq \cos\Big[ \arccos(\abs{\hat{\xi}_j^T\xi_j})+\arccos(\abs{\xi_{0,j}^T\hat{\xi}_j})\Big]\\
    &\leq\cos\Big[\arccos(\sqrt{1-\epsilon^2})+ \arccos(\dfrac{1}{\sqrt{1+\bar{d}d_j}}+O_p(\zeta_j)) \Big]\\
    &=\cos\Big[\arccos \Big(\sqrt{1-\epsilon^2}\sqrt{\dfrac{1}{1+\bar{d}d_j}+O_p(\zeta_j)}-\epsilon\sqrt{1-\big(\dfrac{1}{1+\bar{d}d_j}+O_p(\zeta_j)\big)^2 }\Big)\Big]\\
    &=\cos\Big[\arccos \big(\sqrt{\dfrac{1-\epsilon^2}{1+\bar{d}d_j}}-\sqrt{\dfrac{\bar{d}d_j}{1+\bar{d}d_j}\epsilon^2}  + O_p(\zeta_j)\big)\Big]\\
    &\leq (1+\bar{d}d_j)^{-1/2}+ O_p(\zeta_j).
\end{align*}

Similarly, applying the lower bound using $\theta(\xi_{0,j},\xi_j)\geq \theta(\xi_{0,j},\hat{\xi}_j) - \theta(\hat{\xi}_j,\xi_j)$.
\begin{align*}
    \abs{\xi_{0,j}^T\xi_j}&\geq \cos\Big[ \arccos(\abs{\xi_{0,j}^T\hat{\xi}_j})-\arccos(\abs{\hat{\xi}_j^T\xi_j})\Big]\\
    &\geq \cos\Big[ \arccos(\dfrac{1}{\sqrt{1+\bar{d}d_j}})-\arccos(\sqrt{1-\epsilon^2}) \Big]\\
    &=\cos\Big[\arccos \Big(\sqrt{1-\epsilon^2}\sqrt{\dfrac{1}{1+\bar{d}d_j}+O_p(\zeta_j)}+\epsilon\sqrt{1-\big(\dfrac{1}{1+\bar{d}d_j}+O_p(\zeta_j)\big)^2 }\Big)\Big]\\
    & = \cos\Big[\arccos\big(\sqrt{\dfrac{1-\epsilon^2}{1+\bar{d}d_j}}+\sqrt{\dfrac{\bar{d}d_j}{1+\bar{d}d_j}\epsilon^2} +O_p(\zeta_j)\big) \Big]\\
    & = \sqrt{\dfrac{1-\epsilon^2}{1+\bar{d}d_j}}+\epsilon\sqrt{\dfrac{\bar{d}d_j}{1+\bar{d}d_j}} + O_p(\zeta_j).
\end{align*}

Therefore, we obtain the expectation:
\begin{align*}
    \bbE\Big[1-(\xi_{0,j}^T\xi_j)^2 \big\vert \bfX_n \Big]
    &= \sup_{D_\epsilon}(1-\abs{\xi_{0,j}^T\xi_j}^2)\\
    &= 1- \Big(\dfrac{1-\epsilon^2}{1+\bar{d}d_j}+O_p(\zeta_j)\Big)- \Big(\dfrac{\bar{d}d_j\epsilon^2}{1+\bar{d}d_j}\Big) - 2\sqrt{\Big(\dfrac{1-\epsilon^2}{1+\bar{d}d_j}\Big) \Big(\dfrac{\bar{d}d_j\epsilon^2}{1+\bar{d}d_j}\Big)}+ O_p(\zeta_j)\\
    &= \dfrac{\bar{d}d_j}{1+\bar{d}d_j} + O(\dfrac{p}{n\lambda_{0,j}}\epsilon^2) + O\Big(\epsilon \sqrt{\dfrac{p}{n\lambda_{0,j}}}\Big)+ O_p(\zeta_j)\\
    &= \dfrac{p\bar{d}}{n\lambda_{0,j}+p\bar{d}} + O(\dfrac{p}{n\lambda_{0,j}}\epsilon^2) + O\Big(\epsilon \sqrt{\dfrac{p}{n\lambda_{0,j}}}\Big)+ O_p(\zeta_j).
\end{align*}

If $\lambda_{0,1} \preccurlyeq p / n^{1/4}$, then by proof of Lemma \ref{lem:post_exp_diff_ordered}, we obtain
\begin{align*}
    \abs{\bbE\Big[1-(\xi_{0,j}^T\xi_{j})^2 \vert \bfX_n \Big] - \bbE\Big[1-(\xi_{0,j}^T\xi_{(j)})^2 \vert \bfX_n \Big]} &= \abs{\bbE\Big[\big((\xi_{0,j}^T\xi_{(j)})^2-(\xi_{0,j}^T\xi_{j})^2) \big)\cdot I(\lambda_{(i)})\neq \lambda_i)\vert \bfX_n \Big]}\\
    &\leq 2\bbE\Big[I(\lambda_{(i)})\neq \lambda_i)\vert \bfX_n \Big]\\
    & = 2\bbE\Big[\bbE[I(\lambda_{(i)})\neq \lambda_i)\vert \Gamma]\cdot (I(\Gamma\in D_\epsilon)+I(\Gamma\notin D_\epsilon))\vert \bfX_n \Big]\\
    &\leq 2\bbE\Big[\bbE[I(\lambda_{(i)})\neq \lambda_i)\vert \Gamma]\cdot I(\Gamma\in D_\epsilon)\vert \bfX_n \Big]+2\bbE\Big[\bbE[1\vert \Gamma]\cdot I(\Gamma\notin D_\epsilon)\vert \bfX_n \Big]\\
    &\leq 2(1-\prod_{i=1}^p (1-p_i)) + \pi(\Gamma\notin D_\epsilon\vert \bfX_n)\\
    & \leq  2\sum_{i=1}^p p_i + O\Big( \exp(-\frac{\kappa}{2} \epsilon) \Big) + O\Big( \Big( \frac{n\lambda_{0,k} + p}{p} \Big)^{-\epsilon^2 a_n} \Big).
\end{align*}

Similarly, if $\lambda_{0,1} \succ p / n^{1/4}$, then by proof of Lemma \ref{lem:post_exp_equal_ordered}, we obtain
\begin{align*}
    \abs{\bbE\Big[1-(\xi_{0,j}^T\xi_{j})^2 \vert \bfX_n \Big] - \bbE\Big[1-(\xi_{0,j}^T\xi_{(j)})^2 \vert \bfX_n \Big]}\leq 2\sum_{i=1}^p p_i  + O\Big( \exp(-\frac{\tau}{2}n \epsilon^2)\cdot\Big(\dfrac{2\lambda_{0,1}}{\lambda_{0,k}}\Big)^{kq} \Big) + O\Big( \Big( \frac{n\lambda_{0,k} + p}{p} \Big)^{-\epsilon^2 a_n} \Big).
\end{align*}

By the assumptions of the theorem, in particular $\lambda_{0,1} / \lambda_{0,k} = O(1)$ and $\epsilon \asymp n^{-1/4}$, we have $\sum_{i=1}^p p_i = O(1/n)$. Therefore,
$$\abs{\bbE\Big[1-(\xi_{0,j}^T\xi_{j})^2 \vert \bfX_n \Big] - \bbE\Big[1-(\xi_{0,j}^T\xi_{(j)})^2 \vert \bfX_n \Big]} \preccurlyeq n^{-1/2+\delta} + \dfrac{1}{\lambda_{0,i}}\sqrt{\dfrac{p}{n}}.$$

Hence, we obtain the final posterior bound
$$\bbE\Big[1-\abs{\xi_{0,j}^T\xi_{(j)}}^2 \big\vert \bfX_n \Big] = O (\dfrac{p}{n\lambda_{0,j}})+ O_p(\zeta_j),$$
for $j=1,\ldots,k$, where $\zeta_j = \dfrac{1}{\lambda_{0,j}}\sqrt{\dfrac{p}{n}}+ \dfrac{p}{n^{3/2}\lambda_{0,j}}+\dfrac{1}{n}$.

\end{proof}

We consider the spectral decomposition $\Sigma = \Gamma \Lambda \Gamma^T$, where
$$\Lambda = \begin{pmatrix}
    \Lambda_1 & 0 \\ 0 & \Lambda_2
\end{pmatrix},$$
with $\Lambda_1$ being a $k\times k$ matrix and $\Lambda_2$ an $(p-k)\times (p-k)$ matrix. Furthermore, we write $\Gamma = (\Gamma_1 ,\Gamma_2)$ where $\Gamma_1$ is a $p\times k$ matrix and $\Gamma_2$ is a $p\times(p-k)$ matrix. 

\begin{lemma}[Spectral difference of eigenvectors]\label{lem:spec_diff_eigenvec}
    For the spiked eigenvector matrix $\Gamma_1$, the following inequality holds on the set $D_\epsilon$:
    $$\abs{\abs{\Gamma_{0,1}-\Gamma_1}}^2 = O(d_k)+O_p(\zeta_k),$$
    where $\zeta_k = \dfrac{1}{\lambda_{0,k}}\sqrt{\dfrac{p}{n}}+ \dfrac{p}{n^{3/2}\lambda_{0,k}}+\dfrac{1}{n}$.
\end{lemma}
\begin{proof}

\begin{align*}
    \abs{\abs{\Gamma_{0,1}-\Gamma_1}}^2 &\leq \abs{\abs{\Gamma_{0,1}-\Gamma_1}}_F^2\\
    &\leq \sum\limits_{i=1}^k \abs{\abs{\xi_{0,i}-\xi_i}}_2^2\\
    &= 2\sum\limits_{i=1}^k (1-\abs{\xi_{0,i}^T\xi_i})   \\
    &\leq 2\sum\limits_{i=1}^k \bigg(1-\big(\sqrt{\dfrac{1-\epsilon^2}{1+d_i}}+\epsilon\sqrt{\dfrac{d_i}{1+d_i}} + O_p(\zeta_i) \big)\bigg),
\end{align*}
where $\zeta_i = \dfrac{1}{\lambda_{0,i}}\sqrt{\dfrac{p}{n}}+ \dfrac{p}{n^{3/2}\lambda_{0,i}}+\dfrac{1}{n}$. The last inequality follows from proof of Theorem \ref{main-cor:conv_rate_eigenvector}.

\begin{align*}
    2\sum\limits_{i=1}^k \bigg(1-\big(\sqrt{\dfrac{1-\epsilon^2}{1+d_i}}+\epsilon\sqrt{\dfrac{d_i}{1+d_i}} \big)\bigg)
    &\leq 2\sum\limits_{i=1}^k \bigg(1-\sqrt{\dfrac{1-\epsilon^2}{1+d_i}}\bigg)\\
    &\leq 2k \bigg(1-\sqrt{\dfrac{1-\epsilon^2}{1+d_k}}\bigg)\\
    &\leq 2k \bigg(\sqrt{1+d_k}-\sqrt{1-\epsilon^2}\bigg)\\
    &= 2k \dfrac{d_k+\epsilon^2}{\sqrt{1+d_k}+\sqrt{1-\epsilon^2}}\\
    &= O(d_k).
\end{align*}

Since $\sum_{j=1}^k \zeta_j \leq k \zeta_k$, it follows that $\abs{\abs{\Gamma_{0,1}-\Gamma_1}}^2 = O(d_k)+O_p(\zeta_k).$

\end{proof}

\subsubsection*{Proof of Corollary \ref{main-cor:conv_rate_cov}}

\begin{proof}

\begin{align}
    \abs{\abs{\Sigma-\Sigma_0}}_F^2 &\leq 2\sum\limits_{i=1}^2 \abs{\abs{\Gamma_i\Lambda_i\Gamma_i^T - \Gamma_{0,i}\Lambda_{0,i}\Gamma_{0,i}^T}}_F^2\nonumber\\
    &\leq 2\sum\limits_{i=1}^2 \big[\abs{\abs{(\Gamma_{i}-\Gamma_{0,i})\Lambda_{0,i}\Gamma_{i}^T}}_F +\abs{\abs{\Gamma_{i}(\Lambda_{i}-\Lambda_{0,i})\Gamma_{i}^T}}_F +
    \abs{\abs{\Gamma_{0,i}\Lambda_{0,i}(\Gamma_{i}-\Gamma_{0,i})^T}}_F \big]^2\nonumber\\
    &\leq 6\sum\limits_{i=1}^2 \big[\abs{\abs{(\Gamma_{i}-\Gamma_{0,i})\Lambda_{0,i}\Gamma_{i}^T}}_F^2 +\abs{\abs{\Gamma_{i}(\Lambda_{i}-\Lambda_{0,i})\Gamma_{i}^T}}_F^2 +
    \abs{\abs{\Gamma_{0,i}\Lambda_{0,i}(\Gamma_{i}-\Gamma_{0,i})^T}}_F^2\big]\label{ineq:third}\\
    &\leq 6\sum\limits_{i=1}^2 \big[2\abs{\abs{\Gamma_{0,i}-
    \Gamma_i}}^2\abs{\abs{\Lambda_{0,i}}}_F^2  + \abs{\abs{\Lambda_{0,i}-\Lambda_i}}_F^2  \big]\nonumber\\
    &\leq 6\Big[8(\sum\limits_{i=1}^k\lambda_{0,i}^2)\cdot \big(O(d_k) + O_p(\zeta_k)\big) + \sum\limits_{i=1}^k (\lambda_{0,i}-\lambda_i)^2 \Big] + 6\Big[8 \sum\limits_{i>k} \lambda_{0,i}^2 + \sum\limits_{i>k}(\lambda_{0,i}-\lambda_i)^2 \Big],\label{ineq:fifth}\\
    &\leq 6\Big[8\sum\limits_{i=1}^k\lambda_{0,i}^2 \cdot O(d_k) + \sum\limits_{i=1}^k (\lambda_{0,i}-\lambda_i)^2 \Big] + 6\Big[8\sum\limits_{i>k} \lambda_{0,i}^2 + \sum\limits_{i>k}(\lambda_{0,i}-\lambda_i)^2 \Big] + O_p(\sum\limits_{i=1}^k\lambda_{0,i}^2\zeta_k)\nonumber
\end{align}

The inequality \eqref{ineq:third} follows from the Cauchy-Schwarz inequality, and \eqref{ineq:fifth} follows from Lemma \ref{lem:spec_diff_eigenvec}.

The posterior expectation is given by:
\begin{align*}
    \bbE\Big[(\lambda_i-\lambda_{0,i})^2\big\vert \bfX_n \Big]  &=\dfrac{\displaystyle\int(\lambda_i-\lambda_{0,i})^2 \prod \lambda_i^{-a_i-\frac{n}{2}}\exp(-\dfrac{c_i}{2\lambda_i})d\Lambda d\Gamma}{\displaystyle\int\prod \lambda_i^{-a_i-\frac{n}{2}}\exp(-\dfrac{c_i}{2\lambda_i})d\Lambda d\Gamma}\nonumber \\
    &=\dfrac{\displaystyle\int f(c)\Big(\prod c_{i}\Big)^{-a_i-\frac{n}{2}+1} d\Gamma }{\displaystyle\int \Big(\prod c_{i}\Big)^{-a_i-\frac{n}{2}+1}d\Gamma}\nonumber\\
    &\leq \sup_{D_\epsilon}f(c) + \dfrac{\displaystyle\int_{D_\epsilon^c} f(c)\Big(\prod c_{i}\Big)^{-a_i-\frac{n}{2}+1} d\Gamma }{\displaystyle\int \Big(\prod c_{i}\Big)^{-a_i-\frac{n}{2}+1}d\Gamma},  
\end{align*}
where $$f(c)=\bbE_{\lambda_i}[(\lambda_i-\lambda_{0,i})^2]= \Big(\lambda_{0,i}-\dfrac{c_i}{n+2a_i-4}\Big)^2 + \dfrac{2c_i^2}{(n+2a_i-4)^2(n+2a_i-6)^2},$$
with $\lambda_i\overset{iid}{\sim}InvGam(a_i+\dfrac{n}{2}-1,\dfrac{c_{i}}{2})$.

Therefore, we obtain for $i=1,\ldots,k$:
\begin{align*}
    \sup_{D_\epsilon} f(c) & = \Big(O(\sqrt{\dfrac{p}{n}}+ \lambda_{0,i} n^{-\frac{1}{2}+\delta})\Big)^ 2 +  O(\dfrac{\lambda_{0,i}^2}{n^3})\\
    &= O\Big(\dfrac{p}{n}+ \dfrac{\lambda_{0,i}^2}{n^{1-2\delta}}\Big),
\end{align*}
for all small $\delta>0$.

The posterior expectation is given by:
$$\bbE\Big[(\lambda_i-\lambda_{0,i})^2\big\vert \bfX_n \Big]= O\Big(\dfrac{p}{n}+ \dfrac{\lambda_{0,i}^2}{n^{1-2\delta}}\Big), $$
for $i=1,\ldots,k$.

Similarly, we obtain the upper bound for $f(c)$ over $D_\epsilon$. Using the assumption that $a_i\succ p$ for $i>k$, we obtain:
\begin{align*}
    \sup_{D_\epsilon} f(c) & = \Big(\lambda_{0,i}-\dfrac{c_i}{n+2a_i-4}\Big)^2 + \dfrac{2}{(n+2a_i-6)^2} \cdot \Big(\dfrac{c_i}{(n+2a_i-4)}\Big)^2\\
    &\leq \lambda_{0,i}^2 + \dfrac{1}{a_i^2} \cdot \lambda_{0,i}^2\\
    &\leq 2\lambda_{0,i}^2.
\end{align*}

The posterior expectation is given by:
$$\bbE\Big[(\lambda_i-\lambda_{0,i})^2\big\vert \bfX_n \Big]= O(\lambda_{0,i}^2), $$
for all $i>k$.

Therefore, we obtain the following inequality of posterior expectation:
\begin{align*}
    \bbE\Big[\abs{\abs{\Sigma-\Sigma_0}}_F^2\vert \bfX_n\Big] &= O(d_k\sum\limits_{i=1}^k\lambda_{0,i}^2+ \dfrac{p}{n}+ n^{-1+2\delta}\sum\limits_{i=1}^k\lambda_{0,i}^2) + O(\sum\limits_{i>k}\lambda_{0,i}^2)+O_p(\sum\limits_{i=1}^k\lambda_{0,i}^2\zeta_k)\\
    &= O\Big(\big(n^{-1+2\delta}+ \dfrac{p}{n\lambda_{0,k}}\big)\sum\limits_{i=1}^k\lambda_{0,i}^2 + \dfrac{p}{n}+  \sum\limits_{i=k+1}^p\lambda_{0,i}^2\Big)+O_p(\sum\limits_{i=1}^k\lambda_{0,i}^2\zeta_k)\\
    &= O\Big(\big(n^{-1+2\delta}+ \dfrac{p}{n\lambda_{0,k}}\big)\sum\limits_{i=1}^k\lambda_{0,i}^2  + p\Big) +O_p(\sum\limits_{i=1}^k\lambda_{0,i}^2\zeta_k).
\end{align*}

Since $\lambda_{0,1} / \lambda_{0,k}$ is bounded by a positive constant, the posterior expectation of the covariance is given by:
\begin{align*}
    \bbE\Big[\dfrac{\abs{\abs{\Sigma-\Sigma_0}}_F^2}{\abs{\abs{\Sigma_0}}_F^2}\big\vert \bfX_n \Big] &= O\Big(\dfrac{p }{\sum\limits_{i=1}^p\lambda_{0,i}^2} + n^{-1+2\delta}\dfrac{\sum\limits_{i=1}^k \lambda_{0,i}^2}{\sum\limits_{i=1}^p\lambda_{0,i}^2}\Big)+ O \Big(\dfrac{p}{n\lambda_{0,k}}\dfrac{\sum\limits_{i=1}^k\lambda_{0,i}^2}{\sum\limits_{i=1}^p\lambda_{0,i}^2}\Big)+O_p(\dfrac{\sum\limits_{i=1}^k \lambda_{0,i}^2}{\sum\limits_{i=1}^p\lambda_{0,i}^2}\zeta_k)\\
    &= O\Big(\dfrac{p }{\lambda_{0,1}^2 + p} + n^{-1+2\delta}\dfrac{\lambda_{0,1}^2}{\lambda_{0,1}^2 + p}\Big)+ O \Big(\dfrac{p}{n\lambda_{0,k}}\dfrac{\lambda_{0,1}^2}{\lambda_{0,1}^2 + p}\Big)+O_p(\dfrac{\lambda_{0,1}^2}{\lambda_{0,1}^2 + p}\zeta_k).
\end{align*}

\hfill\break

We consider the spectral decomposition $S = \hat{\Gamma} \hat{\Lambda} \hat{\Gamma}^T$, where
$$\hat{\Lambda}= \begin{pmatrix}
    \hat{\Lambda}_1 & 0 \\ 0 & \hat{\Lambda}_2
\end{pmatrix},$$
with $\hat{\Lambda}_1$ being a $k\times k$ matrix and $\hat{\Lambda}_2$ an $(p-k)\times (p-k)$ matrix. Furthermore, we write $\hat{\Gamma} = (\hat{\Gamma}_1 ,\hat{\Gamma}_2)$ where $\hat{\Gamma}_1$ is a $p\times k$ matrix and $\hat{\Gamma}_2$ is a $p\times(p-k)$ matrix.

By Theorem 3.2 of \cite{wang2017asymptotics}, the following limit holds for each $j=1,\ldots,k$
$$\abs{\xi_{0,j}^T\hat{\xi}_j} = (1+\bar{d}d_j)^{-\frac{1}{2}}+ O_p(\zeta_j).$$

Using this, we bound the spectral norm between the true and sample eigenvectors as
\begin{align*}
    \abs{\abs{\Gamma_{0,1}-\hat{\Gamma}_1}}^2 &\leq \abs{\abs{\Gamma_{0,1}-\hat{\Gamma}_1}}_F^2\\
    &\leq \sum_{j=1}^k \abs{\abs{\xi_{0,j}-\hat{\xi}_j}}_2^2\\
    &= 2\sum_{j=1}^k (1- \abs{\xi_{0,j}^T\hat{\xi}_j})\\
    &= 2\sum_{j=1}^k \Big(1- ((1+\bar{d}d_j)^{-\frac{1}{2}}+ O_p(\zeta_j))\Big)\\
    &=  2\sum_{j=1}^k \dfrac{\bar{d}d_j}{\sqrt{1+\bar{d}d_j}(\sqrt{1+\bar{d}d_j}+1)} + O_p(\zeta_k)\\
    &=  2\sum_{j=1}^k O(d_j) + O_p(\zeta_k)\\
    &= O(d_k) + O_p(\zeta_k).
\end{align*}

By the same procedure, we have
\begin{align}
    \abs{\abs{\Sigma-\Sigma_0}}_F^2 
    &\leq 6\sum\limits_{i=1}^2 \big[2\abs{\abs{\Gamma_{0,i}-
    \hat{\Gamma}_i}}^2\abs{\abs{\Lambda_{0,i}}}_F^2  + \abs{\abs{\Lambda_{0,i}-\hat{\Lambda}_i}}_F^2  \big]\nonumber\\
    &\leq 6\Big[8(\sum\limits_{i=1}^k\lambda_{0,i}^2)\cdot \big(O(d_k) + O_p(\zeta_k)\big) + \sum\limits_{i=1}^k (\lambda_{0,i}-\hat{\lambda}_i)^2 \Big] + 6\Big[8 \sum\limits_{i>k} \lambda_{0,i}^2 + \sum\limits_{i>k}(\lambda_{0,i}-\hat{\lambda}_i)^2 \Big].\label{ineq:sample_cov}
\end{align}

By Lemma \ref{lem:eig}, we obtain 
$$(\lambda_{0,i}-\hat{\lambda}_i)^2 = \begin{cases}
    O(\dfrac{p^2}{n^2}+\dfrac{\lambda_{0,i}^2}{n^{1-2\delta}})\quad &\text{for }i=1,\ldots,k\\
     O(\dfrac{p^2}{n^2})\quad &\text{for }i=k+1,\ldots,n\\
     O(1)\quad &\text{for }i=n+1,\ldots,p.
\end{cases}$$

Substituting into \eqref{ineq:sample_cov}, we obtain
\begin{align*}
    \abs{\abs{\Gamma_{0,1}-\hat{\Gamma}_1}}^2 &\leq 6\Big[8(\sum\limits_{i=1}^k\lambda_{0,i}^2)\cdot \big(O(d_k) + O_p(\zeta_k)\big) + \sum\limits_{i=1}^k O(\dfrac{p^2}{n^2}+\dfrac{\lambda_{0,i}^2}{n^{1-2\delta}}) \Big] + 6\Big[8 \sum\limits_{i>k} O(1) + \sum\limits_{i=k+1}^n O(\dfrac{p^2}{n^2})  + \sum\limits_{i=n+1}^p O(1) \Big]\\
    &= O\Big(d_k \sum_{i=1}^k\lambda_{0,i}^2 + \dfrac{p^2}{n^2}+n^{1-2\delta}\sum_{i=1}^k\lambda_{0,i}^2 +p + \dfrac{p^2}{n}\Big) + O_p\Big(\sum_{i=1}^k\lambda_{0,i}^2\zeta_k\Big)\\
    &= O\Big(\dfrac{p}{n\lambda_{0,k}} \sum_{i=1}^k\lambda_{0,i}^2 +\dfrac{p^2}{n}+n^{1-2\delta}\sum_{i=1}^k\lambda_{0,i}^2\Big) + O_p(\sum_{i=1}^k\lambda_{0,i}^2\zeta_k) \lambda_{0,i}^2 +\dfrac{p^2}{n}\Big) + O_p\Big(\sum_{i=1}^k\lambda_{0,i}^2\zeta_k\Big).
\end{align*}

Since $\lambda_{0,1} / \lambda_{0,k}$ is bounded by a positive constant, it follows that:
\begin{align*}
\dfrac{\abs{\abs{S-\Sigma_0}}_F^2}{\abs{\abs{\Sigma_0}}_F^2} &= O\Big(\dfrac{p}{n\lambda_{0,k}} \dfrac{\sum_{i=1}^k\lambda_{0,i}^2}{\sum_{i=1}^p\lambda_{0,i}^2} +\dfrac{p^2}{n\sum_{i=1}^p\lambda_{0,i}^2}+n^{-1+2\delta}\dfrac{\sum\limits_{i=1}^k \lambda_{0,i}^2}{\sum\limits_{i=1}^p\lambda_{0,i}^2}\Big) + O_p\Big(\dfrac{\sum_{i=1}^k\lambda_{0,i}^2}{\sum_{i=1}^p\lambda_{0,i}^2}\zeta_k\Big) \\
&= O\Big(\dfrac{p^2}{n(\lambda_{0,1}^2 + p)}+n^{-1+2\delta}\dfrac{ \lambda_{0,1}^2}{\lambda_{0,1}^2+p}\Big) + O\Big(\dfrac{p}{n\lambda_{0,k}} \dfrac{\lambda_{0,1}^2}{\lambda_{0,1}^2 + p} \Big) + O_p\Big(\dfrac{\lambda_{0,1}^2}{\lambda_{0,1}^2 + p}\zeta_k\Big).
\end{align*}

\end{proof}

\bibliographystyle{dcu}
\bibliography{SIW}